\documentclass[preprint,11pt]{elsarticle}

\usepackage{graphicx}
\usepackage{amssymb}
\usepackage{amsmath}
\usepackage{amscd}
\newtheorem{theo}{Theorem}
\newtheorem{theorem}[theo]{Theorem}
\newdefinition{definition}[theo]{Definition}

\newtheorem{proposition}[theo]{Proposition}
\newdefinition{remark}[theo]{Remark}
\newproof{proof}{Proof}

\makeatletter \@addtoreset{equation}{section}
\@addtoreset{theo}{section} \makeatother

\begin{document}

\begin{frontmatter}


\title{Regular Reduction of Controlled Hamiltonian System with
Symplectic Structure and Symmetry }

\author{Jerrold E. Marsden}
\address{Control and Dynamical Systems,
California Institute of Technology, Pasadena, CA 91125 USA}
\author{Hong Wang\corref{cor1}}
\ead{hongwang@nankai.edu.cn}
\author{Zhenxing Zhang}
\address{School of Mathematical Sciences and LPMC, Nankai University, Tianjin 300071, P.R.China}
\cortext[cor1]{Corresponding author. Tel.: 0086-022-23501233.
Address: School of Mathematical Sciences, Nankai University, Tianjin
300071, P.R.China.}

\markboth{J.E.Marsden, H.Wang and Z.X.Zhang }{Regular Reduction of
Controlled Hamiltonian System }

\begin{abstract}
In this paper, our goal is to study the regular reduction theory of
a regular controlled Hamiltonian (RCH) system with symplectic
structure and symmetry, and this reduction is an extension of
the regular symplectic reduction theory of a Hamiltonian system under
regular controlled Hamiltonian equivalence conditions. Thus, in
order to describe uniformly RCH systems defined on a cotangent
bundle and on the regular reduced spaces, we first define a kind of
RCH systems on a symplectic fiber bundle. Then we introduce regular
point and regular orbit reducible RCH systems with symmetries by using
momentum maps and the associated reduced symplectic forms. Moreover,
we give regular point and regular orbit reduction theorems for RCH
systems to explain the relationships between RpCH-equivalence,
RoCH-equivalence for the reducible RCH systems with symmetries and
RCH-equivalence for the associated reduced RCH systems. Finally, as the
applications we regard rigid body and heavy top as well as them with
internal rotors as the regular point reducible RCH systems on the
rotation group $\textmd{SO}(3)$ and on the Euclidean group
$\textmd{SE}(3)$, as well as on their generalizations, respectively,
and discuss their RCH-equivalences. We also describe the RCH system
and RCH-equivalence from the viewpoint of port Hamiltonian system
with a symplectic structure.
\end{abstract}


\end{frontmatter}

\tableofcontents

\section{Introduction}

Symmetry is a general phenomenon in the natural world, but it is
widely used in the study of mathematics and mechanics. The reduction
theory for mechanical system with symmetry has its origin in the
classical work of Euler, Lagrange, Hamilton, Jacobi, Routh,
Liouville and Poincar\'{e} and its modern geometric formulation in
the general context of symplectic manifolds and equivariant momentum
maps is developed by Meyer, Marsden and Weinstein; see Abraham and
Marsden \cite{abma78} or Marsden and Weinstein \cite{mawe74} and
Meyer \cite{me73}. The main goal of reduction theory in mechanics is
to use conservation laws and the associated symmetries to reduce the
number of dimensions of a mechanical system required to be
described. So, such reduction theory is regarded as a useful tool
for simplifying and studying concrete mechanical systems. Reduction
is a very general procedure that is applied to arbitrary dynamical
systems with symmetries. However, it is particularly powerful for
conservative systems whose symmetries are induced by a momentum map;
see Abraham and Marsden \cite{abma78}, Arnold \cite{ar89}, Marsden
\cite{ma92}, Marsden et al. \cite{mamiorpera07}, Marsden and Ratiu
\cite{mara99} and Ortega and Ratiu \cite{orra04} for more details.\\

It is well-known that Hamiltonian reduction theory is one of the
most active subjects in the study of modern analytical mechanics and
applied mathematics, in which a lot of deep and beautiful results
have been obtained, see the studies given by Abraham and Marsden
\cite{abma78}, Arnold \cite{ar89}, Leonard and Marsden
\cite{lema97}, Marsden et al.
\cite{ma92,mamiorpera07,mara99,mawe74}, Ortega and Ratiu
\cite{orra04} etc. on regular point reduction and regular orbit
reduction, singular point reduction and singular orbit reduction,
optimal reduction and reduction by stages for Hamiltonian systems
and so on; and there is still much to be done in this
subject.\\

On the other hand, just as we have known that the theory of
mechanical control systems presents a challenging and promising
research area between the study of classical mechanics and modern
nonlinear geometric control theory and there have been a lot of
interesting results. Such as Bloch et al. in
\cite{blchlema01,blkrmaal92,blle02,bllema00}, referred to the use of
feedback control to realize a modification to the structure of a
given mechanical system; Blankenstein et al. in \cite{blorvds02},
Crouch and Van der Schaft in \cite{crvds87}, Nijmeijer and Van der
Schaft in \cite{nivds90}, Van der Schaft in
\cite{vds82,vds86,vdsma95,vds00,vds06}, referred to the reduction
and control of implicit (port) Hamiltonian systems, and to the use
of feedback control to stabilize mechanical systems; and Chang et
al. in \cite{chbllemawo02, chma04}, studied the controlled
Hamiltonian (CH) system by using almost Poisson tensor.
However, we found that the authors in \cite{chbllemawo02, chma04},
didn't consider the phase spaces of CH system and
the reduced CH system, the change of geometrical
structures of the phase spaces of the CH systems,
as well as the momentum map of the CH system with
symmetry, and hence cannot determine precisely the geometrical
structures of phase spaces of the reduced CH systems.
Moreover, it is impossible to give precisely the relations of the
reduced controlled Hamiltonian equivalences, if don't consider
the different Lie group actions and momentum maps. Thus,
we think that there are a lot of serious wrong of rigor
for the definitions of CH system and the reduced CH system,
as well as CH-equivalence and the reduced CH-equivalence
in Chang et al. \cite{chbllemawo02, chma04},
and we want to correct their work.\\

In this paper, our goal is to study regular reduction theory of
a regular controlled Hamiltonian (RCH) system with symplectic
structure and symmetry, by combining with momentum map and the regular
symplectic reduction theory of a Hamiltonian system. In order to do
this, our idea in this paper is that we first define a CH system on
$T^*Q$ by using a symplectic form, and such system is called a RCH
system, and then regard the associated Hamiltonian system on $T^*Q$
as a spacial case of the RCH system without external force and
control. Thus, the set of Hamiltonian systems on $T^*Q$ is a subset
of the set of RCH systems on $T^*Q$. We hope to study regular
reduction theory of the RCH system with symplectic structure and
symmetry, as an extension of the regular symplectic reduction theory of
a Hamiltonian system under regular controlled Hamiltonian equivalence
conditions. The main contributions in this paper is given as
follows. (1) In order to describe uniformly RCH systems defined on a
cotangent bundle and on the regular reduced spaces, we define a kind
of RCH systems on a symplectic fiber bundle by using its symplectic
form; (2) We give the regular point and the regular orbit reducible RCH
systems by using momentum maps and the associated reduced symplectic
forms, and prove regular point and regular orbit reduction theorems
for the RCH systems (see Theorems 4.4 and 5.4); (3) We prove that rigid
body with external force torque, rigid body with internal rotors and
heavy top with internal rotors are all RCH systems, and as a pair of
the regular point reduced RCH systems, rigid body with internal rotors
(or external force torque) and heavy top with internal rotors are
RCH-equivalent; (4) We describe the RCH system from the viewpoint of
port Hamiltonian system with a symplectic structure, and state the
relationship between RCH-equivalence of RCH systems and equivalence
of port Hamiltonian systems.\\

A brief of outline of this paper is as follows. In the second
section, we review some relevant definitions and basic facts about
momentum map, symplectic fiber bundle, Lie group lifted actions on
(co-)tangent bundles and reduction, which will be used in subsequent
sections. The RCH systems are defined by using the symplectic forms
on a symplectic fiber bundle and on the cotangent bundle of a
configuration manifold, respectively, and RCH-equivalence is
introduced in the third section. From the fourth section we begin to
discuss the RCH system with symmetry by combining with the regular
symplectic reduction theory. The regular point and regular orbit
reducible RCH systems are considered respectively in the fourth
section and the fifth section, and give the regular point and
regular orbit reduction theorems for the RCH systems to explain the
relationships between the RpCH-equivalence, RoCH-equivalence for
reducible RCH systems with symmetries and the RCH-equivalence for
the associated reduced RCH systems. As the applications of the
theoretical results, in sixth section, we first give the regular
point reduced RCH systems on a Lie group $G$ and on its
generalization $G\times V$, which are the RCH systems on a coadjoint
orbit $\mathcal{O}_\mu$ of $G$ and on its generalization
$\mathcal{O}_\mu \times V\times V^\ast$. Then we describe uniformly
the rigid body and heavy top as well as them with internal rotors as
the regular point reducible RCH systems on the rotation group
$\textmd{SO}(3)$ and on the Euclidean group $\textmd{SE}(3)$, as
well as on their generalizations, respectively, and give their
regular point reduced RCH systems and discuss their RCH-equivalences.
In order to understand well the abstract definition of RCH system,
we also describe the RCH system and RCH-equivalence from the
viewpoint of port Hamiltonian system with a symplectic structure.
These research work develop the theory of Hamiltonian reduction for
the regular controlled Hamiltonian systems with symmetries and make us
have much deeper understanding and recognition for the structures of
controlled Hamiltonian systems.

\section{Preliminaries}

In order to study the regular reduction theory of RCH systems, we
first give some relevant definitions and basic facts about momentum
maps, symplectic fiber bundle, Lie group lifted actions on
(co-)tangent bundles and reduction, which will be used in subsequent
sections, we shall follow the notations and conventions introduced
in Abraham et al. \cite{abma78,abmara88}, Marsden \cite{ma92},
Marsden et al. \cite{mamiorpera07}, Marsden and Ratiu \cite{mara99},
Ortega and Ratiu \cite{orra04}, Kobayashi and Nomizu \cite{kono63}.
In this paper, we assume that all manifolds are real, smooth and
finite dimensional and all actions are smooth left actions.

\noindent \subsection{Momentum map}

Let $(M,\omega)$ be a symplectic manifold, $G$ a Lie group with Lie
algebra $\mathfrak{g}$. We say that $G$ acts on $M$ and the action
of any $g\in G$ on $z \in M$ will be denoted by $\Phi:G \times
M\rightarrow M: \Phi(g,z)=g\cdot z$. For any $g\in G$, the map
$\Phi_g:=\Phi(g,\cdot): M \rightarrow M$ is a diffeomorphism of $M$
and if the map $\Phi_g $ satisfies $\Phi_g^\ast\omega=\omega, \;
\forall g \in G,$ we say that $G$ acts symplectically on a
symplectic manifold $(M,\omega)$. The isotropy subgroup of a point
$z\in M$ is $G_z=\{g\in G|\; g \cdot z=z\}. $ An action is free if
all the isotropy subgroups $G_z$ are trivial; and is proper if the
map $(g,z)\rightarrow  (g, g\cdot z)$ is proper (i.e., the pre-image
of every compact set is compact). For a proper action, all isotropy
subgroups are compact. The $G$-orbit of $z\in M$ is denoted by
$\mathcal{O}_z=G\cdot z=\{\Phi_g(z)|\; g\in G\},$ and the orbit
space by $M/G= \{\mathcal{O}_z |\; z\in M \}. $ If $G$ acts freely
and properly on $M$, then $M/G$ has a unique smooth structure such
that $\pi_G: M \rightarrow M/G$ is a surjective submersion. If $G$
acts only properly on $M$, does not act freely, then $M/G$ is not
necessarily smooth manifold, but just a quotient topological
space.\\

For each $\xi\in \mathfrak{g}$, the infinitesimal generator of $\xi$
is the vector field $\xi_M$ defined by
$\xi_M(z)=\left.\frac{\mathrm{d}}{\mathrm{d} t}\right|_{t=0}\exp
(t\xi)\cdot z,\forall z\in M$. We will also write $\xi_M(z)$ as
$\xi\cdot z$, and refer to the map $(\xi, z)\mapsto \xi\cdot z$ as
the infinitesimal action of $\mathfrak{g}$ on $M$. A momentum map
$\mathbf{J}:M \rightarrow \mathfrak{g}^\ast$ is defined by
$<\mathbf{J}(z),\xi>= J_\xi(z)$, for every $\xi\in \mathfrak{g}$,
where the function $J_\xi: M \rightarrow \mathbb{R}$ satisfies
$X_{J_\xi}=\xi_M$, and $\mathfrak{g}^\ast$ is the dual of Lie
algebra $\mathfrak{g}$, and $<,>: \mathfrak{g}^\ast \times
\mathfrak{g}\rightarrow \mathbb{R}$ is the duality pairing between
the dual $\mathfrak{g}^\ast $ and $\mathfrak{g}$. If the adjoint
action of $G$ on $\mathfrak{g}$ is denoted by $\operatorname{Ad}$,
and the infinitesimal adjoint action by $\operatorname{ad}$, then
the coadjoint action of $G$ on $g^\ast$ is the inverse dual to the
adjoint action, given by $g\cdot \nu=\operatorname{Ad}_{g^{-1}}^\ast
\nu=(\operatorname{Ad}_{g^{-1}})^\ast\nu,\forall\; \nu\in
\mathfrak{g}^\ast$. The infinitesimal coadjoint action is given by
$\xi\cdot \nu=-\operatorname{ad}_\xi^\ast\nu,\forall\; \nu\in
\mathfrak{g}^\ast$. For $\mu \in \mathfrak{g}^\ast$, a value of
$\mathbf{J}:M \rightarrow \mathfrak{g}^\ast$, $G_\mu$ denotes the
isotropy subgroup of $G$ with respect to the coadjoint $G$-action
$\operatorname{Ad}_{g^{-1}}^\ast$ at the point $\mu$, and
$\mathcal{O}_\mu$ denotes the $G$-orbit of through the point $\mu$
in $\mathfrak{g}^\ast$. The momentum map $\mathbf{J}$ is
$\operatorname{Ad}^\ast$-equivariant if
$\mathbf{J}(\Phi_g(z))=\operatorname{Ad}_{g^{-1}}^\ast
\mathbf{J}(z)$, for any $z\in M$.\\

The following proposition is very important for the regular
reduction and singular reduction of a Hamiltonian system with
symmetry; see Marsden \cite{ma92} and Ortega and Ratiu
\cite{orra04}.

\begin{proposition}
  (Bifurcation Lemma) Let $(M,\omega)$ be a symplectic manifold and
  $G$ a Lie group acting symplectically on $M$ (not necessarily
  freely). Suppose that the action has an associated momentum map
  $\mathbf{J}: M\to \mathfrak{g}^\ast$. Then for any $z\in M$,
  $(\mathfrak{g}_z)^0=\mbox{range}(T_z\mathbf{J})$, where $\mathfrak{g}_z=\{\xi\in
  \mathfrak{g}| \; \xi_M(z)=0\}$ is the Lie algebra of the isotropy
  subgroup $G_z =\{g\in G| \; g \cdot z=z\}$ and
  $(\mathfrak{g}_z)^0=\{\mu\in
  \mathfrak{g}^\ast| \; \mu|_{\mathfrak{g}_z}=0\}$ denotes the
  annihilator of $\mathfrak{g}_z$ in $\mathfrak{g}^\ast$.
\end{proposition}

An immediate consequence of this proposition is the fact that when
the action of $G$ is free, each value $\mu\in \mathfrak{g}^\ast$ of
the momentum map $\mathbf{J}$ is regular. Thus, if $\mu$ is a
singular value of $\mathbf{J}$, then the $G$-action is not free. In
addition, if $\mu$ is a regular value of $\mathbf{J}$ and
$\mathcal{O}_\mu$ is an embedded submanifold of $\mathfrak{g}^\ast$,
the $\mathbf{J}$ is transverse to $\mathcal{O}_\mu$ and hence
$\mathbf{J}^{-1}(\mathcal{O}_\mu)$ is automatically an embedded
submanifold of $M$. In this paper, we consider only that the
$G$-action is free, and the Hamiltonian reductions are regular.

\noindent \subsection{Symplectic fiber bundles }

Let $E$ and $M$ be two smooth manifolds, Lie group $G$ acts freely
on $E$ from the left side. Denote by $(E,M,\pi,G)$ a (left)
principal fiber bundle over $M$ with group $G$, where $E$ is the
bundle space, $M$ is the base space, $G$ is the structure group and
the projection $\pi: E \rightarrow M$ is a surjective submersion.
For each $x\in M$, $\pi^{-1}(x)$ is a closed submanifold of $E$,
which is called the fiber over $x$. Each fiber of the principal
bundle $(E,M,\pi,G)$ is diffeomorphic to $G$. In the following we
shall give a construction of the associated bundle of $G$-principal
bundle. Assume that $F$ is another smooth manifold and Lie group $G$
acts on $F$ from the left side. We can define a fiber bundle
associated to principal bundle $(E,M,\pi,G)$ with fiber $F$ as
follows. Consider the left action of $G$ on the product manifold
$E\times F$, $\Phi: G\times (E\times F)\rightarrow E\times F$ is
given by $\Phi(g,(z,y))= (gz, g^{-1}y),\; \forall\; g\in G,\; z\in
E,\; y\in F. $ Denote by $E\times_G F$ that is the orbit space
$(E\times F)/G$, and the map $\rho: E\times_G F \rightarrow M$ is
uniquely determined by the condition $\rho\cdot \pi_{/G} =\pi\cdot
\pi_{E}$, that is, the following commutative Diagram-1,
$$
\begin{CD}
 E\times F @> \pi_{/G}  >> E\times_G F \\
@V \pi_E VV @VV \rho  V \\
E @> \pi >> M
\end{CD}
$$
$$\mbox{Diagram-1}$$
where $\pi_{/G}: E\times F \rightarrow E\times_G F$ is the canonical
projection and $\pi_E: E\times F \rightarrow E$ is the projection
onto the first factor. Then $(E\times_G F,M,F,\rho,G)$, simply
written as $(E,M,F,\pi,G)$, is a fiber bundle with fiber $F$ and
structure group $G$ associated to principal bundle $(E,M,\pi,G)$. In
particular, if $F=V$ is a vector space, then $(E,M,V,\pi,G)$ is a
vector bundle associated to principal bundle $(E,M,\pi,G)$.

A bundle of symplectic manifolds is such a fiber bundle
$(E,M,F,\pi,G)$, all of whose fibers are symplectic and whose
structure group $G$ preserves the symplectic structure on $F$. From
Gotay et al. \cite{golasnwe83} we know that there exists a
presymplectic form $\omega_E$ on $E$ under some topological
conditions, whose pull-back to each fiber is the given fiber
symplectic form. We assume that if a symplectic form $\omega_E$ is
given on $E$, then $(E, \omega_E)$ is called a symplectic fiber
bundle. In particular, if $E$ is a vector bundle, then $(E,
\omega_E)$ is called a symplectic vector bundle; see Libermann and
Marle \cite{lima87}.

\noindent \subsection{Lie group lifted action on (co-)tangent
bundles and reduction }

For a smooth manifold $Q$, its cotangent bundle $T^\ast Q$ has a
canonical symplectic form $\omega_0$, which is given in natural
cotangent bundle coordinates $(q^i,p_i)$ by
$\omega_0=\mathbf{d}q^i\wedge \mathbf{d}p_i$, so $T^\ast Q$ is a
symplectic vector bundle. Let $\Phi: G\times Q \rightarrow Q$ be a
left smooth action of a Lie group $G$ on the manifold $Q$. The
tangent lift of this action $\Phi:G\times Q\rightarrow Q$ is the
action of $G$ on $TQ$, $\Phi^T:G\times TQ\rightarrow TQ$ given by
$g\cdot v_q =T\Phi_g(v_q),\;\forall\; v_q\in T_qQ, q\in Q$. The
cotangent lift is the action of $G$ on $T^\ast Q$,
$\Phi^{T^\ast}:G\times T^\ast Q\rightarrow T^\ast Q$ given by
$g\cdot \alpha_q=(T\Phi_{g^{-1}})^\ast\cdot
\alpha_q,\;\forall\;\alpha_q\in T^\ast_qQ,\; q\in Q$. The tangent or
cotangent lift of any proper (resp. free) $G$-action is proper
(resp. free). Each cotangent lift action is symplectic with respect
to the canonical symplectic form $\omega_0$, and has an
$\operatorname{Ad}^\ast$-equivariant momentum map $\mathbf{J}:T^\ast
Q\to \mathfrak{g}^\ast$ given by
$<\mathbf{J}(\alpha_q),\xi>=\alpha_q(\xi_Q(q)), $ where $\xi\in
\mathfrak{g}$, $\xi_Q(q)$ is the value of the infinitesimal
generator $\xi_Q$ of the $G$-action at $q\in Q$, $<,>:
\mathfrak{g}^\ast \times \mathfrak{g}\rightarrow \mathbb{R}$ is the
duality pairing between the dual $\mathfrak{g}^\ast $ and
$\mathfrak{g}$.\\

The reduction theory of cotangent bundle is a very important special
case of general reduction theory. Let $\mu\in \mathfrak{g}^\ast $ is
a regular value of the momentum map $\mathbf{J}$, the simplest case
of symplectic reduction of cotangent bundle $T^\ast Q$ is regular
point reduction at zero, in this case the symplectic reduced space
formed at $\mu=0$ is given by $((T^\ast Q)_\mu, \omega_\mu)=
(T^\ast(Q/G), \omega_0)$, where $\omega_0$ is the canonical
symplectic form of cotangent bundle $T^\ast (Q/G)$. Thus, the
reduced space $((T^\ast Q)_\mu, \omega_\mu)$ at $\mu=0$ is a
symplectic vector bundle. If $\mu\neq0$, from Marsden et al.
\cite{mamiorpera07} we know that, when $G_\mu=G$, the regular point
reduced space $((T^*Q)_\mu, \omega_\mu)$ is symplectically
diffeomorphic to symplectic vector bundle $(T^\ast (Q/G),
\omega_0-B_\mu)$, where $B_\mu$ is a magnetic term; If $G$ is not
Abelian and $G_\mu\neq G$, the regular point reduced space
$((T^*Q)_\mu, \omega_\mu)$ is symplectically diffeomorphic to a
symplectic fiber bundle over $T^\ast (Q/G_\mu)$ with fiber to be the
coadjoint orbit $\mathcal{O}_\mu$. In the case of regular orbit
reduction, from Ortega and Ratiu \cite{orra04} and the regular
reduction diagram, we know that the regular orbit reduced space
$((T^\ast Q)_{\mathcal{O}_\mu},\omega_{\mathcal{O}_\mu})$ is
symplectically diffeomorphic to the regular point reduced space
$((T^*Q)_\mu, \omega_\mu)$, and hence is symplectically
diffeomorphic to a symplectic fiber bundle.
Thus, the symplectic reduced space on a cotangent bundle may not be a
cotangent bundle, and hence the symplectic reduced system of a
Hamiltonian system with symmetry defined on the cotangent bundle
$T^*Q$ may not be a Hamiltonian system on a cotangent bundle.
To sum up above
discussion, if we may define a RCH system on a symplectic fiber
bundle, then it is possible to describe uniformly the RCH systems on
$T^*Q$ and their regular reduced RCH systems on the associated
reduced spaces.

\section{Regular Controlled Hamiltonian Systems}

In this paper, our goal is to study regular reduction theory of a RCH
system with symplectic structure and symmetry, as an extension of
the regular symplectic reduction theory of a Hamiltonian system under
regular controlled Hamiltonian equivalence conditions. Thus, in
order to describe uniformly RCH systems defined on a cotangent
bundle and on its regular reduced spaces, in this section we first
define a RCH system on a symplectic fiber bundle. In particular, we
obtain the RCH system by using the symplectic structure on the
cotangent bundle of a configuration manifold as a special case, and
discuss RCH-equivalence. In consequence, we can study the RCH
systems with symmetries by combining with regular symplectic reduction
theory of Hamiltonian systems. For convenience, we assume that all
controls appearing in this paper are the admissible controls.\\

Let $(E,M,N,\pi,G)$ be a fiber bundle and $(E, \omega_E)$ be a
symplectic fiber bundle. If for any function $H: E \rightarrow
\mathbb{R}$, we have a Hamiltonian vector field $X_H$ defined by
$i_{X_H}\omega_E=\mathbf{d}H$, then $(E, \omega_E, H )$ is a
Hamiltonian system. Moreover, if considering the external force and
control, we can define a kind of regular controlled Hamiltonian
(RCH) system on the symplectic fiber bundle $E$ as follows.

\begin{definition}
(RCH System) A RCH system on $E$ is a 5-tuple $(E, \omega_E, H, F,
W)$, where $(E, \omega_E, H )$ is a Hamiltonian system, and the
function $H: E \rightarrow \mathbb{R}$ is called the Hamiltonian, a
fiber-preserving map $F: E\rightarrow E$ is called the (external)
force map, and a fiber submanifold $W$ of $E$ is called the control
subset.
\end{definition}
Sometimes, $W$ also denotes the set of fiber-preserving maps from
$E$ to $W$. When a feedback control law $u:E\rightarrow W$ is
chosen, the 5-tuple $(E, \omega_E, H, F, u)$ denotes a closed-loop
dynamic system. In particular, when $Q$ is a smooth manifold, and
$T^\ast Q$ its cotangent bundle with a symplectic form $\omega$ (not
necessarily canonical symplectic form), then $(T^\ast Q, \omega )$
is a symplectic vector bundle. If we take that $E= T^* Q$, from
above definition we can obtain a RCH system on the cotangent bundle
$T^\ast Q$, that is, 5-tuple $(T^\ast Q, \omega, H, F, W)$. Where
the fiber-preserving map $F: T^*Q\rightarrow T^*Q$ is the (external)
force map, that is the reason that the fiber-preserving map $F:
E\rightarrow E$ is called an (external) force map in above
definition.\\

In order to describe the dynamics of the RCH system
$(E,\omega_E,H,F,W)$ with a control law $u$, we need to give a good
expression of the dynamical vector field of RCH system. At first, we
introduce a notations of vertical lift map of a vector along a
fiber. For a smooth manifold $E$, its tangent bundle $TE$ is a
vector bundle, and for the fiber bundle $\pi: E \rightarrow M$, we
consider the tangent mapping $T\pi: TE \rightarrow TM$ and its
kernel $ker (T\pi)=\{\rho\in TE| T\pi(\rho)=0\}$, which is a vector
subbundle of $TE$. Denote $VE:= ker(T\pi)$, which is called a
vertical bundle of $E$. Assume that there is a metric on $E$, and we
take a Levi-Civita connection $\mathcal{A}$ on $TE$, and denote
$HE:= ker(\mathcal{A})$, which is called a horizontal bundle of $E$,
such that $TE= HE \oplus VE. $ For any $x\in M, \; a_x, b_x \in E_x,
$ any tangent vector $\rho(b_x)\in T_{b_x}E$ can be split into
horizontal and vertical parts, that is, $\rho(b_x)=
\rho^h(b_x)\oplus \rho^v(b_x)$, where $\rho^h(b_x)\in H_{b_x}E$ and
$\rho^v(b_x)\in V_{b_x}E$. Let $\gamma$ be a geodesic in $E_x$
connecting $a_x$ and $b_x$, and denote $\rho^v_\gamma(a_x)$ a
tangent vector at $a_x$, which is a parallel displacement of the
vertical vector $\rho^v(b_x)$ along the geodesic $\gamma$ from $b_x$
to $a_x$. Since the angle between two vectors is invariant under a
parallel displacement along a geodesic, then
$T\pi(\rho^v_\gamma(a_x))=0, $ and hence $\rho^v_\gamma(a_x) \in
V_{a_x}E. $ Now, for $a_x, b_x \in E_x $ and tangent vector
$\rho(b_x)\in T_{b_x}E$, we can define the vertical lift map of a
vector along a fiber given by
$$\mbox{vlift}: TE_x \times E_x \rightarrow TE_x; \;\; \mbox{vlift}(\rho(b_x),a_x) = \rho^v_\gamma(a_x). $$
It is easy to check from the basic fact in differential geometry
that this map does not depend on the choice of $\gamma$. If $F: E
\rightarrow E$ is a fiber-preserving map, for any $x\in M$, we have
that $F_x: E_x \rightarrow E_x$ and $TF_x: TE_x \rightarrow TE_x$,
then for any $a_x \in E_x$ and $\rho\in TE_x$, the vertical lift of
$\rho$ under the action of $F$ along a fiber is defined by
$$(\mbox{vlift}(F_x)\rho)(a_x)=\mbox{vlift}((TF_x\rho)(F_x(a_x)), a_x)= (TF_x\rho)^v_\gamma(a_x), $$
where $\gamma$ is a geodesic in $E_x$ connecting $F_x(a_x)$ and
$a_x$.\\

In particular, when $\pi: E \rightarrow M$ is a vector bundle, for
any $x\in M$, the fiber $E_x=\pi^{-1}(x)$ is a vector space. In this
case, we can choose the geodesic $\gamma$ to be a straight line, and
the vertical vector is invariant under a parallel displacement along
a straight line, that is, $\rho^v_\gamma(a_x)= \rho^v(b_x).$
Moreover, when $E= T^*Q$, by using the local trivialization of
$TT^*Q$, we have that $TT^*Q\cong TQ \times T^*Q$. Because of $\pi:
T^*Q \rightarrow Q$, and $T\pi: TT^*Q \rightarrow TQ$, then in this
case, for any $\alpha_x, \; \beta_x \in T^*_x Q, \; x\in Q, $ we
know that $(0, \beta_x) \in V_{\beta_x}T^*_x Q, $ and hence we can
get that
$$ \mbox{vlift}((0, \beta_x)(\beta_x), \alpha_x) = (0, \beta_x)(\alpha_x)
= \left.\frac{\mathrm{d}}{\mathrm{d}s}\right|_{s=0}(\alpha_x+s\beta_x), $$
which is consistent with the definition of vertical lift map along
fiber in Marsden and Ratiu \cite{mara99}.\\

For a given RCH System $(T^\ast Q, \omega, H, F, W)$, the dynamical
vector field of the associated Hamiltonian system $(T^\ast Q,
\omega, H) $ is that $X_H= (\mathbf{d}H)^\sharp$, where, $\sharp:
T^\ast T^\ast Q \rightarrow TT^\ast Q; \mathbf{d}H \mapsto
(\mathbf{d}H)^\sharp$, such that
$i_{(\mathbf{d}H)^\sharp}\omega=\mathbf{d}H$. If considering the
external force $F: T^*Q \rightarrow T^*Q, $ by using the above
notations of vertical lift map of a vector along a fiber, the
change of $X_H$ under the action of $F$ is that
$$\mbox{vlift}(F)X_H(\alpha_x)= \mbox{vlift}((TFX_H)(F(\alpha_x)), \alpha_x)= (TFX_H)^v_\gamma(\alpha_x),$$
where $\alpha_x \in T^*_x Q, \; x\in Q $ and $\gamma$ is a straight
line in $T^*_x Q$ connecting $F_x(\alpha_x)$ and $\alpha_x$. In the
same way, when a feedback control law $u: T^\ast Q \rightarrow W$ is
chosen, the change of $X_H$ under the action of $u$ is that
$$\mbox{vlift}(u)X_H(\alpha_x)= \mbox{vlift}((TuX_H)(u(\alpha_x)), \alpha_x)= (TuX_H)^v_\gamma(\alpha_x).$$
In consequence, the dynamical vector field of a RCH system $(T^\ast
Q,\omega,H,F,W)$ with a control law $u$ is the synthetic of
Hamiltonian vector field $X_H$ and its changes under the actions of
the external force $F$ and control $u$, that is,
$$X_{(T^\ast Q,\omega,H,F,u)}(\alpha_x)= X_H(\alpha_x)+ \mbox{vlift}(F)X_H(\alpha_x)+ \mbox{vlift}(u)X_H(\alpha_x),$$
for any $\alpha_x \in T^*_x Q, \; x\in Q $. For convenience, it is
simply written as
\begin{equation}X_{(T^\ast Q,\omega,H,F,u)}
=(\mathbf{d}H)^\sharp+\textnormal{vlift}(F)+\textnormal{vlift}(u).\label{3.1}\end{equation}
We also denote that $\mbox{vlift}(W)= \bigcup\{\mbox{vlift}(u)X_H |
\; u\in W\}$. For the RCH system $(E,\omega_E,H,F,W)$ with a control
law $u$, we have also a similar expression of its dynamical vector
field. It is worthy of note that in order to deduce and calculate
easily, we always use the simple expression of dynamical vector
field $X_{(T^\ast Q,\omega,H,F,u)}$. Moreover, we also use the
simple expressions for $R_P$-reduced vector field $X_{((T^\ast
Q)_\mu, \omega_\mu, h_\mu, f_\mu, u_\mu)}$ and $R_O$-reduced vector
field $X_{((T^\ast Q)_{\mathcal{O}_\mu}, \omega_{\mathcal{O}_\mu},
h_{\mathcal{O}_\mu},f_{\mathcal{O}_\mu},u_{\mathcal{O}_\mu})}$ in
Section 4 and Section 5.\\

Next, we note that when a RCH system is given, the force map $F$ is
determined, but the feedback control law $u: T^\ast Q\rightarrow W$
could be chosen. In order to describe the feedback control law to
modify the structure of RCH system, the Hamiltonian matching
conditions and RCH-equivalence are induced as follows.
\begin{definition}
(RCH-equivalence) Suppose that we have two RCH systems $(T^\ast
Q_i,\omega_i,H_i,F_i,W_i),\; i= 1,2,$ we say them to be
RCH-equivalent, or simply, $(T^\ast
Q_1,\omega_1,H_1,F_1,W_1)\stackrel{RCH}{\sim}(T^\ast
Q_2,\omega_2,H_2,F_2,W_2)$, if there exists a diffeomorphism
$\varphi: Q_1\rightarrow Q_2$, such that the following Hamiltonian
matching conditions hold:

\noindent {\bf RHM-1:} The cotangent lift map of $\varphi$, that is,
$\varphi^\ast= T^\ast \varphi:T^\ast Q_2\rightarrow T^\ast Q_1$ is
symplectic, and $W_1=\varphi^\ast (W_2).$

\noindent {\bf RHM-2:}
$Im[(\mathbf{d}H_1)^\sharp+\textnormal{vlift}(F_1)-((\varphi_\ast)^\ast
\mathbf{d}H_2)^\sharp-\textnormal{vlift}(\varphi^\ast
F_2\varphi_\ast)]\subset \textnormal{vlift}(W_1)$, where the map
$\varphi_\ast=(\varphi^{-1})^\ast: T^\ast Q_1\rightarrow T^\ast
Q_2$, and $(\varphi^\ast)_\ast=(\varphi_\ast)^\ast=T^\ast
\varphi_\ast: T^\ast T^\ast Q_2\rightarrow T^\ast T^\ast Q_1$, and
$Im$ means the pointwise image of the map in brackets.
\end{definition}
It is worthy of noting that our RCH system is defined by using the
symplectic structure on the cotangent bundle of a configuration
manifold, we must keep with the symplectic structure when we define
the RCH-equivalence, that is, the induced equivalent map $\varphi^*$
is symplectic on the cotangent bundle. In the same way, for the RCH
systems on the symplectic fiber bundles, we can also define the
RCH-equivalence by replacing $T^\ast Q_i$ and $\varphi:
Q_1\rightarrow Q_2$ by $E_i$ and $\varphi^\ast: E_2\rightarrow E_1$,
respectively. Moreover, the following Theorem 3.3 explains the
significance of the above RCH-equivalence relation.
\begin{theorem}
Suppose that two RCH systems $(T^\ast Q_i,\omega_i,H_i,F_i,W_i)$,
$i=1,2,$ are RCH-equivalent, then there exist two control laws $u_i:
T^\ast Q_i \rightarrow W_i, \; i=1,2, $ such that the two
closed-loop systems produce the same equations of motion, that is,
$X_{(T^\ast Q_1,\omega_1,H_1,F_1,u_1)}\cdot \varphi^\ast
=T(\varphi^\ast) X_{(T^\ast Q_2,\omega_2,H_2,F_2,u_2)}$, where the
map $T(\varphi^\ast):TT^\ast Q_2\rightarrow TT^\ast Q_1$ is the
tangent map of $\varphi^\ast$. Moreover, the explicit relation
between the two control laws $u_i, i=1,2$ is given by
\begin{equation}\textnormal{vlift}(u_1) -\textnormal{vlift}(\varphi^\ast
u_2\varphi_\ast)=-(\mathbf{d}H_1)^\sharp
-\textnormal{vlift}(F_1)+((\varphi_\ast)^\ast
\mathbf{d}H_2)^\sharp+\textnormal{vlift}(\varphi^\ast F_2
\varphi_\ast)\label{3.2}\end{equation}
\end{theorem}

{\bf Proof:} From (\ref{3.1}), we have that $X_{(T^\ast
Q_1,\omega_1,H_1,F_1,u_1)}
=(\mathbf{d}H_1)^\sharp+\textnormal{vlift}(F_1)+\textnormal{vlift}(u_1)$
and
\begin{align*}
T(\varphi^\ast) X_{(T^\ast
Q_2,\omega_2,H_2,F_2,u_2)}&=T(\varphi^\ast)[(\mathbf{d}H_2)^\sharp
+\textnormal{vlift}(F_2)+\textnormal{vlift}(u_2)]\\
&=T(\varphi^\ast)(\mathbf{d}H_2)^\sharp
+T(\varphi^\ast)\textnormal{vlift}(F_2)+T(\varphi^\ast)\textnormal{vlift}(u_2)
\end{align*}
\[
\begin{CD}
T^\ast Q_2 @>\textnormal{vlift}>>TT^\ast Q_2@<\sharp<<T^\ast T^\ast Q_2\\
@V\varphi^\ast VV @VT \varphi^\ast  VV @V(\varphi_\ast)^\ast VV\\
T^\ast Q_1 @>>\textnormal{vlift}> TT^\ast Q_1 @<<\sharp<T^\ast
T^\ast Q_1
\end{CD}
\]
$$\mbox{Diagram-2}$$
From the commutative Diagram-2 and the definition of the vertical
lift operator $\mathrm{vlift}$, we have that for $\alpha \in T^\ast
Q_2$,
\begin{align*}
&T(\varphi^\ast)\textnormal{vlift}(F_2)(\alpha)
=T(\varphi^\ast)\left.\frac{\mathrm{d}}{\mathrm{d}s}\right|_{s=0}(\alpha+sF_2(\alpha))\\
&=\left.\frac{\mathrm{d}}{\mathrm{d}s}\right|_{s=0}(\varphi^\ast\alpha+s\varphi^\ast
F_2\varphi_\ast(\varphi^\ast
\alpha))=\textnormal{vlift}(\varphi^\ast
F_2\varphi_\ast)(\varphi^\ast\alpha).
\end{align*}
In the same way, we have that
$T(\varphi^\ast)\textnormal{vlift}(u_2)=\textnormal{vlift}(\varphi^\ast
u_2\varphi_\ast)\cdot \varphi^\ast$. Since $\varphi^\ast:T^\ast
Q_2\rightarrow T^\ast Q_1$ is symplectic, and
$i_{(\mathbf{d}H_i)^\sharp}\omega_i=\mathbf{d}H_i$, we have that
$T(\varphi^\ast)(\mathbf{d}H_2)^\sharp=((\varphi_\ast)^\ast
\mathbf{d}H_2)^\sharp\cdot \varphi^\ast $. Thus,
$$T(\varphi^\ast) X_{(T^\ast
Q_2,\omega_2,H_2,F_2,u_2)}=((\varphi_\ast)^\ast
\mathbf{d}H_2)^\sharp\cdot
\varphi^\ast+\textnormal{vlift}(\varphi^\ast F_2\varphi_\ast)\cdot
\varphi^\ast+\textnormal{vlift}(\varphi^\ast u_2\varphi_\ast)\cdot
\varphi^\ast.$$ From $X_{(T^\ast Q_1,\omega_1,H_1,F_1,u_1)}\cdot
\varphi^\ast=T(\varphi^\ast) X_{(T^\ast Q_2,\omega_2,H_2,F_2,u_2)}$,
we have that (\ref{3.2}) holds. \hskip 1cm $\blacksquare$ \\

In the following we shall introduce the regular point and regular
orbit reducible RCH systems with symplectic forms and symmetries, and
show a variety of relationships of their regular reduced
RCH-equivalences.

\section{Regular Point Reduction of RCH Systems }

Let $Q$ be a smooth manifold and $T^\ast Q$ its cotangent bundle
with the symplectic form $\omega$. Let $\Phi:G\times Q\rightarrow Q$
be a smooth left action of the Lie group $G$ on $Q$, which is free
and proper. Then the cotangent lifted left action
$\Phi^{T^\ast}:G\times T^\ast Q\rightarrow T^\ast Q$ is also
free and proper. Assume that this $G$-action is symplectic
and admits an $\operatorname{Ad}^\ast$-equivariant
momentum map $\mathbf{J}:T^\ast Q\rightarrow \mathfrak{g}^\ast$. If
$\mu\in\mathfrak{g}^\ast$ is a regular value of $\mathbf{J}$ and
$G_\mu=\{g\in G|\operatorname{Ad}_g^\ast \mu=\mu \}$ is the isotropy
subgroup of coadjoint $G$-action at the point $\mu$. Since $G_\mu
(\subset G)$ acts freely and properly on $Q$ and on $T^\ast Q$, then
$Q_\mu=Q/G_\mu$ is a smooth manifold and that the canonical
projection $\rho_\mu:Q\rightarrow Q_\mu$ is a surjective submersion.
It follows that $G_\mu$ acts also freely and properly on
$\mathbf{J}^{-1}(\mu)$, so that the space $(T^\ast
Q)_\mu=\mathbf{J}^{-1}(\mu)/G_\mu$ is a symplectic manifold with
symplectic form $\omega_\mu$ uniquely characterized by the relation
\begin{equation}\pi_\mu^\ast \omega_\mu=i_\mu^\ast
\omega. \label{4.1}\end{equation} The map
$i_\mu:\mathbf{J}^{-1}(\mu)\rightarrow T^\ast Q$ is the inclusion
and $\pi_\mu:\mathbf{J}^{-1}(\mu)\rightarrow (T^\ast Q)_\mu$ is the
projection. The pair $((T^\ast Q)_\mu,\omega_\mu)$ is called the
symplectic point reduced space of $(T^\ast Q,\omega)$ at $\mu$.

\begin{remark}
If $(T^\ast Q, \omega)$ is a connected symplectic manifold, and
$\mathbf{J}:T^\ast Q\rightarrow \mathfrak{g}^\ast$ is a
non-equivariant momentum map with a non-equivariance group
one-cocycle $\sigma: G\rightarrow \mathfrak{g}^\ast$, which is
defined by $\sigma(g):=\mathbf{J}(g\cdot
z)-\operatorname{Ad}^\ast_{g^{-1}}\mathbf{J}(z)$, where $g\in G$ and
$z\in T^\ast Q$. Then we know that $\sigma$ produces a new affine
action $\Theta: G\times \mathfrak{g}^\ast \rightarrow
\mathfrak{g}^\ast $ defined by
$\Theta(g,\mu):=\operatorname{Ad}^\ast_{g^{-1}}\mu + \sigma(g)$,
where $\mu \in \mathfrak{g}^\ast$, with respect to which the given
momentum map $\mathbf{J}$ is equivariant. Assume that $G$ acts
freely and properly on $T^\ast Q$, and $\tilde{G_\mu}$ denotes the
isotropy subgroup of $\mu \in \mathfrak{g}^\ast$ relative to this
affine action $\Theta$ and $\mu$ is a regular value of $\mathbf{J}$.
Then the quotient space $(T^\ast
Q)_\mu=\mathbf{J}^{-1}(\mu)/\tilde{G_\mu}$ is also a symplectic
manifold with symplectic form $\omega_\mu$ uniquely characterized by
(\ref{4.1}), see Ortega and Ratiu \cite{orra04}.
\end{remark}

Let $H: T^\ast Q\rightarrow \mathbb{R}$ be a $G$-invariant
Hamiltonian, the flow $F_t$ of the Hamiltonian vector field $X_H$
leaves the connected components of $\mathbf{J}^{-1}(\mu)$ invariant
and commutes with the $G$-action, then it induces a flow $f_t^\mu$
on $(T^\ast Q)_\mu$, defined by $f_t^\mu\cdot \pi_\mu=\pi_\mu \cdot
F_t\cdot i_\mu$, and the vector field $X_{h_\mu}$ generated by the
flow $f_t^\mu$ on $((T^\ast Q)_\mu,\omega_\mu)$ is Hamiltonian with
the associated regular point reduced Hamiltonian function
$h_\mu:(T^\ast Q)_\mu\rightarrow \mathbb{R}$ defined by
$h_\mu\cdot\pi_\mu=H\cdot i_\mu$, and the Hamiltonian vector fields
$X_H$ and $X_{h_\mu}$ are $\pi_\mu$-related. On the other hand, from
section 2, we know that the regular point reduced space
$((T^*Q)_\mu, \omega_\mu)$ is symplectically diffeomorphic to a
symplectic fiber bundle. Thus, we can introduce a regular point
reducible RCH systems as follows.

\begin{definition}
(Regular Point Reducible RCH System) A 6-tuple $(T^\ast Q, G,
\omega, H, F, W)$, where the Hamiltonian $H:T^\ast Q\rightarrow
\mathbb{R}$, the fiber-preserving map $F:T^\ast Q\rightarrow T^\ast
Q$ and the fiber submanifold $W$ of\; $T^\ast Q$ are all
$G$-invariant, is called a regular point reducible RCH system, if
there exists a point $\mu\in\mathfrak{g}^\ast$, which is a regular
value of the momentum map $\mathbf{J}$, such that the regular point
reduced system, that is, the 5-tuple $((T^\ast Q)_\mu,
\omega_\mu,h_\mu,f_\mu,W_\mu)$, where $(T^\ast
Q)_\mu=\mathbf{J}^{-1}(\mu)/G_\mu$, $\pi_\mu^\ast
\omega_\mu=i_\mu^\ast\omega$, $h_\mu\cdot \pi_\mu=H\cdot i_\mu$,
$F(\mathbf{J}^{-1}(\mu))\subset \mathbf{J}^{-1}(\mu) $, $f_\mu\cdot
\pi_\mu=\pi_\mu \cdot F\cdot i_\mu$, $W \cap
\mathbf{J}^{-1}(\mu)\neq \emptyset $, $W_\mu=\pi_\mu(W\cap
\mathbf{J}^{-1}(\mu))$, is a RCH system, which is simply written as
$R_P$-reduced RCH system. Where $((T^\ast Q)_\mu,\omega_\mu)$ is the
$R_P$-reduced space, the function $h_\mu:(T^\ast Q)_\mu\rightarrow
\mathbb{R}$ is called the reduced Hamiltonian, the fiber-preserving
map $f_\mu:(T^\ast Q)_\mu\rightarrow (T^\ast Q)_\mu$ is called the
reduced (external) force map, $W_\mu$ is a fiber submanifold of
\;$(T^\ast Q)_\mu$ and is called the reduced control subset.
\end{definition}

It is worthy of noting that for the regular point reducible RCH system
$(T^\ast Q,G,\omega,H,F,W)$, the $G$-invariant external force map
$F: T^*Q \rightarrow T^*Q $ has to satisfy the conditions
$F(\mathbf{J}^{-1}(\mu))\subset \mathbf{J}^{-1}(\mu), $ and
$f_\mu\cdot \pi_\mu=\pi_\mu \cdot F\cdot i_\mu, $ such that we can
define the reduced external force map $f_\mu:(T^\ast
Q)_\mu\rightarrow (T^\ast Q)_\mu. $ The condition $W \cap
\mathbf{J}^{-1}(\mu)\neq \emptyset $ in above definition makes that
the $G$-invariant control subset $W\cap \mathbf{J}^{-1}(\mu)$ can be
reduced and the reduced control subset is $W_\mu= \pi_\mu(W\cap
\mathbf{J}^{-1}(\mu))$. If the control subset cannot be reduced, we
cannot get the $R_P$-reduced RCH system. The study of RCH system
which is not regular point reducible is beyond the limits in this
paper, it may be a topic in future study.\\

Denote by $X_{(T^\ast Q,G,\omega,H,F,u)}$ the vector field of a
regular point reducible RCH system $(T^\ast Q,G,\omega,H,\\F,W)$
with a control law $u$. Assume that we have that
\begin{equation}X_{(T^\ast Q,G,\omega,H,F,u)}
=(\mathbf{d}H)^\sharp+\textnormal{vlift}(F)+\textnormal{vlift}(u).\label{4.2}\end{equation}
Then, for the regular point reducible RCH system we can also
introduce the regular point reduced controlled Hamiltonian
equivalence (RpCH-equivalence) as follows.

\begin{definition}(RpCH-equivalence)
Suppose that we have two regular point reducible RCH systems
$(T^\ast Q_i, G_i,\omega_i,H_i, F_i, W_i),\; i=1,2$, we say them to
be RpCH-equivalent, or simply,\\ $(T^\ast Q_1,
G_1,\omega_1,H_1,F_1,W_1)\stackrel{RpCH}{\sim}(T^\ast
Q_2,G_2,\omega_2,H_2,F_2,W_2)$, if there exists a diffeomorphism
$\varphi:Q_1\rightarrow Q_2$, such that the following Hamiltonian
matching conditions hold:

\noindent {\bf RpHM-1:} The cotangent lift map
$\varphi^\ast:T^\ast Q_2\rightarrow T^\ast Q_1$ is symplectic.

\noindent {\bf RpHM-2:} For $\mu_i\in \mathfrak{g}^\ast_i $, the
regular reducible points of RCH systems $(T^\ast Q_i, G_i,\omega_i,
H_i, F_i, W_i),\; i=1,2$, the map
$\varphi_\mu^\ast=i_{\mu_1}^{-1}\cdot\varphi^\ast\cdot i_{\mu_2}:
\mathbf{J}_2^{-1}(\mu_2)\rightarrow \mathbf{J}_1^{-1}(\mu_1)$ is
$(G_{2\mu_2},G_{1\mu_1})$-equivariant and $W_1\cap
\mathbf{J}_1^{-1}(\mu_1)=\varphi_\mu^\ast (W_2\cap
\mathbf{J}_2^{-1}(\mu_2))$, where $\mu=(\mu_1, \mu_2)$, and denote
by $i_{\mu_1}^{-1}(S)$ the pre-image of a subset $S\subset T^\ast
Q_1$ for the map $i_{\mu_1}:\mathbf{J}_1^{-1}(\mu_1)\rightarrow
T^\ast Q_1$.

\noindent {\bf RpHM-3:}
$Im[(\mathbf{d}H_1)^\sharp+\textnormal{vlift}(F_1)-((\varphi_\ast)^\ast
\mathbf{d}H_2)^\sharp-\textnormal{vlift}(\varphi^\ast
F_2\varphi_\ast)]\subset\textnormal{vlift}(W_1)$.
\end{definition}

It is worthy of noting that for the regular point reducible RCH
system, the induced equivalent map $\varphi^*$ not only keeps the
symplectic structure, but also keeps the equivariance of $G$-action
at the regular point. If a feedback control law $u_\mu:(T^\ast
Q)_\mu\rightarrow W_\mu$ is chosen, the $R_P$-reduced RCH system
$((T^\ast Q)_\mu, \omega_\mu, h_\mu, f_\mu, u_\mu)$ is a closed-loop
regular dynamic system with a control law $u_\mu$. Assume that its
vector field $X_{((T^\ast Q)_\mu, \omega_\mu, h_\mu, f_\mu, u_\mu)}$
can be expressed by
\begin{equation}X_{((T^\ast Q)_\mu, \omega_\mu, h_\mu, f_\mu, u_\mu)}
=(\mathbf{d}h_\mu)^\sharp+\textnormal{vlift}(f_\mu)+\textnormal{vlift}(u_\mu),
\label{4.3}\end{equation} where $(\mathbf{d}h_\mu)^\sharp =
X_{h_\mu}$, $\textnormal{vlift}(f_\mu)=
\textnormal{vlift}(f_\mu)X_{h_\mu}$, $\textnormal{vlift}(u_\mu)=
\textnormal{vlift}(u_\mu)X_{h_\mu}$, and satisfies the condition
\begin{equation}X_{((T^\ast Q)_\mu, \omega_\mu, h_\mu, f_\mu,
u_\mu)}\cdot \pi_\mu=T\pi_\mu\cdot X_{(T^\ast
Q,G,\omega,H,F,u)}\cdot i_\mu. \label{4.4}\end{equation} Then we can
obtain the following regular point reduction theorem for the RCH system,
which explains the relationship between the RpCH-equivalence for
the regular point reducible RCH system with symmetry and the
RCH-equivalence for the associated $R_P$-reduced RCH system. This
theorem can be regarded as an extension of regular point reduction
theorem of a Hamiltonian system under regular controlled Hamiltonian
equivalence condition.

\begin{theorem}
Two regular point reducible RCH systems $(T^\ast Q_i, G_i, \omega_i,
H_i, F_i,W_i)$, $i=1,2,$ are RpCH-equivalent if and only if the
associated $R_P$-reduced RCH systems $((T^\ast
Q_i)_{\mu_i},\omega_{i\mu_i},h_{i\mu_i},f_{i\mu_i},\\ W_{i\mu_i}),
i=1,2,$ are RCH-equivalent.
\end{theorem}

{\bf Proof:} If $(T^\ast Q_1, G_1, \omega_1, H_1, F_1,
W_1)\stackrel{RpCH}{\sim}(T^\ast Q_2, G_2, \omega_2, H_2, F_2,
W_2)$, then there exists a diffeomorphism $\varphi:Q_1\rightarrow
Q_2$ such that $\varphi^\ast:T^\ast Q_2\rightarrow T^\ast Q_1$ is
symplectic and for $\mu_i\in \mathfrak{g}^\ast_i, i=1,2$,
$\varphi_\mu^\ast=i_{\mu_1}^{-1}\cdot \varphi^\ast \cdot
i_{\mu_2}:\mathbf{J}_2^{-1}(\mu_2)\rightarrow
\mathbf{J}_1^{-1}(\mu_1)$ is $(G_{2\mu_2},G_{1\mu_1})$-equivariant,
$W_1\cap \mathbf{J}_1^{-1}(\mu_1)=\varphi_\mu^\ast (W_2\cap
\mathbf{J}_2^{-1}(\mu_2))$ and RpHM-3 holds. From the following
commutative Diagram-3:
\[
\begin{CD}
T^\ast Q_2 @<i_{\mu_2}<< \mathbf{J}_2^{-1}(\mu_2) @>\pi_{\mu_2}>> (T^\ast Q_2)_{\mu_2}\\
@V\varphi^\ast VV @V\varphi^\ast_\mu VV @V\varphi^\ast_{\mu/G}VV\\
T^\ast Q_1 @<i_{\mu_1}<< \mathbf{J}_1^{-1}(\mu_1)
@>\pi_{\mu_1}>>(T^\ast Q_1)_{\mu_1}
\end{CD}
\]
$$\mbox{Diagram-3}$$
We can define a map $\varphi_{\mu/G}^\ast:(T^\ast
Q_2)_{\mu_2}\rightarrow (T^\ast Q_1)_{\mu_1}$ such that
$\varphi_{\mu/G}^\ast \cdot
\pi_{\mu_2}=\pi_{\mu_1}\cdot\varphi^\ast_\mu$. Because
$\varphi_\mu^\ast: \mathbf{J}_2^{-1}(\mu_2)\rightarrow
\mathbf{J}_1^{-1}(\mu_1)$ is $(G_{2\mu_2},G_{1\mu_1})$-equivariant,
$\varphi_{\mu/G}^\ast$ is well-defined. We shall show that
$\varphi_{\mu/G}^\ast$ is symplectic and
$W_{1\mu_1}=\varphi_{\mu/G}^\ast (W_{2\mu_2})$. In fact, since
$\varphi^\ast: T^\ast Q_2\rightarrow T^\ast Q_1$ is symplectic, the
map $(\varphi^\ast)^\ast:\Omega^2(T^\ast Q_1)\rightarrow
\Omega^2(T^\ast Q_2)$ satisfies $(\varphi^\ast)^\ast
\omega_1=\omega_2$. By (\ref{4.1}),
$i_{\mu_i}^\ast\omega_i=\pi_{\mu_i}^\ast\omega_{i\mu_i}, i=1,2$,
from the following commutative Diagram-4,
\[
\begin{CD}
\Omega^2(T^\ast Q_1) @ >i_{\mu_1}^\ast>>
\Omega^2(\mathbf{J}_1^{-1}(\mu_1)) @ <\pi_{\mu_1}^\ast<< \Omega^2((T^\ast Q_1)_{\mu_1})\\
@V(\varphi^\ast)^\ast VV @V(\varphi^\ast_\mu)^\ast VV @V(\varphi^\ast_{\mu/G})^\ast VV\\
\Omega^2(T^\ast Q_2) @>i_{\mu_2}^\ast>>
\Omega^2(\mathbf{J}_2^{-1}(\mu_2)) @<\pi_{\mu_2}^\ast
<<\Omega^2((T^\ast Q_2)_{\mu_2})
\end{CD}
\]
$$\mbox{Diagram-4}$$
we have that
\begin{align*}
\pi_{\mu_2}^\ast \cdot(\varphi_{\mu/G}^\ast)^\ast\omega_{1\mu_1}&
=(\varphi_{\mu/G}^\ast\cdot \pi_{\mu_2})^\ast\omega_{1\mu_1}
=(\pi_{\mu_1}\cdot \varphi_{\mu}^\ast)^\ast\omega_{1\mu_1}
=(i_{\mu_1}^{-1}\cdot \varphi^\ast \cdot i_{\mu_2})^\ast
\cdot\pi_{\mu_1}^\ast\omega_{1\mu_1}\\
&=i_{\mu_2}^\ast\cdot(\varphi^\ast)^\ast\cdot(i_{\mu_1}^{-1})^\ast
\cdot i_{\mu_1}^\ast \omega_1=i_{\mu_2}^\ast\cdot
(\varphi^\ast)^\ast\omega_1
=i_{\mu_2}^\ast\omega_2=\pi_{\mu_2}^\ast\omega_{2\mu_2}.
\end{align*}
Notice that $\pi_{\mu_2}^\ast$ is a surjective, thus,
$(\varphi_{\mu/G}^\ast)^\ast\omega_{1\mu_1}=\omega_{2\mu_2}$.
Because by hypothesis $W_i\cap \mathbf{J}_i^{-1}(\mu_i)\neq
\emptyset $, $W_{i\mu_i}=\pi_{\mu_i}(W_i\cap
\mathbf{J}_i^{-1}(\mu_i)),\; i=1,2$ and $W_1\cap
\mathbf{J}_1^{-1}(\mu_1)=\varphi_\mu^\ast (W_2\cap
\mathbf{J}_2^{-1}(\mu_2))$, we have that
\begin{align*}
W_{1\mu_1} & =\pi_{\mu_1}(W_1\cap \mathbf{J}_1^{-1}(\mu_1))
=\pi_{\mu_1}\cdot \varphi_\mu^\ast (W_2\cap
\mathbf{J}_2^{-1}(\mu_2))\\ & =\varphi_{\mu/G}^\ast\cdot
\pi_{\mu_2}(W_2\cap \mathbf{J}_2^{-1}(\mu_2))
=\varphi_{\mu/G}^\ast(W_{2\mu_2}).
\end{align*}
Next, from
(\ref{4.2}) and (\ref{4.3}), we know that for $i=1,2$,
$$X_{(T^\ast Q_i, G_i, \omega_i, H_i, F_i, u_i)}
=(\mathbf{d}H_i)^\sharp+\textnormal{vlift}(F_i)+\textnormal{vlift}(u_i),$$
$$X_{((T^\ast Q_i)_{\mu_i}, \omega_{i\mu_i}, h_{i\mu_i}, f_{i\mu_i},u_{i\mu_i})}
=(\mathbf{d}h_{i\mu_i})^\sharp+\textnormal{vlift}(f_{i\mu_i})+\textnormal{vlift}(u_{i\mu_i}),$$
and from (\ref{4.4}), we have that
$$X_{((T^\ast Q_i)_{\mu_i}, \omega_{i\mu_i}, h_{i\mu_i},
f_{i\mu_i},u_{i\mu_i})}\cdot \pi_{\mu_i}=T\pi_{\mu_i}\cdot
X_{(T^\ast Q_i, G_i, \omega_i, H_i, F_i, u_i)}\cdot i_{\mu_i}.$$
Since $H_i,F_i$ and $W_i$ are all $G_i$-invariant, $i=1,2$ and
$$h_{i\mu_i}\cdot \pi_{\mu_i}=H_i\cdot
i_{\mu_i},\;\;f_{i\mu_i}\cdot\pi_{\mu_i}=\pi_{\mu_i}\cdot F_i\cdot
i_{\mu_i},\;\;u_{i\mu_i}\cdot\pi_{\mu_i}=\pi_{\mu_i}\cdot u_i\cdot
i_{\mu_i}, \;\;i=1,2.$$ From the following commutative Diagram-5,
\[
\begin{CD}
T^\ast T^\ast Q_2 @>i_{\mu_2}^\ast>>
T^\ast \mathbf{J}_2^{-1}(\mu_2) @<\pi_{\mu_2}^\ast<< T^\ast ((T^\ast Q_2)_{\mu_2})\\
@V(\varphi^\ast)_\ast VV @V(\varphi^\ast_\mu)_\ast VV @V(\varphi^\ast_{\mu/G})_\ast VV\\
T^\ast T^\ast Q_1 @>i_{\mu_1}^\ast>> T^\ast
\mathbf{J}_1^{-1}(\mu_1) @<\pi_{\mu_1}^\ast <<T^\ast ((T^\ast
Q_1)_{\mu_1})
\end{CD}
\]
$$\mbox{Diagram-5}$$
we have that $\pi_{\mu_1}^\ast\cdot(\varphi_{\mu/G}^\ast)_\ast
\mathbf{d}h_{2\mu_2}=i_{\mu_1}^\ast\cdot(\varphi^\ast)_\ast
\mathbf{d}H_2$, then $$((\varphi_{\mu/G}^\ast)_\ast
\mathbf{d}h_{2\mu_2})^\sharp \cdot
\pi_{\mu_1}=T\pi_{\mu_1}\cdot((\varphi^\ast)_\ast
\mathbf{d}H_2)^\sharp \cdot i_{\mu_1},$$
$$\textnormal{vlift}(\varphi_{\mu/G}^\ast\cdot
f_{2\mu_2}\cdot \varphi_{\mu/G\ast})\cdot \pi_{\mu_1}
=T\pi_{\mu_1}\cdot \textnormal{vlift}(\varphi^\ast
F_2\varphi_\ast)\cdot i_{\mu_1},$$
$$\textnormal{vlift}(\varphi_{\mu/G}^\ast\cdot
u_{2\mu_2}\cdot \varphi_{\mu/G\ast})\cdot \pi_{\mu_1}
=T\pi_{\mu_1}\cdot \textnormal{vlift}(\varphi^\ast
u_2\varphi_\ast)\cdot i_{\mu_1},$$ where
$\varphi_{\mu/G\ast}=(\varphi^{-1})^\ast_{\mu/G}:(T^\ast
Q_1)_{\mu_1}\rightarrow (T^\ast Q_2)_{\mu_2}$ and
$(\varphi_{\mu/G}^\ast)_\ast=(\varphi_{\mu/G\ast})^\ast:
T^\ast((T^\ast Q_2)_{\mu_2})\rightarrow T^\ast((T^\ast
Q_1)_{\mu_1})$. From Hamiltonian matching condition RpHM-3 we have
that
\begin{equation}\begin{aligned}&Im[(\mathrm{d}h_{1\mu_1})^\sharp
+\textnormal{vlift}(f_{1\mu_1})-((\varphi_{\mu/G}^\ast)_\ast
\mathrm{d}h_{2\mu_2})^\sharp-\textnormal{vlift}(\varphi_{\mu/G}^\ast
\cdot f_{2\mu_2}\cdot \varphi_{\mu/G \ast})]\subset
\textnormal{vlift}(W_{1\mu_1}).\end{aligned}\label{4.5}\end{equation}
So,
$$((T^\ast
Q_1)_{\mu_1},\omega_{1\mu_1},h_{1\mu_1},f_{1\mu_1},W_{1\mu_1})\stackrel{RCH}{\sim}((T^\ast
Q_2)_{\mu_2},\omega_{2\mu_2},h_{2\mu_2},f_{2\mu_2},W_{2\mu_2}).$$

Conversely, assume that the $R_P$-reduced RCH systems $((T^\ast
Q_i)_{\mu_i},\omega_{i\mu_i},h_{i\mu_i},f_{i\mu_i},W_{i\mu_i})$,
$i=1,2,$ are RCH-equivalent. Then there exists a diffeomorphism
$\varphi_{\mu/G}^\ast:(T^\ast Q_2)_{\mu_2}\rightarrow (T^\ast
Q_1)_{\mu_1}$, which is symplectic,
$W_{1\mu_1}=\varphi_{\mu/G}^\ast(W_{2\mu_2}),\;
\mu_i\in\mathfrak{g}_i^\ast, \; i=1,2$ and (\ref{4.5}) holds. We can
define a map $\varphi_\mu^\ast:\mathbf{J}^{-1}_2(\mu_2)\rightarrow
\mathbf{J}^{-1}_1(\mu_1)$ such that $\pi_{\mu_1}\cdot
\varphi_\mu^\ast=\varphi_{\mu/G}^\ast\cdot \pi_{\mu_2}, $ and the
map $\varphi^\ast: T^\ast Q_2\rightarrow T^\ast Q_1$ such that
$\varphi^\ast\cdot i_{\mu_2}=i_{\mu_1}\cdot \varphi_{\mu}^\ast, $
see the commutative Diagram-3, as well as a diffeomorphism $\varphi:
Q_1\rightarrow Q_2, $ whose cotangent lift is just
$\varphi^\ast:T^\ast Q_2\rightarrow T^\ast Q_1$. From definition of
$\varphi_\mu^\ast$, we know that $\varphi_\mu^\ast$ is
$(G_{2\mu_2},G_{1\mu_1})$-equivariant. In fact, for any $z_i\in
\mathbf{J}_i^{-1}(\mu_i)$, $g_i\in G_{i\mu_i}$, $i=1,2$ such that
$z_1=\varphi_\mu^\ast(z_2)$, $[z_1]=\varphi^\ast_{\mu/G}[z_2]$, then
we have that
\begin{align*} \pi_{\mu_1}\cdot\varphi_\mu^\ast(\Phi_{2g_2}(z_2))&
=\pi_{\mu_1}\cdot \varphi_\mu^\ast(g_2z_2)
=\varphi_{\mu/G}^\ast\cdot \pi_{\mu_2}(g_2z_2)=\varphi_{\mu/G}^\ast[z_2]=[z_1]\\
&=\pi_{\mu_1}(g_1z_1)=\pi_{\mu_1}(\Phi_{1g_1}(z_1))
=\pi_{\mu_1}\cdot \Phi_{1g_1}\cdot \varphi_\mu^\ast(z_2).
\end{align*}
Since $\pi_{\mu_1}$ is surjective, so, $\varphi_\mu^\ast\cdot
\Phi_{2g_2}=\Phi_{1g_1}\cdot \varphi_\mu^\ast$. Moreover, we have
that
\begin{align*}
\pi_{\mu_1}(W_1\cap \mathbf{J}_1^{-1}(\mu_1)) & =W_{1\mu_1}
=\varphi_{\mu/G}^\ast(W_{2\mu_2})\\ & =\varphi_{\mu/G}^\ast\cdot
\pi_{2\mu_2}(W_2\cap \mathbf{J}_2^{-1}(\mu_2)) =\pi_{\mu_1}\cdot
\varphi_\mu^\ast (W_2\cap \mathbf{J}_2^{-1}(\mu_2)).
\end{align*}
Since $W_i \cap \mathbf{J}_i^{-1}(\mu_i)\neq \emptyset,\; i=1,2, $
and $\pi_{\mu_1}$ is surjective, then $W_1\cap
\mathbf{J}_1^{-1}(\mu_1)=\varphi_\mu^\ast (W_2\cap
\mathbf{J}_2^{-1}(\mu_2))$. We shall show that $\varphi^\ast$ is
symplectic. Because $\varphi_{\mu/G}^\ast:(T^\ast
Q_2)_{\mu_2}\rightarrow (T^\ast Q_1)_{\mu_1}$ is symplectic, the map
$(\varphi_{\mu/G}^\ast)^\ast:\Omega^2((T^\ast
Q_1)_{\mu_1})\rightarrow \Omega^2((T^\ast Q_2)_{\mu_2})$ satisfies
$(\varphi_{\mu/G}^\ast)^\ast \omega_{1\mu_1}=\omega_{2\mu_2}$. By
(\ref{4.1}), $i_{\mu_i}^\ast
\omega_i=\pi_{\mu_i}^\ast\omega_{i\mu_i}$, $i=1,2$, from the
commutative Diagram-4, we have that
\begin{align*}
i_{\mu_2}^\ast\omega_2&=\pi_{\mu_2}^\ast\omega_{2\mu_2}
=\pi_{\mu_2}^\ast\cdot(\varphi_{\mu/G}^\ast)^\ast\omega_{1\mu_1}
=(\varphi_{\mu/G}^\ast\cdot \pi_{\mu_2})^\ast\omega_{1\mu_1}
=(\pi_{\mu_1}\cdot \varphi_\mu^\ast)^\ast\omega_{1\mu_1}\\
&=(i_{\mu_1}^{-1}\cdot \varphi^\ast\cdot i_{\mu_2})^\ast\cdot
\pi_{\mu_1}^\ast \omega_{1\mu_1}
=i_{\mu_2}^\ast\cdot(\varphi^\ast)^\ast\cdot(i^{-1}_{\mu_1})^\ast\cdot
i_{\mu_1}^\ast\omega_1=i_{\mu_2}^\ast\cdot
(\varphi^\ast)^\ast\omega_1.
\end{align*}
Notice that $i_{\mu_2}^\ast$ is injective, thus,
$\omega_2=(\varphi^\ast)^\ast\omega_1$. Since the vector field
$X_{(T^\ast Q_i,G_i,\omega_i,H_i,F_i,u_i)}$ and\\ $X_{((T^\ast
Q_i)_{\mu_i},\omega_{i\mu_i},h_{i\mu_i},f_{i\mu_i},u_{i\mu_i})}$ is
$\pi_{\mu_i}$-related, $i=1,2,$ and $H_i, F_i$ and $W_i$ are all
$G_i$-invariant, $i=1,2$, in the same way, from (\ref{4.5}), we have
that
$$Im[(\mathbf{d}H_1)^\sharp+\textnormal{vlift}(F_1)-((\varphi_\ast)^\ast
\mathbf{d}H_2)^\sharp-\textnormal{vlift}(\varphi^\ast
F_2\varphi_\ast)]\subset \textnormal{vlift}(W_1),$$ that is,
Hamiltonian matching condition RpHM-3 holds. Thus,
$$(T^\ast Q_1, G_1, \omega_1, H_1, F_1,
W_1)\stackrel{RpCH}{\sim}(T^\ast Q_2, G_2, \omega_2, H_2, F_2, W_2).
 \hskip 1cm \blacksquare $$

\begin{remark}
If $(T^\ast Q, \omega)$ is a connected symplectic manifold, and
$\mathbf{J}:T^\ast Q\rightarrow \mathfrak{g}^\ast$ is a
non-equivariant momentum map with a non-equivariance group
one-cocycle $\sigma: G\rightarrow \mathfrak{g}^\ast$, in this case,
we can also define the regular point reducible RCH system
$(T^*Q,G,\omega,H,F,W)$ and RpCH-equivalence, and prove the regular
point reduction theorem for the RCH system by using the above same way,
where the reduced space $((T^\ast Q)_\mu, \omega_\mu )$ is
determined by the affine action given in Remark 4.1.
\end{remark}

\section{Regular Orbit Reduction of RCH Systems }

Let $\mu\in \mathfrak{g}^\ast$ be a regular value of the momentum
map $\mathbf{J}$ and $\mathcal{O}_\mu=G\cdot \mu\subset
\mathfrak{g}^\ast$ be the $G$-orbit of the coadjoint $G$-action
through the point $\mu$. Since $G$ acts freely, properly and
symplectically on $T^\ast Q$, then the quotient space $(T^\ast
Q)_{\mathcal{O}_\mu}= \mathbf{J}^{-1}(\mathcal{O}_\mu)/G$ is a
regular quotient symplectic manifold with the symplectic form
$\omega_{\mathcal{O}_\mu}$ uniquely characterized by the relation
\begin{equation}i_{\mathcal{O}_\mu}^\ast \omega=\pi_{\mathcal{O}_{\mu}}^\ast
\omega_{\mathcal{O}
_\mu}+\mathbf{J}_{\mathcal{O}_\mu}^\ast\omega_{\mathcal{O}_\mu}^+,
\label{5.1}\end{equation} where $\mathbf{J}_{\mathcal{O}_\mu}$ is
the restriction of the momentum map $\mathbf{J}$ to
$\mathbf{J}^{-1}(\mathcal{O}_\mu)$, that is,
$\mathbf{J}_{\mathcal{O}_\mu}=\mathbf{J}\cdot i_{\mathcal{O}_\mu}$
and $\omega_{\mathcal{O}_\mu}^+$ is the $(+)$-symplectic structure
on the orbit $\mathcal{O}_\mu$ given by
\begin{equation}\omega_{\mathcal{O}_\mu}^
+(\nu)(\xi_{\mathfrak{g}^\ast}(\nu),\eta_{\mathfrak{g}^\ast}(\nu))
=<\nu,[\xi,\eta]>,\;\; \forall\;\nu\in\mathcal{O}_\mu, \;
\xi,\eta\in \mathfrak{g}. \label{5.2}\end{equation} The maps
$i_{\mathcal{O}_\mu}:\mathbf{J}^{-1}(\mathcal{O}_\mu)\rightarrow
T^\ast Q$ and
$\pi_{\mathcal{O}_\mu}:\mathbf{J}^{-1}(\mathcal{O}_\mu)\rightarrow
(T^\ast Q)_{\mathcal{O}_\mu}$ are natural injection and the
projection, respectively. The pair $((T^\ast
Q)_{\mathcal{O}_\mu},\omega_{\mathcal{O}_\mu})$ is called the
symplectic orbit reduced space of $(T^\ast Q,\omega)$.

\begin{remark}
If $(T^\ast Q, \omega)$ is a connected symplectic manifold, and
$\mathbf{J}:T^\ast Q\rightarrow \mathfrak{g}^\ast$ is a
non-equivariant momentum map with a non-equivariance group
one-cocycle $\sigma: G\rightarrow \mathfrak{g}^\ast$, which is
defined by $\sigma(g):=\mathbf{J}(g\cdot
z)-\operatorname{Ad}^\ast_{g^{-1}}\mathbf{J}(z)$, where $g\in G$ and
$z\in T^\ast Q$. Then we know that $\sigma$ produces a new affine
action $\Theta: G\times \mathfrak{g}^\ast \rightarrow
\mathfrak{g}^\ast $ defined by
$\Theta(g,\mu):=\operatorname{Ad}^\ast_{g^{-1}}\mu + \sigma(g)$,
where $\mu \in \mathfrak{g}^\ast$, with respect to which the given
momentum map $\mathbf{J}$ is equivariant. Assume that $G$ acts
freely and properly on $T^\ast Q$, and $\mathcal{O}_\mu= G\cdot \mu
\subset \mathfrak{g}^\ast$ denotes the G-orbit of the point $\mu \in
\mathfrak{g}^\ast$ with respect to this affine action $\Theta$, and
$\mu$ is a regular value of $\mathbf{J}$. Then the quotient space
$(T^\ast Q)_{\mathcal{O}_\mu}=\mathbf{J}^{-1}(\mathcal{O}_\mu)/ G $
is also a symplectic manifold with symplectic form
$\omega_{\mathcal{O}_\mu}$ uniquely characterized by $(5.1)$, see
Ortega and Ratiu \cite{orra04}.
\end{remark}

Let $H:T^\ast Q\rightarrow \mathbb{R}$ be a $G$-invariant
Hamiltonian, the flow $F_t$ of the Hamiltonian vector field $X_H$
leaves the connected components of
$\mathbf{J}^{-1}(\mathcal{O}_\mu)$ invariant and commutes with the
$G$-action, then it induces a flow $f_t^{\mathcal{O}_\mu}$ on
$(T^\ast Q)_{\mathcal{O}_\mu}$, defined by
$f_t^{\mathcal{O}_\mu}\cdot
\pi_{\mathcal{O}_\mu}=\pi_{\mathcal{O}_\mu} \cdot F_t\cdot
i_{\mathcal{O}_\mu}$, and the vector field $X_{h_{\mathcal{O}_\mu}}$
generated by the flow $f_t^{\mathcal{O}_\mu}$ on $((T^\ast
Q)_{\mathcal{O}_\mu},\omega_{\mathcal{O}_\mu})$ is Hamiltonian with
the associated regular orbit reduced Hamiltonian function
$h_{\mathcal{O}_\mu}:(T^\ast Q)_{\mathcal{O}_\mu}\rightarrow
\mathbb{R}$ defined by $h_{\mathcal{O}_\mu}\cdot
\pi_{\mathcal{O}_\mu}= H\cdot i_{\mathcal{O}_\mu}$ and the
Hamiltonian vector fields $X_H$ and $X_{h_{\mathcal{O}_\mu}}$ are
$\pi_{\mathcal{O}_\mu}$-related.

When $Q=G$ is a Lie group with Lie algebra $\mathfrak{g}$, and the
$G$-action is the cotangent lift of left translation, then the
associated momentum map $\mathbf{J}_L:T^\ast G\rightarrow
\mathfrak{g}^\ast$ is right invariant. In the same way, the momentum
map $\mathbf{J}_R:T^\ast G\rightarrow \mathfrak{g}^\ast$ for the
cotangent lift of right translation is left invariant. For regular
value $\mu \in \mathfrak{g}^\ast$, $\mathcal{O}_\mu=G\cdot
\mu=\{\operatorname{Ad}^\ast_{g^{-1}}\mu|g\in G\}$ and the
Kostant-Kirilllov-Sourian (KKS) symplectic forms on coadjoint orbit
$\mathcal{O}_\mu (\subset \mathfrak{g}^\ast)$ are given by
$$\omega_{\mathcal{O}_\mu}^-(\nu)(\operatorname{ad}_\xi^\ast(\nu),\operatorname{ad}_\eta^\ast(\nu))
=-<\nu,[\xi,\eta]>, \;\;\forall\; \nu \in \mathcal{O}_\mu, \;
\xi,\eta \in \mathfrak{g}.$$ From Ortega and Ratiu \cite{orra04}, we
know that by using the momentum map $\mathbf{J}_R$ one can induce a
symplectic diffeomorphism from the symplectic point reduced space
$((T^\ast G)_\mu,\omega_\mu)$ to the symplectic orbit space
$(\mathcal{O}_\mu,\omega_{\mathcal{O}_\mu}^-)$. In general case, we
maybe thought that the structure of the symplectic orbit reduced
space $((T^\ast Q)_{\mathcal{O}_\mu},\omega_{\mathcal{O}_\mu})$ is
more complex than that of the symplectic point reduced space
$((T^\ast Q)_\mu,\omega_\mu)$, but, from the regular reduction
diagram, we know that the regular orbit reduced space $((T^\ast
Q)_{\mathcal{O}_\mu},\omega_{\mathcal{O}_\mu})$ is symplectically
diffeomorphic to the regular point reduced space $((T^*Q)_\mu,
\omega_\mu)$, and hence is also symplectically diffeomorphic to a
symplectic fiber bundle. Thus, we can introduce a kind of the
regular orbit reducible RCH systems as follows.

\begin{definition}
(Regular Orbit Reducible RCH System) A 6-tuple $(T^\ast Q, G,
\omega,H,F,W)$, where the Hamiltonian $H: T^\ast Q\rightarrow
\mathbb{R}$, the fiber-preserving map $F: T^\ast Q\rightarrow T^\ast
Q$ and the fiber submanifold $W$ of $T^\ast Q$ are all
$G$-invariant, is called a regular orbit reducible RCH system, if
there exists an orbit $\mathcal{O}_\mu, \; \mu\in\mathfrak{g}^\ast$,
where $\mu$ is a regular value of the momentum map $\mathbf{J}$,
such that the regular orbit reduced system, that is, the 5-tuple
$((T^\ast
Q)_{\mathcal{O}_\mu},\omega_{\mathcal{O}_\mu},h_{\mathcal{O}_\mu},f_{\mathcal{O}_\mu},
W_{\mathcal{O}_\mu})$, where $(T^\ast
Q)_{\mathcal{O}_\mu}=\mathbf{J}^{-1}(\mathcal{O}_\mu)/G$,
$\pi_{\mathcal{O}_\mu}^\ast \omega_{\mathcal{O}_\mu}
=i_{\mathcal{O}_\mu}^\ast\omega-\mathbf{J}_{\mathcal{O}_\mu}^\ast\omega_{\mathcal{O}_\mu}^+$,
$h_{\mathcal{O}_\mu}\cdot \pi_{\mathcal{O}_\mu} =H\cdot
i_{\mathcal{O}_\mu}$, $F(\mathbf{J}^{-1}(\mathcal{O}_\mu))\subset
\mathbf{J}^{-1}(\mathcal{O}_\mu)$, $f_{\mathcal{O}_\mu}\cdot
\pi_{\mathcal{O}_\mu}=\pi_{\mathcal{O}_\mu}\cdot F\cdot
i_{\mathcal{O}_\mu}$, and $W \cap
\mathbf{J}^{-1}(\mathcal{O}_\mu)\neq \emptyset $,
$W_{\mathcal{O}_\mu}=\pi_{\mathcal{O}_\mu}(W \cap
\mathbf{J}^{-1}(\mathcal{O}_\mu))$, is a RCH system, which is simply
written as $R_O$-reduced RCH system. Where $((T^\ast
Q)_{\mathcal{O}_\mu},\omega_{\mathcal{O}_\mu})$ is the $R_O$-reduced
space, the function $h_{\mathcal{O}_\mu}:(T^\ast
Q)_{\mathcal{O}_\mu}\rightarrow \mathbb{R}$ is called the reduced
Hamiltonian, the fiber-preserving map $f_{\mathcal{O}_\mu}:(T^\ast
Q)_{\mathcal{O}_\mu} \rightarrow (T^\ast Q)_{\mathcal{O}_\mu}$ is
called the reduced (external) force map, $W_{\mathcal{O}_\mu}$ is a
fiber submanifold of $(T^\ast Q)_{\mathcal{O}_\mu}$, and is called
the reduced control subset.
\end{definition}

It is worthy of noting that for the regular orbit reducible RCH system
$(T^\ast Q,G,\omega,H,F,W)$, the $G$-invariant external force map
$F: T^*Q \rightarrow T^*Q $ has to satisfy the conditions
$F(\mathbf{J}^{-1}(\mathcal{O}_\mu))\subset
\mathbf{J}^{-1}(\mathcal{O}_\mu), $ and $f_{\mathcal{O}_\mu}\cdot
\pi_{\mathcal{O}_\mu}=\pi_{\mathcal{O}_\mu}\cdot F\cdot
i_{\mathcal{O}_\mu}$, such that we can define the reduced external
force map $f_{\mathcal{O}_\mu}:(T^\ast Q)_{\mathcal{O}_\mu}
\rightarrow (T^\ast Q)_{\mathcal{O}_\mu}. $ The condition $W \cap
\mathbf{J}^{-1}(\mathcal{O}_\mu)\neq \emptyset $ in above definition
makes that the  $G$-invariant control subset $W \cap
\mathbf{J}^{-1}(\mathcal{O}_\mu)$ can be reduced and the reduced
control subset is $W_{\mathcal{O}_\mu}= \pi_{\mathcal{O}_\mu}(W \cap
\mathbf{J}^{-1}(\mathcal{O}_\mu))$. If the control subset cannot be
reduced, we cannot get the $R_O$-reduced RCH system. The study of
RCH system which is not regular orbit reducible is beyond the limits
in this paper, it may be a topic in future study.\\

Denote by $X_{(T^\ast Q,G,\omega,H,F,u)}$ the vector field of a
regular orbit reducible RCH system $(T^\ast Q, G,\omega,\\ H,F,W)$
with a control law $u$. Assume that we have that
\begin{equation}X_{(T^\ast Q,G,\omega,H,F,u)}
=(\mathbf{d}H)^\sharp+\textnormal{vlift}(F)+\textnormal{vlift}(u).\label{5.3}
\end{equation}
Then, for the regular orbit reducible RCH system we can also
introduce the regular orbit reduced controlled Hamiltonian
equivalence (RoCH-equivalence) as follows.

\begin{definition}
(RoCH-equivalence) Suppose that we have two regular orbit reducible
RCH systems $(T^\ast Q_i, G_i, \omega_i, H_i, F_i, W_i)$, $i=1,2$,
we say them to be RoCH-equivalent, or simply,\\ $(T^\ast Q_1, G_1,
\omega_1, H_1, F_1, W_1)\stackrel{RoCH}{\sim}(T^\ast Q_2, G_2,
\omega_2, H_2, F_2, W_2)$, if there exists a diffeomorphism
$\varphi:Q_1\rightarrow Q_2$ such that the following Hamiltonian
matching conditions hold:

\noindent {\bf RoHM-1:} The cotangent lift map $\varphi^\ast:
T^\ast Q_2\rightarrow T^\ast Q_1$ is symplectic.

\noindent {\bf RoHM-2:} For $\mathcal{O}_{\mu_i},\; \mu_i\in
\mathfrak{g}^\ast_i$, the regular reducible orbits of RCH systems
$(T^\ast Q_i, G_i, \omega_i, H_i, F_i, W_i)$, $i=1,2$, the map
$\varphi^\ast_{\mathcal{O}_\mu}=i_{\mathcal{O}_{\mu_1}}^{-1}\cdot\varphi^\ast\cdot
i_{\mathcal{O}_{\mu_2}}:\mathbf{J}_2^{-1}(\mathcal{O}_{\mu_2})\rightarrow
\mathbf{J}_1^{-1}(\mathcal{O}_{\mu_1})$ is $(G_2,G_1)$-equivariant,
$W_1\cap
\mathbf{J}_1^{-1}(\mathcal{O}_{\mu_1})=\varphi_{\mathcal{O}_\mu}^\ast
(W_2\cap \mathbf{J}_2^{-1}(\mathcal{O}_{\mu_2}))$, and
$\mathbf{J}_{2\mathcal{O}_{\mu_2}}^\ast
\omega_{2\mathcal{O}_{\mu_2}}^{+}=(\varphi_{\mathcal{O}_\mu}^\ast)^\ast
\cdot\mathbf{J}_{1\mathcal{O}_{\mu_1}}^\ast\omega_{1\mathcal{O}_{\mu_1}}^{+},$
where $\mu=(\mu_1, \mu_2)$, and denote by
$i_{\mathcal{O}_{\mu_1}}^{-1}(S)$ the pre-image of a subset
$S\subset T^\ast Q_1$ for the map
$i_{\mathcal{O}_{\mu_1}}:\mathbf{J}_1^{-1}(\mathcal{O}_{\mu_1})\rightarrow
T^\ast Q_1$.

\noindent {\bf RoHM-3:}
$Im[(\mathbf{d}H_1)^\sharp+\textnormal{vlift}(F_1)-((\varphi_\ast)^\ast\mathbf{d}H_2)^\sharp
-\textnormal{vlift}(\varphi^\ast F_2\varphi_\ast)]\subset
\textnormal{vlift}(W_1).$
\end{definition}

It is worthy of noting that for the regular orbit reducible RCH
system, the induced equivalent map $\varphi^*$ not only keeps the
symplectic structure and the restriction of the $(+)$-symplectic
structure on the regular orbit to
$\mathbf{J}^{-1}(\mathcal{O}_\mu)$, but also keeps the equivariance
of $G$-action on the regular orbit. If a feedback control law
$u_{\mathcal{O}_\mu}:(T^\ast Q)_{\mathcal{O}_\mu}\rightarrow
W_{\mathcal{O}_\mu}$ is chosen, the $R_O$-reduced RCH system
$((T^\ast Q)_{\mathcal{O}_\mu},
\omega_{\mathcal{O}_\mu},h_{\mathcal{O}_\mu},f_{\mathcal{O}_\mu},u_{\mathcal{O}_\mu})$
is a closed-loop regular dynamic system with a control law
$u_{\mathcal{O}_\mu}$. Assume that its vector field $X_{((T^\ast
Q)_{\mathcal{O}_\mu}, \omega_{\mathcal{O}_\mu},
h_{\mathcal{O}_\mu},f_{\mathcal{O}_\mu},u_{\mathcal{O}_\mu})}$ can
be expressed by
\begin{equation}X_{((T^\ast Q)_{\mathcal{O}_\mu},
\omega_{\mathcal{O}_\mu},h_{\mathcal{O}_\mu},f_{\mathcal{O}_\mu},u_{\mathcal{O}_\mu})}=
(\mathbf{d}h_{\mathcal{O}_\mu})^\sharp+\textnormal{vlift}(f_{\mathcal{O}_\mu})
+\textnormal{vlift}(u_{\mathcal{O}_\mu}), \label{5.4}\end{equation}
where $(\mathbf{d}h_{\mathcal{O}_\mu})^\sharp =
X_{h_{\mathcal{O}_\mu}}$, $\textnormal{vlift}(f_{\mathcal{O}_\mu})=
\textnormal{vlift}(f_{\mathcal{O}_\mu})X_{h_{\mathcal{O}_\mu}}$,
$\textnormal{vlift}(u_{\mathcal{O}_\mu})=
\textnormal{vlift}(u_{\mathcal{O}_\mu})X_{h_{\mathcal{O}_\mu}}$, and
satisfies the condition
\begin{equation}X_{((T^\ast Q)_{\mathcal{O}_\mu},
\omega_{\mathcal{O}_\mu},h_{\mathcal{O}_\mu},f_{\mathcal{O}_\mu},u_{\mathcal{O}_\mu})}\cdot
\pi_{\mathcal{O}_\mu} =T\pi_{\mathcal{O}_\mu} \cdot X_{(T^\ast
Q,G,\omega,H,F,u)}\cdot
i_{\mathcal{O}_\mu}.\label{5.5}\end{equation} Then we can obtain the
following regular orbit reduction theorem for the RCH system, which
explains the relationship between the RoCH-equivalence for the regular
orbit reducible RCH system with symmetry and the RCH-equivalence
for the associated $R_O$-reduced RCH system. This theorem can be
regarded as an extension of regular orbit reduction theorem of a
Hamiltonian system under regular controlled Hamiltonian equivalence
conditions.

\begin{theorem}
If two regular orbit reducible RCH systems $(T^\ast Q_i, G_i,
\omega_i, H_i, F_i,W_i)$, $i=1,2,$ are RoCH-equivalent, then their
associated $R_O$-reduced RCH systems $((T^\ast
Q)_{\mathcal{O}_{\mu_i}}, \omega_{i\mathcal{O}_{\mu_i}},
h_{i\mathcal{O}_{\mu_i}}, f_{i\mathcal{O}_{\mu_i}},
W_{i\mathcal{O}_{\mu_i}})$, $i=1,2,$ must be RCH-equivalent.
Conversely, if the $R_O$-reduced RCH systems $((T^\ast
Q)_{\mathcal{O}_{\mu_i}}, \omega_{i\mathcal{O}_{\mu_i}},
h_{i\mathcal{O}_{\mu_i}},\\ f_{i\mathcal{O}_{\mu_i}},
W_{i\mathcal{O}_{\mu_i}})$, $i=1,2,$ are RCH-equivalent and the
induced map
$\varphi^\ast_{\mathcal{O}_\mu}:\mathbf{J}_2^{-1}(\mathcal{O}_{\mu_2})\rightarrow
\mathbf{J}_1^{-1}(\mathcal{O}_{\mu_1})$, such that
$\mathbf{J}_{2\mathcal{O}_{\mu_2}}^\ast
\omega_{2\mathcal{O}_{\mu_2}}^{+}=(\varphi_{\mathcal{O}_\mu}^\ast)^\ast
\cdot\mathbf{J}_{1\mathcal{O}_{\mu_1}}^\ast\omega_{1\mathcal{O}_{\mu_1}}^{+},$
then the regular orbit reducible RCH systems $(T^\ast Q_i, G_i,
\omega_i,\\ H_i, F_i, W_i)$, $i=1,2,$ are RoCH-equivalent.
\end{theorem}

{\bf Proof:} If $(T^\ast
Q_1,G_1,\omega_1,H_1,F_1,W_1)\stackrel{RoCH}{\sim}(T^\ast
Q_2,G_2,\omega_2,H_2,F_2,W_2)$, then there exists a diffeomorphism
$\varphi:Q_1\rightarrow Q_2$, such that $\varphi^\ast:T^\ast
Q_2\rightarrow T^\ast Q_1$ is symplectic and for $\mu_i\in
\mathfrak{g}_i^\ast,i=1,2$,
$\varphi_{\mathcal{O}_\mu}^\ast=i_{\mathcal{O}_{\mu_1}}^{-1}\cdot
\varphi^\ast\cdot
i_{\mathcal{O}_{\mu_2}}:\mathbf{J}_2^{-1}(\mathcal{O}_{\mu_2})\rightarrow
\mathbf{J}_1^{-1}(\mathcal{O}_{\mu_1})$ is $(G_2,G_1)$-equivariant,
$W_1\cap
\mathbf{J}_1^{-1}(\mathcal{O}_{\mu_1})=\varphi_{\mathcal{O}_\mu}^\ast
(W_2\cap \mathbf{J}_2^{-1}(\mathcal{O}_{\mu_2}))$,
$\mathbf{J}_{2\mathcal{O}_{\mu_2}}^\ast\omega_{2\mathcal{O}_{\mu_2}}^+
=(\varphi_{\mathcal{O}_\mu}^\ast)^\ast
\cdot\mathbf{J}_{1\mathcal{O}_{\mu_1}}^\ast\omega_{1\mathcal{O}_{\mu_1}}^+$,
and RoHM-3 holds. From the following commutative Diagram-6,
\[
\begin{CD}
T^\ast Q_2 @<i_{\mathcal{O}_{\mu_2}}<<
\mathbf{J}_2^{-1}(\mathcal{O}_{\mu_2}) @>\pi_{\mathcal{O}_{\mu_2}}
>> (T^\ast Q_2)_{\mathcal{O}_{\mu_2}}\\
@V\varphi^\ast VV @V\varphi^\ast_{\mathcal{O}_\mu} VV @V\varphi^\ast_{\mathcal{O}_{\mu/G}}VV\\
T^\ast Q_1 @<i_{\mathcal{O}_{\mu_1}}<<
\mathbf{J}_1^{-1}(\mathcal{O}_{\mu_1})
@>\pi_{\mathcal{O}_{\mu_1}}>>(T^\ast Q_1)_{\mathcal{O}_{\mu_1}}
\end{CD}
\]
$$\mbox{Diagram-6}$$
we can define a map $\varphi_{\mathcal{O}_\mu/G}^\ast:(T^\ast
Q_2)_{\mathcal{O}_{\mu_2}}\rightarrow (T^\ast
Q_1)_{\mathcal{O}_{\mu_1}}$, such that
$\varphi_{\mathcal{O}_\mu/G}^\ast \cdot
\pi_{\mathcal{O}_{\mu_2}}=\pi_{\mathcal{O}_{\mu_1}}\cdot
\varphi_{\mathcal{O}_\mu}^\ast$. Because
$\varphi_{\mathcal{O}_\mu}^\ast:
\mathbf{J}_2^{-1}(\mathcal{O}_{\mu_2})\rightarrow
\mathbf{J}_1^{-1}(\mathcal{O}_{\mu_1})$ is $(G_2,G_1)$-equivariant,
$\varphi_{\mathcal{O}_\mu/G}^\ast$ is well-defined. We can prove
that $\varphi_{\mathcal{O}_\mu/G}^\ast$ is symplectic, that is,
$(\varphi_{\mathcal{O}_\mu/G}^\ast)^\ast\omega_{1\mathcal{O}_{\mu_1}}=\omega_{2\mathcal{O}_{\mu_2}}$
and $W_{1\mathcal{O}_{\mu_1}}=\varphi_{\mathcal{O}_\mu/G}^\ast
(W_{2\mathcal{O}_{\mu_2}})$. In fact, since $\varphi^\ast:T^\ast
Q_1\to T^\ast Q_2$ is symplectic, the map
$(\varphi^\ast)^\ast:\Omega^2(T^\ast Q_1)\rightarrow \Omega^2(T^\ast
Q_2)$ satisfies $(\varphi^\ast)^\ast\omega_1=\omega_2$. By
(\ref{5.1}), $i_{\mathcal{O}_{\mu_i}}^\ast
\omega_i=\pi_{\mathcal{O}_{\mu_i}}^\ast\omega_{i\mathcal{O}_{\mu_i}}
+\mathbf{J}_{i\mathcal{O}_{\mu_i}}^\ast\omega_{i\mathcal{O}_{\mu_i}}^{+}$,
$i=1,2,$ and
$\mathbf{J}_{2\mathcal{O}_{\mu_2}}^\ast\omega_{2\mathcal{O}_{\mu_2}}^+=(\varphi_{\mathcal{O}_\mu}^\ast)^\ast
\cdot\mathbf{J}_{1\mathcal{O}_{\mu_1}}^\ast\omega_{1\mathcal{O}_{\mu_1}}^+$,
from the following commutative Diagram-7,
\[
\begin{CD}
\Omega^2(T^\ast Q_1) @>i_{\mathcal{O}_{\mu_1}}^\ast>>
\Omega^2(\mathbf{J}_1^{-1}(\mathcal{O}_{\mu_1}))
@<\pi_{\mathcal{O}_{\mu_1}}^\ast<< \Omega^2((T^\ast Q_1)_{\mathcal{O}_{\mu_1}})\\
@V(\varphi^\ast)^\ast VV @V(\varphi^\ast_{\mathcal{O}_\mu})^\ast
VV @V(\varphi^\ast_{\mathcal{O}_\mu/G})^\ast VV\\
\Omega^2(T^\ast Q_2) @>i_{\mathcal{O}_{\mu_2}}^\ast>>
\Omega^2(\mathbf{J}_2^{-1}(\mathcal{O}_{\mu_2}))
@<\pi_{\mathcal{O}_{\mu_2}}^\ast <<\Omega^2((T^\ast
Q_2)_{\mathcal{O}_{\mu_2}})
\end{CD}\]
$$\mbox{Diagram-7}$$
we have that
\begin{align*}
&\pi_{\mathcal{O}_{\mu_2}}^\ast\cdot
(\varphi_{\mathcal{O}_\mu/G}^\ast)^\ast
\omega_{1\mathcal{O}_{\mu_1}}
=(\varphi_{\mathcal{O}_\mu/G}^\ast\cdot
\pi_{\mathcal{O}_{\mu_2}})^\ast\omega_{1\mathcal{O}_{\mu_1}}
=(\pi_{\mathcal{O}_{\mu_1}}\cdot \varphi_{\mathcal{O}_\mu}^\ast)^\ast\omega_{1\mathcal{O}_{\mu_1}}\\
&=(\varphi_{\mathcal{O}_\mu}^\ast)^\ast\cdot
\pi_{\mathcal{O}_{\mu_1}}^\ast\omega_{1\mathcal{O}_{\mu_1}}
=(i_{\mathcal{O}_{\mu_1}}^{-1}\cdot\varphi^\ast\cdot
i_{\mathcal{O}_{\mu_2}})^\ast\cdot
i_{\mathcal{O}_{\mu_1}}^\ast\omega_1-(\varphi_{\mathcal{O}_\mu}^\ast)^\ast
\cdot\mathbf{J}_{1\mathcal{O}_{\mu_1}}^\ast \omega_{1\mathcal{O}_{\mu_1}}^+\\
&=i_{\mathcal{O}_{\mu_2}}^\ast\cdot(\varphi^\ast)^\ast\omega_1
-\mathbf{J}_{2\mathcal{O}_{\mu_2}}^\ast\omega_{2\mathcal{O}_{\mu_2}}^+
=i_{\mathcal{O}_{\mu_2}}^\ast\omega_2-\mathbf{J}_{2\mathcal{O}_{\mu_2}}^\ast\omega_{2\mathcal{O}_{\mu_2}}^+
=\pi_{\mathcal{O}_{\mu_2}}^\ast\omega_{2\mathcal{O}_{\mu_2}}.
\end{align*}
Because $\pi_{\mathcal{O}_{\mu_2}}^\ast$ is surjective, thus
$(\varphi_{\mathcal{O}_\mu/G}^\ast)^\ast\omega_{1\mathcal{O}_{\mu_1}}=\omega_{2\mathcal{O}_{\mu_2}}$.
Notice that $W_i\cap \mathbf{J}_i^{-1}(\mathcal{O}_{\mu_i})\neq
\emptyset $,
$W_{i\mathcal{O}_{\mu_i}}=\pi_{\mathcal{O}_{\mu_i}}(W_i\cap
\mathbf{J}_i^{-1}(\mathcal{O}_{\mu_i}))$, $i=1,2,$ and $W_1\cap
\mathbf{J}_1^{-1}(\mathcal{O}_{\mu_1})=\varphi_{\mathcal{O}_\mu}^\ast
(W_2\cap \mathbf{J}_2^{-1}(\mathcal{O}_{\mu_2}))$, we have that
\begin{align*}
W_{1\mathcal{O}_{\mu_1}}& =\pi_{\mathcal{O}_{\mu_1}}(W_1\cap
\mathbf{J}_1^{-1}(\mathcal{O}_{\mu_1}))
=\pi_{\mathcal{O}_{\mu_1}}\cdot \varphi_{\mathcal{O}_\mu}^\ast
(W_2\cap \mathbf{J}_2^{-1}(\mathcal{O}_{\mu_2}))\\
& =\varphi_{\mathcal{O}_\mu/G}^\ast\cdot
\pi_{\mathcal{O}_{\mu_2}}(W_2\cap
\mathbf{J}_2^{-1}(\mathcal{O}_{\mu_2}))
=\varphi_{\mathcal{O}_\mu/G}^\ast (W_{2\mathcal{O}_{\mu_2}}).
\end{align*} Next, from (\ref{5.3}) and
(\ref{5.4}), we know that for $i=1,2,$
$$X_{(T^\ast Q_i, G_i, \omega_i, H_i, F_i, u_i)}
=(\mathbf{d}H_i)^\sharp+\textnormal{vlift}(F_i)+\textnormal{vlift}(u_i),$$
$$X_{((T^\ast Q_i)_{\mathcal{O}_{\mu_i}},
\omega_{i\mathcal{O}_{\mu_i}},h_{i\mathcal{O}_{\mu_i}},f_{i\mathcal{O}_{\mu_i}},u_{i\mathcal{O}_{\mu_i}})}=
(\mathbf{d}h_{i\mathcal{O}_{\mu_i}})^\sharp+\textnormal{vlift}(f_{i\mathcal{O}_{\mu_i}})
+\textnormal{vlift}(u_{i\mathcal{O}_{\mu_i}}),
$$ and from (\ref{5.5}), we have that
$$X_{((T^\ast
Q_i)_{\mathcal{O}_{\mu_i}},\omega_{i\mathcal{O}_{\mu_i}},h_{i\mathcal{O}_{\mu_i}},
f_{i\mathcal{O}_{\mu_i}},u_{i\mathcal{O}_{\mu_i}})} \cdot
\pi_{\mathcal{O}_{\mu_i}} =T\pi_{\mathcal{O}_{\mu_i}} \cdot
X_{(T^\ast Q_i,G_i,\omega_i,H_i,F_i,u_i)} \cdot
i_{\mathcal{O}_{\mu_i}}.$$
Since $H_i$, $F_i$ and $W_i$ are all
$G_i$-invariant, $i=1,2,$ and
\begin{align*}h_{i\mathcal{O}_{\mu_i}}\cdot
\pi_{\mathcal{O}_{\mu_i}}&=H_i\cdot i_{\mathcal{O}_{\mu_i}},\;\;
f_{i\mathcal{O}_{\mu_i}}\cdot
\pi_{\mathcal{O}_{\mu_i}}=\pi_{\mathcal{O}_{\mu_i}}\cdot F_i\cdot
i_{\mathcal{O}_{\mu_i}},\;\; u_{i\mathcal{O}_{\mu_i}}\cdot
\pi_{\mathcal{O}_{\mu_i}}&=\pi_{\mathcal{O}_{\mu_i}}\cdot u_i\cdot
i_{\mathcal{O}_{\mu_i}},\;\; i=1,2.\end{align*} From the following
commutative Diagram-8,
\[
\begin{CD}
T^\ast T^\ast Q_2 @>i_{\mathcal{O}_{\mu_2}}^\ast>> T^\ast
\mathbf{J}_2^{-1}(\mathcal{O}_{\mu_2})
@<\pi_{\mathcal{O}_{\mu_2}}^\ast<< T^\ast ((T^\ast Q_2)_{\mathcal{O}_{\mu_2}})\\
@V(\varphi^\ast)_\ast VV @V(\varphi^\ast_{\mathcal{O}_\mu})_\ast VV
@V(\varphi^\ast_{\mathcal{O}_\mu/G})_\ast VV\\
T^\ast T^\ast Q_1 @>i_{\mathcal{O}_{\mu_1}}^\ast>> T^\ast
\mathbf{J}_1^{-1}(\mathcal{O}_{\mu_1})
@<\pi_{\mathcal{O}_{\mu_1}}^\ast <<T^\ast ((T^\ast
Q_1)_{\mathcal{O}_{\mu_1}})
\end{CD}\]
$$\mbox{Diagram-8}$$
we have that
$\pi_{\mathcal{O}_{\mu_1}}^\ast\cdot(\varphi_{\mathcal{O}_\mu/G}^\ast)_\ast
\mathbf{d}h_{2\mathcal{O}_{\mu_2}}=i_{\mathcal{O}_{\mu_1}}^\ast
\cdot(\varphi^\ast)_\ast \mathbf{d}H_2$, then
$$((\varphi_{\mathcal{O}_\mu/G}^\ast)_\ast
\mathbf{d}h_{2\mathcal{O}_{\mu_2}})^\sharp\cdot
\pi_{\mathcal{O}_{\mu_1}}=T\pi_{\mathcal{O}_{\mu_1}}\cdot((\varphi^\ast)_\ast
\mathbf{d}H_2)^\sharp \cdot i_{\mathcal{O}_{\mu_1}},$$
$$\textnormal{vlift}(\varphi_{\mathcal{O}_\mu/G}^\ast\cdot
f_{2\mathcal{O}_{\mu_2}}\cdot \varphi_{\mathcal{O}_\mu/G\ast} )\cdot
\pi_{\mathcal{O}_{\mu_1}}=T\pi_{\mathcal{O}_{\mu_1}}\cdot
\textnormal{vlift}(\varphi^\ast F_2\varphi_\ast)\cdot
i_{\mathcal{O}_{\mu_1}},$$
$$\textnormal{vlift}(\varphi_{\mathcal{O}_\mu/G}^\ast\cdot
u_{2\mathcal{O}_{\mu_2}}\cdot \varphi_{\mathcal{O}_\mu/G\ast})\cdot
\pi_{\mathcal{O}_{\mu_1}} =T\pi_{\mathcal{O}_{\mu_1}}\cdot
\textnormal{vlift}(\varphi^\ast u_2\varphi_\ast)\cdot
i_{\mathcal{O}_{\mu_1}},$$ where the map
$\varphi_{\mathcal{O}_\mu/G\ast}=(\varphi^{-1})_{\mathcal{O}_\mu/G}^\ast:
(T^\ast Q_1)_{\mathcal{O}_{\mu_1}}\rightarrow (T^\ast
Q_2)_{\mathcal{O}_{\mu_2}}$ and
$(\varphi_{\mathcal{O}_\mu/G}^\ast)_\ast=(\varphi_{\mathcal{O}_\mu/G\ast})^\ast:
T^\ast((T^\ast Q_2)_{\mathcal{O}_{\mu_2}})\rightarrow T^\ast
((T^\ast Q_1)_{\mathcal{O}_{\mu_1}})$. From the Hamiltonian matching
condition RoHM-3 we have that
\begin{equation}Im[(\mathbf{d}h_{1\mathcal{O}_{\mu_1}})^\sharp+\textnormal{vlift}(f_{1\mathcal{O}_{\mu_1}})
-((\varphi_{\mathcal{O}_\mu/G}^\ast)_\ast
\mathbf{d}h_{2\mathcal{O}_{\mu_2}})^\sharp-\textnormal{vlift}(\varphi_{\mathcal{O}_\mu/G}^\ast\cdot
f_{2\mathcal{O}_{\mu_2}}\cdot
\varphi_{\mathcal{O}_\mu/G\ast})]\subset
\textnormal{vlift}(W_{1\mathcal{O}_{\mu_1}}). \label{5.6}
\end{equation}
So,
\begin{align*}((T^\ast Q_1)_{\mathcal{O}_{\mu_1}},
\omega_{1\mathcal{O}_{\mu_1}}, h_{1\mathcal{O}_{\mu_1}},
f_{1\mathcal{O}_{\mu_1}},
W_{1\mathcal{O}_{\mu_1}})\stackrel{RCH}{\sim}((T^\ast
Q_2)_{\mathcal{O}_{\mu_2}}, \omega_{2\mathcal{O}_{\mu_2}},
h_{2\mathcal{O}_{\mu_2}}, f_{2\mathcal{O}_{\mu_2}},
W_{2\mathcal{O}_{\mu_2}}).\end{align*}

Conversely, assume that the $R_O$-reduced RCH systems $((T^\ast
Q_i)_{\mathcal{O}_{\mu_i}}, \omega_{i\mathcal{O}_{\mu_i}},
h_{i\mathcal{O}_{\mu_i}}, f_{i\mathcal{O}_{\mu_i}},
W_{i\mathcal{O}_{\mu_i}})$, $i=1,2,$ are RCH-equivalent, then there
exists a diffeomorphism $\varphi_{\mathcal{O}_\mu/G}^\ast:(T^\ast
Q_2)_{\mathcal{O}_{\mu_2}}\rightarrow (T^\ast
Q_1)_{\mathcal{O}_{\mu_1}}$, which is symplectic, $
W_{1\mathcal{O}_{\mu_1}}=\varphi_{\mathcal{O}_\mu/G}^\ast(W_{2\mathcal{O}_\mu})$
and (\ref{5.6}) holds. Thus, we can define a map
$\varphi_{\mathcal{O}_\mu}^\ast:\mathbf{J}_2^{-1}(\mathcal{O}_{\mu_2})\rightarrow
\mathbf{J}_1^{-1}(\mathcal{O}_{\mu_1})$ such that
$\pi_{\mathcal{O}_{\mu_1}}\cdot
\varphi_{\mathcal{O}_\mu}^\ast=\varphi_{\mathcal{O}_\mu/G}^\ast
\cdot \pi_{\mathcal{O}_{\mu_2}}, $ and map $\varphi^\ast: T^\ast
Q_2\rightarrow T^\ast Q_1$ such that $i_{\mathcal{O}_{\mu_1}}\cdot
\varphi_{\mathcal{O}_\mu}^\ast=\varphi^\ast\cdot
i_{\mathcal{O}_{\mu_2}}, $ see the commutative Diagram-6, as well as
a diffeomorphism $\varphi: Q_1\rightarrow Q_2$, whose cotangent lift
is just $\varphi^\ast:T^\ast Q_2\rightarrow T^\ast Q_1$. At first,
from definition of $\varphi_{\mathcal{O}_\mu}^\ast$ we know that
$\varphi_{\mathcal{O}_\mu}^\ast$ is $(G_2,G_1)$-equivariant. In
fact, for any $z_i\in \mathbf{J}_i^{-1}(\mathcal{O}_{\mu_i})$,
$g_i\in G_{i}$, $i=1,2$ such that
$z_1=\varphi_{\mathcal{O}_\mu}^\ast(z_2)$,
$[z_1]=\varphi^\ast_{\mathcal{O}_\mu/G}[z_2]$, then we have that
\begin{align*}
&\pi_{\mathcal{O}_{\mu_1}}\cdot
\varphi_{\mathcal{O}_\mu}^\ast(\Phi_{2g_2}(z_2))
=\pi_{\mathcal{O}_{\mu_1}}\cdot
\varphi_{\mathcal{O}_\mu}^\ast(g_2z_2)
=\varphi_{\mathcal{O}_\mu/G}^\ast\cdot
\pi_{\mathcal{O}_{\mu_2}}(g_2z_2)
=\varphi_{\mathcal{O}_\mu/G}^\ast[z_2]\\
&=[z_1]=\pi_{\mathcal{O}_{\mu_1}}(g_1z_1)=\pi_{\mathcal{O}_{\mu_1}}(\Phi_{1g_1}(z_1))
=\pi_{\mathcal{O}_{\mu_1}}\cdot \Phi_{1g_1}\cdot
\varphi_{\mathcal{O}_\mu}^\ast(z_2).
\end{align*}
Since $\pi_{\mathcal{O}_{\mu_1}}$ is surjective, so,
$\varphi_{\mathcal{O}_\mu}^\ast\cdot \Phi_{2g_2}=\Phi_{1g_1}\cdot
\varphi_{\mathcal{O}_\mu}^\ast$. Moreover, we have that
\begin{align*}
\pi_{\mathcal{O}_{\mu_1}}(W_1\cap
\mathbf{J}_1^{-1}(\mathcal{O}_{\mu_1}))& =W_{1\mathcal{O}_{\mu_1}}
=\varphi_{\mathcal{O}_\mu/G}^\ast(W_{2\mathcal{O}_{\mu_2}})\\
& =\varphi_{\mathcal{O}_\mu/G}^\ast \cdot
\pi_{\mathcal{O}_{\mu_2}}(W_2\cap
\mathbf{J}_2^{-1}(\mathcal{O}_{\mu_2}))=\pi_{\mathcal{O}_{\mu_1}}\cdot
\varphi_{\mathcal{O}_\mu}^\ast(W_2\cap
\mathbf{J}_2^{-1}(\mathcal{O}_{\mu_2})).
\end{align*} Since $W_i\cap
\mathbf{J}_i^{-1}(\mathcal{O}_{\mu_i})\neq \emptyset,\; i=1,2,$ and
$\pi_{\mathcal{O}_{\mu_1}}$ is surjective, then $W_1\cap
\mathbf{J}_1^{-1}(\mathcal{O}_{\mu_1})=\varphi_{\mathcal{O}_\mu}^\ast
(W_2\cap \mathbf{J}_2^{-1}(\mathcal{O}_{\mu_2}))$. Now we shall show
that $\varphi^\ast$ is symplectic, that is,
$\omega_2=(\varphi^\ast)^\ast\omega_1$. In fact, since
$\varphi_{\mathcal{O}_\mu/G}^\ast:(T^\ast
Q_2)_{\mathcal{O}_{\mu_2}}\rightarrow (T^\ast
Q_1)_{\mathcal{O}_{\mu_1}}$ is symplectic, the map
$(\varphi_{\mathcal{O}_\mu/G}^\ast)^\ast: \Omega^2((T^\ast
Q_1)_{\mathcal{O}_{\mu_1}})\rightarrow \Omega^2((T^\ast
Q_2)_{\mathcal{O}_{\mu_2}})$ satisfies
$(\varphi_{\mathcal{O}_\mu/G}^\ast)^\ast
\omega_{1\mathcal{O}_{\mu_1}}=\omega_{2\mathcal{O}_{\mu_2}}$. By
(\ref{5.1}),
$i_{\mathcal{O}_{\mu_i}}^\ast\omega_i=\pi_{\mathcal{O}_{\mu_i}}^\ast
\omega_{i\mathcal{O}_{\mu_i}}+\mathbf{J}_{i\mathcal{O}_{\mu_i}}^\ast\omega_{i\mathcal{O}_{\mu_i}}^+$,
$i=1,2$, from the commutative Diagram-7, we have that
\begin{align*}
&i_{\mathcal{O}_{\mu_2}}^\ast\omega_2=\pi_{\mathcal{O}_{\mu_2}}^\ast
\omega_{2\mathcal{O}_{\mu_2}}+\mathbf{J}_{2\mathcal{O}_{\mu_2}}^\ast\omega_{2\mathcal{O}_{\mu_2}}^+
=\pi_{2\mathcal{O}_{\mu_2}}^\ast\cdot
(\varphi_{\mathcal{O}_\mu/G}^\ast)^\ast\omega_{1\mathcal{O}_{\mu_1}}
+\mathbf{J}_{2\mathcal{O}_{\mu_2}}^\ast\omega_{2\mathcal{O}_{\mu_2}}^+\\
&=(\varphi_{\mathcal{O}_\mu/G}^\ast\cdot
\pi_{\mathcal{O}_{\mu_2}})^\ast
\omega_{1\mathcal{O}_{\mu_1}}+\mathbf{J}_{2\mathcal{O}_{\mu_2}}^\ast\omega_{2\mathcal{O}_{\mu_2}}^{+}
=(\pi_{\mathcal{O}_{\mu_1}}\cdot
\varphi_{\mathcal{O}_\mu}^\ast)^\ast
\omega_{1\mathcal{O}_{\mu_1}}+\mathbf{J}_{2\mathcal{O}_{\mu_2}}^\ast\omega_{2\mathcal{O}_{\mu_2}}^{+}\\
&=(i_{\mathcal{O}_{\mu_1}}^{-1}\cdot\varphi^\ast\cdot
i_{\mathcal{O}_{\mu_2}})^\ast \cdot
\pi_{\mathcal{O}_{\mu_1}}^\ast\omega_{1\mathcal{O}_{\mu_1}}
+\mathbf{J}_{2\mathcal{O}_{\mu_2}}^\ast\omega_{2\mathcal{O}_{\mu_2}}^+\\
&=i_{\mathcal{O}_{\mu_2}}^\ast\cdot(\varphi^\ast)^\ast\cdot
(i_{\mathcal{O}_{\mu_1}}^{-1})^\ast \cdot
[i_{\mathcal{O}_{\mu_1}}^{\ast}\omega_1
-\mathbf{J}_{1\mathcal{O}_{\mu_1}}^{\ast}\omega_{1\mathcal{O}_{\mu_1}}^+]
+\mathbf{J}_{2\mathcal{O}_{\mu_2}}^{\ast}\omega_{2\mathcal{O}_{\mu_2}}^{+}\\
&=i_{\mathcal{O}_{\mu_2}}^{\ast}\cdot(\varphi^\ast)^\ast\omega_1-
(\varphi_{\mathcal{O}_\mu}^\ast)^\ast \cdot
\mathbf{J}_{1\mathcal{O}_{\mu_1}}^{\ast}\omega_{1\mathcal{O}_{\mu_1}}^+
+\mathbf{J}_{2\mathcal{O}_{\mu_2}}^{\ast}\omega_{2\mathcal{O}_{\mu_2}}^+
\end{align*}
Notice that $i_{\mathcal{O}_{\mu_2}}^{\ast}$ is injective, and by
our hypothesis,
$\mathbf{J}_{2\mathcal{O}_{\mu_2}}^{\ast}\omega_{2\mathcal{O}_{\mu_2}}^+
=(\varphi_{\mathcal{O}_\mu}^\ast)^\ast \cdot
\mathbf{J}_{1\mathcal{O}_{\mu_1}}^{\ast}\omega_{1\mathcal{O}_{\mu_1}}^+$,
then $\omega_2=(\varphi^\ast)^\ast \omega_1$, that is,
$\varphi^\ast$ is symplectic. Since the vector fields $X_{(T^\ast
Q_i,G_i,\omega_i,H_i,F_i,u_i)}$ and\\ $X_{((T^\ast
Q_i)_{\mathcal{O}_{\mu_i}}, \omega_{i\mathcal{O}_{\mu_i}},
h_{i\mathcal{O}_{\mu_i}},f_{i\mathcal{O}_{\mu_i}},u_{i\mathcal{O}_{\mu_i}})}$
is $\pi_{\mathcal{O}_{\mu_i}}$-related, $i=1,2,$ and $H_i, F_i$ and
$W_i$ are all $G_i$-invariant, $i=1,2$, in the same way, from
(\ref{5.6}) we have that Hamiltonian matching condition RoHM-3
holds. Thus,
$$(T^\ast Q_1,G_1,\omega_1,H_1,F_1,W_1)\stackrel{RoCH}{\sim}(T^\ast
Q_2,G_2,\omega_2,H_2,F_2,W_2).   \hskip 1cm \blacksquare $$

\begin{remark}
If $(T^\ast Q, \omega)$ is a connected symplectic manifold, and
$\mathbf{J}:T^\ast Q\rightarrow \mathfrak{g}^\ast$ is a
non-equivariant momentum map with a non-equivariance group
one-cocycle $\sigma: G\rightarrow \mathfrak{g}^\ast$, in this case,
we can also define the regular orbit reducible RCH system
$(T^*Q,G,\omega,H,F,W)$ and RoCH-equivalence, and prove the regular
orbit reduction theorem for the RCH system by using the above same way,
where the reduced space $((T^\ast Q)_{\mathcal{O}_\mu},
\omega_{\mathcal{O}_\mu} )$ is determined by the affine action given
in Remark 5.1.
\end{remark}

\section{Applications }

As the applications of regular point reduction theory of RCH systems
with symmetries, in this section, we first study the regular point
reducible RCH systems on a Lie group and on its generalization,
respectively, and give their $R_P$-reduced RCH systems, which are
the RCH systems on a coadjoint orbit and on its generalization,
respectively. Next, we describe uniformly the rigid body and heavy
top, as well as them with internal rotors (or the external force
torques) as the regular point reducible RCH systems on the rotation
group $\textmd{SO}(3)$ and on the Euclidean group $\textmd{SE}(3)$,
as well as on their generalizations, respectively, and give their
$R_P$-reduced RCH systems and discuss their RCH-equivalences.
Moreover, in order to understand well the abstract definition of RCH
system and the significance of Theorem 3.3, we describe the RCH
system from the viewpoint of port Hamiltonian system with a
symplectic structure, and state the relationship between
RCH-equivalence and equivalence of port Hamiltonian system.

\subsection{Regular Point Reducible RCH Systems on
a Lie Group and Its Generalization}

Let $G$ be a Lie group with Lie algebra $\mathfrak{g}$ and $T^\ast
G$ its cotangent bundle with the canonical symplectic form
$\omega_0$. A RCH system on $G$ is a 5-tuple $(T^\ast G,\omega_0,
H,F,W)$, where $(T^\ast G,\omega_0, H)$ is a Hamiltonian system and
$H: T^\ast G \rightarrow \mathbb{R}$ is a Hamiltonian, the
fiber-preserving map $F: T^\ast G \rightarrow T^\ast G$ is an
(external) force map and the fiber submanifold $W$ of $T^\ast G$ is
a control subset. In the following we shall give its $R_P$-reduced
RCH system. We know that the left and right translation on $G$
induce the left and right action of $G$ on itself. If $I_g: G \to
G$; $I_g(h)=ghg^{-1}=L_g\cdot R_{g^{-1}}(h)$, for $g,h\in G$, is the
inner automorphism on $G$, then the adjoint representation of a Lie
group $G$ is defined by $\operatorname{Ad}_g=T_eI_g = T_{g^{-1}}L_g
\cdot T_e R_{g^{-1}}:\mathfrak{g}\to \mathfrak{g}$, and the
coadjoint representation is given by
$\operatorname{Ad}_{g^{-1}}^\ast:\mathfrak{g}^\ast\to
\mathfrak{g}^\ast$; $\langle
\operatorname{Ad}_{g^{-1}}^\ast(\mu),\xi\rangle=\langle\mu,\operatorname{Ad}_{g^{-1}}(\xi)\rangle$,
where $\mu\in\mathfrak{g}^\ast$, $\xi\in\mathfrak{g}$ and
$\langle,\rangle$ denotes the pairing between $\mathfrak{g}^\ast$
and $\mathfrak{g}$. Since the coadjoint representation
$\operatorname{Ad}_{g^{-1}}^\ast:\mathfrak{g}^\ast\to
\mathfrak{g}^\ast$ can induce a left coadjoint action of $G$ on
$\mathfrak{g}^\ast$, then the coadjoint orbit $\mathcal{O}_\mu$ of
this action through $\mu\in\mathfrak{g}^\ast$ is an immersed
submanifold of $\mathfrak{g}^\ast$. We know that $\mathfrak{g}^\ast$
is a Poisson manifold with respect to the $(\pm)$-Lie-Poisson
bracket $\{\cdot,\cdot\}_\pm$ defined by
\begin{equation}
\{f,g\}_\pm(\mu):=\pm<\mu,[\frac{\delta f}{\delta \mu},\frac{\delta
g}{\delta\mu}]>,\;\; \forall f,g\in C^\infty(\mathfrak{g}^\ast),\;\;
\mu\in \mathfrak{g}^\ast,\label{6.1}
\end{equation}
where the element $\frac{\delta f}{\delta \mu}\in\mathfrak{g}$ is
defined by the equality $<v,\frac{\delta f}{\delta
\mu}>:=Df(\mu)\cdot v$, for any $v\in \mathfrak{g}^\ast$, see
Marsden and Ratiu \cite{mara99}. Thus, for the coadjoint orbit
$\mathcal{O}_\mu, \; \mu\in\mathfrak{g}^\ast$, the orbit symplectic
structure can be defined by
\begin{equation}\omega_{\mathcal{O}_\mu}^\pm(\nu)(\operatorname{ad}_\xi^\ast(\nu),
\operatorname{ad}_\eta^\ast(\nu))=\pm \langle\nu,[\xi,\eta]\rangle,
\qquad \forall\; \xi,\eta\in\mathfrak{g}, \;\;
\nu\in\mathcal{O}_\mu\subset\mathfrak{g}^\ast,
\label{6.2}\end{equation} which coincide with the restriction of
the Lie-Poisson brackets on $\mathfrak{g}^\ast$ to the coadjoint
orbit $\mathcal{O}_\mu$. From the Symplectic Stratification theorem
we know that a finite dimensional Poisson manifold is the disjoint
union of its symplectic leaves, and its each symplectic leaf is an
injective immersed Poisson submanifold whose induced Poisson
structure is symplectic. In consequence, when $\mathfrak{g}^\ast$ is
endowed one of the Lie-Poisson structures $\{\cdot,\cdot\}_\pm$, the
symplectic leaves of the Poisson manifolds
$(\mathfrak{g}^\ast,\{\cdot,\cdot\}_\pm)$ coincide with the
connected components of the orbits of the elements in
$\mathfrak{g}^\ast$ under the coadjoint action.\\

We now identify $T^\ast G$ and
$G\times\mathfrak{g}^\ast$ by using the left translation. In fact,
the map $\lambda: T^\ast G \rightarrow G \times \mathfrak{g}^\ast,
\; \lambda(\alpha_g):=(g,(T_eL_g)^\ast \alpha_g)$, for any
$\alpha_g \in T^\ast_g G$, which defines a vector bundle
isomorphism usually referred to as the local left trivialization of
$T^\ast G$. If
the left $G$-action $L_g: G \rightarrow G $ is free and proper,
then the cotangent lift of the action to its cotangent bundle
$T^\ast G$, given by $\Phi^{T^*}: G\times T^*G \rightarrow T^*G,
\; g\cdot(h,\nu):=(gh,\nu)$, for any $g,h \in G,\; \nu \in
\mathfrak{g}^\ast$,  is also a free and proper action, and the
orbit space $(T^\ast G)/ G$ is a smooth manifold and $\pi:T^*G
\rightarrow (T^*G )/G $ is a smooth submersion. We note that
$(T^\ast G)/ G$ is diffeomorphic to $(G\times\mathfrak{g}^\ast)/
G$, since $G$ acts trivially on $\mathfrak{g}^\ast$, it follows that
$(G\times\mathfrak{g}^\ast)/ G \cong \mathfrak{g}^\ast, $
and $(T^\ast G)/ G$ is diffeomorphic to $\mathfrak{g}^\ast$.
Moreover, we consider $T^\ast G$ as a symplectic manifold with
the canonical symplectic form $\omega_0$. Assume that the cotangent lifted
left action $\Phi^{T^*}$ is symplectic, and admits an
associated $\operatorname{Ad}^\ast$-equivariant momentum map
$\mathbf{J}_L:T^\ast G \rightarrow \mathfrak{g}^\ast$,
$\mathbf{J}_L(g,\alpha_g)=T_e^\ast R_g(\alpha_g)$, $\forall\; g
\in G, \; \alpha_g \in T^\ast_g G$. If
$\mu\in\mathfrak{g}^\ast$ is a regular value of $\mathbf{J}_L$ and
denote $G_\mu$ the isotropy subgroup of the coadjoint action
$G$ at the point $\mu\in\mathfrak{g}^\ast$, which is defined by
$G_\mu=\{g\in G|\operatorname{Ad}_g^\ast \mu=\mu \}$. It follows
that $G_\mu$ acts also freely and properly on
$\mathbf{J}_L^{-1}(\mu)$, the regular point reduced space $(T^\ast
G)_\mu=\mathbf{J}_L^{-1}(\mu)/G_\mu$ of $(T^\ast G,\omega_0)$ at
$\mu$, is a symplectic manifold with symplectic form $\omega_\mu$
uniquely characterized by the relation $\pi_\mu^\ast
\omega_\mu=i_\mu^\ast \omega_0$, where the map
$i_\mu:\mathbf{J}_L^{-1}(\mu)\rightarrow T^\ast G$ is the
inclusion and $\pi_\mu:\mathbf{J}_L^{-1}(\mu)\rightarrow (T^\ast
G)_\mu$ is the projection. From Abraham and
Marsden \cite{abma78}, we have the following result.

\begin{proposition}
The coadjoint orbit $(\mathcal{O}_\mu,
\omega_{\mathcal{O}_\mu}^{-}), \; \mu\in \mathfrak{g}^\ast,$ is
symplectically diffeomorphic to a regular point reduced space
$((T^\ast G)_\mu,\omega_\mu)$ of $T^*G$.
\end{proposition}

In the same way, we can also identify tangent bundle $TG$ and
$G\times\mathfrak{g}$, by the local left trivialization of $TG$. In consequence,
we can consider the Lagrangian $L(g,\xi):TG \cong G\times
\mathfrak{g}\to \mathbb{R}$, which is usually the kinetic minus the
potential energy of the system, where $(g,\xi)\in
G\times\mathfrak{g}$, and $\xi \in \mathfrak{g}$, regarded as the
velocity of system. If we introduce the conjugate momentum
$p_i=\frac{\partial L}{\partial \xi^i}$, $i=1,\cdots,n,\; n=dim G$,
and by the Legendre transformation $FL: TG \cong
G\times\mathfrak{g}\to T^\ast G \cong G\times\mathfrak{g}^\ast$,
$(g^i,\xi^i)\to (g^i,p_i)$, we have the Hamiltonian $H(g,p):T^\ast G
\cong G\times\mathfrak{g}^\ast \to \mathbb{R}$ given by
\begin{equation}H(g^i,p_i)=\sum_{i=1}^{n}p_i\xi^i-L(g^i,\xi^i).\label{6.3}
\end{equation}
If the Hamiltonian $H(g,p):T^\ast G\cong G\times\mathfrak{g} \to
\mathbb{R}$ is left cotangent lifted $G$-action invariant, for
$\mu\in\mathfrak{g}^\ast$ we have the associated reduced Hamiltonian
$h_\mu: (T^\ast G)_\mu \cong \mathcal{O}_\mu\to \mathbb{R}$, defined
by $h_\mu\cdot \pi_\mu=H\cdot i_\mu$. By the $(\pm)$-Lie-Poisson
brackets on $\mathfrak{g}^\ast$ and the symplectic structure on the
coadjoint orbit $\mathcal{O}_\mu$, we have the associated
Hamiltonian vector field $X_{h_\mu}$ given by
\begin{equation}   X_{h_\mu}(\nu)=\mp
\operatorname{ad}^\ast_{\delta h_\mu/\delta \nu} \nu,\quad \forall
\nu\in \mathcal{O}_\mu.\label{6.4}
\end{equation}
See Marsden and Ratiu \cite{mara99}. Thus, if the Hamiltonian
$H:T^\ast G\to \mathbb{R}$, the fiber-preserving map $F:T^\ast G\to
T^\ast G$ and the fiber submanifold $W$ of $T^\ast G$ are all left
cotangent lifted $G$-action invariant, then we may give the
$R_P$-reduced RCH system as follows.

\begin{theorem}
The 6-tuple $(T^\ast G,G,\omega_0,H,F,W)$ is a regular point
reducible RCH system on Lie group $G$, where the Hamiltonian $H:
T^\ast G \to \mathbb{R}$, the fiber-preserving map $F: T^\ast G\to
T^\ast G$ and the fiber submanifold $W$ of  $T^\ast G$ are all left
cotangent lifted $G$-action invariant. For a point
$\mu\in\mathfrak{g}^\ast$, the regular value of the momentum map
$\mathbf{J}_L: T^\ast G \rightarrow \mathfrak{g}^\ast$, the
$R_P$-reduced system, that is, the 5-tuple
$(\mathcal{O}_\mu,\omega_{\mathcal{O}_\mu}^{-},h_\mu,f_\mu,W_\mu)$,
is a RCH system, where $\mathcal{O}_\mu \subset \mathfrak{g}^\ast$
is the coadjoint orbit, $\omega_{\mathcal{O}_\mu}^{-}$ is orbit
symplectic form, $h_\mu\cdot \pi_\mu=H\cdot i_\mu$,
$F(\mathbf{J}_L^{-1}(\mu))\subset \mathbf{J}_L^{-1}(\mu) $,
$f_\mu\cdot \pi_\mu=\pi_\mu\cdot F\cdot i_\mu$,
$W \cap \mathbf{J}_L^{-1}(\mu)\neq \emptyset $,
$W_\mu=\pi_\mu(W\cap \mathbf{J}_L^{-1}(\mu))$.
Moreover, two regular point reducible RCH system
$(T^\ast G_i, G_i, \omega_{i0}, H_i, F_i, W_i),$ \; $ i=1,2,$ are
RpCH-equivalent if and only if the associated $R_P$-reduced RCH
systems
$(\mathcal{O}_{i\mu_i},\omega_{\mathcal{O}_{i\mu_i}}^{-},h_{i\mu_i},f_{i\mu_i},W_{i\mu_i}),
\; i=1,2,$ are RCH-equivalent.
\end{theorem}

Next, in order to study the regular reduction of rigid body and
heavy top with internal rotors, we need the regular symplectic
reduction theory of the cotangent bundle $T^\ast Q$, where the
configuration space $Q=G\times V$, and $G$ is a Lie group and $V$ is
a $k$-dimensional vector space. Defined the left $G$-action $\Phi:
G\times Q \rightarrow Q, \; \Phi(g,(h,\theta)):=(gh,\theta)$, for
any $g,h \in G,\; \theta \in V$, that is , the $G$-action on $Q$ is by
the left translation on the first factor $G$, and the $G$ trivial action
on the second factor $V$. Because locally,
$T^\ast Q \cong T^\ast G \times T^\ast V$,
and $T^\ast V\cong V\times V^\ast$, by using the left trivialization
of $T^\ast G$, that is, $T^\ast G \cong G \times \mathfrak{g}^\ast, $ where
$\mathfrak{g}^\ast$ is the dual of $\mathfrak{g}$, and hence we have
that locally, $T^\ast Q= G\times \mathfrak{g}^\ast
\times V \times V^\ast$. If the left $G$-action $\Phi: G\times Q
\rightarrow Q $ is free and proper, then the cotangent lift of the
action to its cotangent bundle $T^\ast Q$, given by $\Phi^{T^*}:
G\times T^*Q \rightarrow T^*Q, \;
\Phi^{T^*}(g,(h,\mu,\theta,\lambda)):=(gh,\mu,\theta,\lambda)$, for
any $g,h \in G,\; \mu \in \mathfrak{g}^\ast, \; \theta \in V, \;
\lambda \in V^\ast$, is also a free and proper action, and the orbit
space $(T^\ast Q)/ G $ is a smooth manifold and $\pi: T^*Q
\rightarrow (T^*Q )/G $ is a smooth submersion. Since $G$ acts
trivially on $\mathfrak{g}^\ast$, $V$ and $V^\ast$, it follows that
$(T^\ast Q)/ G$ is diffeomorphic to $\mathfrak{g}^\ast \times V
\times V^\ast$.\\

For $\mu \in \mathfrak{g}^\ast$, the coadjoint orbit
$\mathcal{O}_\mu \subset \mathfrak{g}^\ast$ has the orbit symplectic
forms $\omega^\pm_{\mathcal{O}_\mu}$. Let $\omega_V$ be the
canonical symplectic form on $T^\ast V \cong V \times V^\ast$ given
by $$\omega_V((\theta_1, \lambda_1),(\theta_2,
\lambda_2))=<\lambda_2,\theta_1> -<\lambda_1,\theta_2>,$$ where
$(\theta_i, \lambda_i)\in V\times V^\ast, \; i=1,2$, $<\cdot,\cdot>$
is the natural pairing between $V^\ast$ and $V$. Thus, we can induce
a symplectic forms $\tilde{\omega}^\pm_{\mathcal{O}_\mu \times V
\times V^\ast}= \pi_{\mathcal{O}_\mu}^\ast
\omega^\pm_{\mathcal{O}_\mu}+ \pi_V^\ast \omega_V$ on the smooth
manifold $\mathcal{O}_\mu \times V \times V^\ast$, where the maps
$\pi_{\mathcal{O}_\mu}: \mathcal{O}_\mu \times V \times V^\ast \to
\mathcal{O}_\mu$ and $\pi_V: \mathcal{O}_\mu \times V \times V^\ast
\to V\times V^\ast$ are canonical projections.\\

On the other hand, the cotangent bundle $T^\ast Q$ has a canonical
symplectic form $\omega_Q$, from locally,
$T^\ast Q \cong T^\ast G \times T^\ast V$ we have that
$\omega_Q= \pi^\ast_1 \omega_0 +\pi^\ast_2
\omega_V$ on $T^\ast Q$, where $\omega_0$ is the canonical
symplectic form on $T^\ast G$ and the maps $\pi_1: Q=G\times V \to
G$ and $\pi_2: Q=G\times V \to V$ are canonical projections. Assume that
the cotangent lift of the left $G$-action $\Phi^{T^*}: G \times
T^\ast Q \to T^\ast Q$ is symplectic, and admits an associated
$\operatorname{Ad}^\ast$-equivariant momentum map $\mathbf{J}_Q:
T^\ast Q \to \mathfrak{g}^\ast$ such that $\mathbf{J}_Q\cdot
\pi^\ast_1=\mathbf{J}_L$, where $\mathbf{J}_L:T^\ast G \rightarrow
\mathfrak{g}^\ast$ is a momentum map of left G-action
on $T^\ast G$ and we assume that it exists,
and $\pi^\ast_1: T^\ast G \to T^\ast Q$. If
$\mu\in\mathfrak{g}^\ast$ is a regular value of $\mathbf{J}_Q$, then
$\mu\in\mathfrak{g}^\ast$ is also a regular value of $\mathbf{J}_L$
and $\mathbf{J}_Q^{-1}(\mu)\cong \mathbf{J}_L^{-1}(\mu)\times V
\times V^\ast$. Denote $G_\mu=\{g\in G|\operatorname{Ad}_g^\ast
\mu=\mu \}$ the isotropy subgroup of coadjoint $G$-action at the
point $\mu\in\mathfrak{g}^\ast$. It follows that $G_\mu$ acts also
freely and properly on $\mathbf{J}_Q^{-1}(\mu)$, the regular point
reduced space $(T^\ast Q)_\mu=\mathbf{J}_Q^{-1}(\mu)/G_\mu\cong
(T^\ast G)_\mu \times V \times V^\ast$ of $(T^\ast Q,\omega_Q)$ at
$\mu$, is a symplectic manifold with symplectic form $\omega_\mu$
uniquely characterized by the relation $\pi_\mu^\ast
\omega_\mu=i_\mu^\ast \omega_Q=i_\mu^\ast \pi^\ast_1 \omega_0
+i_\mu^\ast \pi^\ast_2 \omega_V$, where the map
$i_\mu:\mathbf{J}_Q^{-1}(\mu)\rightarrow T^\ast Q$ is the inclusion
and $\pi_\mu:\mathbf{J}_Q^{-1}(\mu)\rightarrow (T^\ast Q)_\mu$ is
the projection. Because $((T^\ast G)_\mu,\omega_\mu)$ is
symplectically diffeomorphic to
$(\mathcal{O}_\mu,\omega_{\mathcal{O}_\mu}^{-})$, we have that
$((T^\ast Q)_\mu,\omega_\mu)$ is symplectically diffeomorphic to
$(\mathcal{O}_\mu \times V\times
V^\ast,\tilde{\omega}_{\mathcal{O}_\mu \times V \times
V^\ast}^{-})$.\\

Now we identify $TG$ and $G\times \mathfrak{g}$ locally,
by using the left translation, and $TV\cong V\times V$,
then locally, $TQ\cong G\times \mathfrak{g} \times V\times V$.
In consequence, we consider the Lagrangian $L(g,\xi,\theta,\dot{\theta}):TQ
\cong G\times \mathfrak{g} \times V\times V \to \mathbb{R}$, which is usually
the total kinetic minus potential energy of the system, where
$(g,\xi) \in G\times \mathfrak{g}$, and $\theta \in V$, $\xi^i$ and
$\dot{\theta}^j=\frac{\mathrm{d} \theta^j}{\mathrm{d} t}$,
($i=1,\cdots,n, \; j=1,\cdots,k$, $n=\dim G$, $k=\dim V$), regarded
as the velocity of system. If we introduce the conjugate momentum
$p_i=\frac{\partial L}{\partial \xi^i}, \; l_j=\frac{\partial
L}{\partial \dot{\theta}^j}$, $i=1,\cdots,n, \; j=1,\cdots,k,$ and
by the Legendre transformation $FL: TQ \cong G\times
\mathfrak{g}\times V \times V \to T^\ast Q \cong G\times
\mathfrak{g}^\ast \times V \times V^\ast$,
$(g^i,\xi^i,\theta^j,\dot{\theta}^j)\to (g^i,p_i,\theta^j,l_j)$, we
have the Hamiltonian $H(g,p,\theta,l):T^\ast Q \cong G\times
\mathfrak{g}^\ast \times V \times V^\ast \to \mathbb{R}$ given by
\begin{equation}
H(g^i,p_i,\theta^j,l_j)=\sum_{i=1}^{n}p_i\xi^i+\sum_{j=1}^{k}l_j\dot{\theta}^j
-L(g^i,\xi^i,\theta^j,\dot{\theta}^j).\label{6.5}
\end{equation}
If the Hamiltonian $H(g,p,\theta,l):T^\ast Q \cong G\times
\mathfrak{g}^\ast \times V \times V^\ast \to \mathbb{R}$ is left
cotangent lifted $G$-action $\Phi^{T^*}$ invariant, for
$\mu\in\mathfrak{g}^\ast$ we have the associated reduced Hamiltonian
$h_\mu(\nu,\theta,l): (T^\ast Q)_\mu \cong\mathcal{O}_\mu \times V
\times V^\ast \to \mathbb{R}$, defined by $h_\mu\cdot \pi_\mu=H\cdot
i_\mu$. Note that for $F,K: T^\ast V \cong V \times V^\ast \to
\mathbb{R}$, by using the canonical symplectic form $\omega_V$ on
$T^\ast V \cong V \times V^\ast$, we can define the Poisson bracket
$\{\cdot,\cdot\}_V$ on $T^\ast V$ as follows
$$\{F,K\}_V(\theta,\lambda)=<\frac{\delta F}{\delta \theta},
\frac{\delta K}{\delta \lambda}>- <\frac{\delta K}{\delta
\theta},\frac{\delta F}{\delta \lambda}>
$$
If $\theta_i, \; i=1,\cdots, k,$ is a base of $V$, and $\lambda_i,
\; i=1,\cdots, k,$ a base of $V^\ast$, then we have that
\begin{equation*}\{F,K\}_V(\theta,\lambda)=\sum_{i=1}^k(\frac{\partial F}{\partial \theta_i}
\frac{\partial K}{\partial \lambda_i}- \frac{\partial K}{\partial
\theta_i}\frac{\partial F}{\partial \lambda_i}).\label{6.6}
\end{equation*}
Thus, by the $(\pm)$-Lie-Poisson brackets on $\mathfrak{g}^\ast$ and
the Poisson bracket $\{\cdot,\cdot\}_V$ on $T^\ast V$, for $F,K:
\mathfrak{g}^\ast \times V \times V^\ast \to \mathbb{R}$, we can
define the Poisson bracket on $\mathfrak{g}^\ast \times V \times
V^\ast$ as follows
\begin{equation} \{F,K\}_\pm(\mu,\theta,\lambda)= \{F,K\}_\pm(\mu)+
\{F,K\}_V(\theta,\lambda)\\ = \pm<\mu,[\frac{\delta F}{\delta
\mu},\frac{\delta K}{\delta\mu}]>+
 \sum_{i=1}^k(\frac{\partial F}{\partial \theta_i} \frac{\partial
K}{\partial \lambda_i}- \frac{\partial K}{\partial
\theta_i}\frac{\partial F}{\partial \lambda_i}).\label{6.6}
\end{equation}
See Krishnaprasad and Marsden \cite{krma87}. In particular, for
$F_\mu,K_\mu: \mathcal{O}_\mu \times V \times V^\ast \to
\mathbb{R}$, we have that $\tilde{\omega}_{\mathcal{O}_\mu \times V
\times V^\ast}^{-}(X_{F_\mu}, X_{K_\mu})=
\{F_\mu,K_\mu\}_{-}|_{\mathcal{O}_\mu \times V \times V^\ast}$.
Moreover, for the reduced Hamiltonian $h_\mu(\nu,\theta,l):
\mathcal{O}_\mu \times V \times V^\ast \to \mathbb{R}$, we have the
Hamiltonian vector field
$X_{h_\mu}(K_\mu)=\{K_\mu,h_\mu\}_{-}|_{\mathcal{O}_\mu \times V
\times V^\ast}. $ Thus, if the Hamiltonian $H: T^\ast Q \to
\mathbb{R}$, the fiber-preserving map $F: T^\ast Q \to T^\ast Q$ and
the fiber submanifold $W$ of $T^\ast Q$ are all left cotangent
lifted $G$-action $\Phi^{T^*}$ invariant, then we have the following
theorem.

\begin{theorem}
The 6-tuple $(T^\ast Q,G,\omega_Q,H,F,W)$ is a regular point
reducible RCH system, where $Q=G\times V$, and $G$ is a Lie group
and $V$ is a $k$-dimensional vector space, and the Hamiltonian
$H:T^\ast Q \to \mathbb{R}$, the fiber-preserving map $F: T^\ast Q
\to T^\ast Q$ and the fiber submanifold $W$ of $T^\ast Q$ are all
left cotangent lifted $G$-action $\Phi^{T^*}$ invariant. For a point
$\mu\in\mathfrak{g}^\ast$, the regular value of the momentum map
$\mathbf{J}_Q: T^\ast Q \rightarrow \mathfrak{g}^\ast$, the
$R_P$-reduced system, that is, the 5-tuple $(\mathcal{O}_\mu \times
V\times V^\ast,\tilde{\omega}_{\mathcal{O}_\mu \times V \times
V^\ast }^{-},h_\mu,f_\mu,W_\mu)$, is a RCH system, where
$\mathcal{O}_\mu \subset \mathfrak{g}^\ast$ is the coadjoint orbit,
$\tilde{\omega}_{\mathcal{O}_\mu \times V \times V^\ast }^{-}$ is
orbit symplectic form on $\mathcal{O}_\mu \times V\times V^\ast $,
$h_\mu\cdot \pi_\mu=H\cdot i_\mu$,
$F(\mathbf{J}_Q^{-1}(\mu))\subset \mathbf{J}_Q^{-1}(\mu) $,
$f_\mu\cdot \pi_\mu=\pi_\mu\cdot F\cdot i_\mu$,
$W \cap \mathbf{J}_Q^{-1}(\mu)\neq \emptyset $,
$W_\mu=\pi_\mu(W\cap \mathbf{J}_Q^{-1}(\mu))
\subset \mathcal{O}_\mu \times V \times V^\ast$.
Moreover, two regular point reducible RCH system $(T^\ast Q_i, G_i,
\omega_{i0}, H_i, F_i, W_i),$ \; $ i=1,2,$ are RpCH-equivalent if
and only if the associated $R_P$-reduced RCH systems
$(\mathcal{O}_{i\mu_i}\times V_i\times
V_i^\ast,\tilde{\omega}_{\mathcal{O}_{i\mu_i}}^{-},h_{i\mu_i},f_{i\mu_i},W_{i\mu_i}),
\; i=1,2,$ are RCH-equivalent.
\end{theorem}

The third, in order to study the regular reduction of heavy top we
need to the theory of Hamiltonian reduction by stages for semidirect
product Lie group. See Marsden et al. \cite{mamiorpera07}. Assume
that $S=G\circledS V$ is a semidirect product Lie group, where $V$
is a vector space and $V^\ast$ its dual space, $G$ is a Lie group
acting on the left by linear maps on $V$, and $\mathfrak{g}$ its Lie
algebra and $\mathfrak{g}^\ast$ the dual of $\mathfrak{g}$. Note
that $G$ also acts on the left on the dual space $V^\ast$ of $V$,
and the action by an element $g$ on $V^\ast$ is the transpose of the
action of $g^{-1}$ on $V$. As a set, the underlying manifold of $S$
is $G\times V$ and the multiplication on $S$ is given by
\begin{equation}
  (g_1,v_1)(g_2,v_2):=(g_1g_2,v_1+\sigma(g_1)v_2),\quad g_1,g_1\in
  G,\quad v_1,v_2\in V\label{6.7}
\end{equation}
where $\sigma:G\to \operatorname{Aut}(V)$ is a representation of the
Lie group $G$ on $V$, $\operatorname{Aut}(V)$ denotes the Lie group
of linear isomorphisms of $V$ onto itself whose Lie algebra is
$\operatorname{End}(V)$, the space of all linear maps of $V$ to
itself.\\

The Lie algebra of $S$ is the semidirect product of Lie algebras
$\mathfrak{s}=\mathfrak{g}\circledS V$, $\mathfrak{s}^\ast$ is the
dual of $\mathfrak{s}$, that is, $\mathfrak{s}^\ast=
(\mathfrak{g}\circledS V)^\ast$. The underlying vector space of
$\mathfrak{s}$ is $\mathfrak{g}\times V$ and the Lie bracket on
$\mathfrak{s}$ is given by
\begin{equation} [(\xi_1,v_1),(\xi_2,v_2)]=([\xi_1,\xi_2],\sigma
'(\xi_1)v_2-\sigma '(\xi_2)v_1),\quad\forall \xi_1,\xi_2\in
\mathfrak{g},\quad v_1,v_2\in V\label{6.8}
\end{equation}
where $\sigma ':\mathfrak{g}\to \operatorname{End}(V)$ is the
induced Lie algebra representation given by \begin{equation}
\sigma'(\xi)v:=\left.\frac{\mathrm{d}}{\mathrm{d}
t}\right|_{t=0}\sigma(\exp t\xi)v,\quad \xi\in\mathfrak{g},\quad
v\in V\label{6.9}
\end{equation}
Identify the underlying vector space of $\mathfrak{s}^\ast$ with
$\mathfrak{g}^\ast\times V^\ast$ by using the duality pairing on
each factor. One can give the formula for the $(\pm)$-Lie-Poisson
bracket on the semidirect product $\mathfrak{s}^\ast $ as follows,
that is, for $F,K:\mathfrak{s}^\ast\to\mathbb{R}$, their semidirect
product bracket is given by \begin{equation}
\{F,K\}_\pm(\mu,a)=\pm\langle\mu,[\frac{\delta F}{\delta
\mu},\frac{\delta K}{\delta \mu}]\rangle\pm\langle a,\frac{\delta
F}{\delta \mu}\cdot\frac{\delta K}{\delta a}-\frac{\delta K}{\delta
\mu}\cdot\frac{\delta F}{\delta a}\rangle\label{6.10}
\end{equation}
where $(\mu,a)\in \mathfrak{s}^\ast$ and $\dfrac{\delta F}{\delta
\mu}\in \mathfrak{g}$, $\dfrac{\delta F}{\delta a}\in V$ are the
functional derivatives. Moreover, the Hamiltonian vector field of
a smooth function $H:\mathfrak{s}^\ast\to\mathbb{R}$ is given by
\begin{equation} X_H(\mu,a)=\mp(\operatorname{ad}_{\delta H/\delta \mu}^\ast
\mu-\rho_{\delta H/\delta a}^\ast a, \; \frac{\delta H}{\delta
\mu}\cdot a),\label{6.11}
\end{equation}
where the infinitesimal action of $\mathfrak{g}$ on $V$ can be
denoted by $\xi\cdot v=\rho_v(\xi)$, for any $\xi\in\mathfrak{g}$,
$v\in V$ and the map $\rho_v:\mathfrak{g}\to V$ is the derivative of
the map $g\mapsto gv$ at the identity and $\rho_v^\ast:
V^\ast\to\mathfrak{g}^\ast$ is its dual.\\

Now we consider a symplectic action of $S$ on a symplectic manifold
$P$ and assume that this action has an
$\operatorname{Ad}^\ast$-equivariant momentum map $\mathbf{J}_S:P\to
\mathfrak{s}^\ast$. On the one hand, we can regard $V$ as a normal
subgroup of $S$, it also acts on $P$ and has a momentum map
$\mathbf{J}_V:P\to V^\ast$ given by $\mathbf{J}_V=i_V^\ast\cdot
\mathbf{J}_S$, where $i_V:V\to \mathfrak{s};\;v\mapsto (0,v)$ is the
inclusion, and $i_V^\ast: \frak{s}^\ast \to V^\ast$ is its dual.
$\mathbf{J}_V$ is called the second component of $\mathbf{J}_S$. On
the other hand, we can also regard $G$ as a subgroup of $S$ by the
inclusion $i_G:G\to S$, $g\mapsto (g,0)$. Thus, $G$ also has a
momentum map $\mathbf{J}_G:P\to \mathfrak{g}^\ast$ given by
$\mathbf{J}_G=i_G^\ast\cdot \mathbf{J}_S$, which is called the first
component of $\mathbf{J}_S$. Moreover, from the
$\operatorname{Ad}^\ast$-equivariance of $\mathbf{J}_S$ under
$G$-action, we know that $\mathbf{J}_V$ is also
$\operatorname{Ad}^\ast$-equivariant under $G$-action. Thus, we can
carry out reduction of $P$ by $S$ at a regular value
$\sigma=(\mu,a)\in \mathfrak{s}^\ast$ of the momentum map
$\mathbf{J}_S$ in two stages using the following procedure. (i)First
reduce $P$ by $V$ at the value $a\in V^\ast$, and get the reduced
space $P_a=\mathbf{J}_V^{-1}(a)/V$. Since the reduction is by the
Abelian group $V$, so the quotient is done using the whole of $V$.
(ii)The isotropy subgroup $G_a\subset G$, consists of elements of
$G$ that leave the point $a\in V^\ast$ fixed using the action of $G$
on $V^\ast$. One can prove that the group $G_a$ leaves the set
$\mathbf{J}_V^{-1}(a)\subset P$ invariant, and acts symplectically
on the reduced space $P_a$ and has a naturally induced momentum map
$\mathbf{J}_a:P_a\to \mathfrak{g}_a^\ast$, where $\mathfrak{g}_a$ is
the Lie algebra of the isotropy subgroup $G_a$ and
$\mathfrak{g}_a^\ast$ is its dual. (iii)Reduce the first reduced
space $P_a$ at the point $\mu_a=\mu|_{\mathfrak{g}^\ast_a}\in
\mathfrak{g}_a^\ast$, one can get the second reduced space
$(P_a)_{\mu_a}=\mathbf{J}_a^{-1}(\mu_a)/(G_a)_{\mu_a}$. Thus, we can
give the following proposition on the reduction by stages for
semidirect products, see Marsden et al. \cite{mamiorpera07}.

\begin{proposition}
  The reduced space $(P_a)_{\mu_a}$ is symplectically
  diffeomorphic to the reduced space $P_\sigma$ obtained by
  reducing $P$ by $S$ at the regular point
  $\sigma=(\mu,a)\in\mathfrak{s}^\ast$.
\end{proposition}

In particular, we can choose that $P=T^\ast S$, where $S=G\circledS
V$ is a semidirect product Lie group, with the cotangent lift action
of $S$ on $T^\ast S$ induced by left translations of $S$ on itself.
Since the reduction of $T^\ast S$ by the action of $V$ can give a
space which is isomorphic to $T^\ast G$, from the above reduction by
stages proposition for semidirect products we can get the following
semidirect product reduction proposition.

\begin{proposition}
  The reduction of $T^\ast G$ by $G_a$ at the regular values
  $\mu_a=\mu|_{\mathfrak{g}^\ast_a}$ gives a space which is isomorphic to
  the coadjoint orbit $\mathcal{O}_\sigma\subset\mathfrak{s}^\ast$
  through the point
  $\sigma=(\mu,a)\in\mathfrak{s}^\ast$, where $\mathfrak{s}^\ast$ is the dual of the Lie
  algebra $\mathfrak{s}$ of $S$.
\end{proposition}

Thus, from the above proposition we know that the reduced space of
the heavy top is obtained by the reduction of $T^\ast
\textmd{SE}(3)$ by left action of
$\textmd{SE}(3)=\textmd{SO}(3)\circledS \mathbb{R}^3$, which is a
coadjoint orbit in $\mathfrak{se}^\ast(3)$. Moreover, the
configuration space of the heavy top with internal rotors is
$Q=\textmd{SE}(3)\times V$, and the reduced space is symplectically
diffeomorphic to a leaf of Poisson manifold
$\mathfrak{se}^\ast(3)\times V\times V^*$. In consequence, we can
deal with uniformly the symplectic reduction of the rigid body,
heavy top, as well as them with internal rotors, such that we can
state that all these systems are the regular point reducible RCH
systems and can give their RCH-equivalences.

\subsection{Rigid Body and Heavy Top }

In this subsection, by using the above method, we describe uniformly
the rigid body and heavy top as well as them with internal rotors
(or external force torques) as the regular point reducible RCH
systems on the rotation group $\textmd{SO}(3)$ and on the Euclidean
group $\textmd{SE}(3)$, as well as on their generalizations,
respectively, and give their $R_P$-reduced RCH systems and discuss
their RCH-equivalence. Note that our description of the motion and
the equations of rigid body and heavy top follows some of the
notations and conventions in
Marsden and Ratiu \cite{mara99}, Marsden \cite{ma92}.\\

\noindent {\bf (1). Rigid Body with External Force Torque.}\\

In the following we take Lie group $G= \textmd{SO}(3), $ and state
the rigid body with external force torque to be a regular point
reducible RCH system. It is well known that, usually, the
configuration space for a $3$-dimensional rigid body moving freely
in space is $\textmd{SE}(3)$, the six dimensional group of Euclidean
(rigid) transformations of three-dimensional space $\mathbb{R}^3$,
that is, all possible rotations and translations. If translations
are ignored and only rotations are considered, then the
configuration space $Q$ is $\textmd{SO}(3)$, consists of all
orthogonal linear transformations of Euclidean three space to
itself, which have determinant one. Its Lie algebra, denoted
$\mathfrak{so}(3)$, consists of all $3\times 3$ skew matrices, and
we can identify the Lie algebra $(\mathfrak{so}(3), [,])$ with
$(\mathbb{R}^3, \times )$. Denote $\mathfrak{so}^\ast(3)$ the
dual of the Lie algebra $\mathfrak{so}(3)$, and we also identity
$\mathfrak{so}^\ast(3)$ with $\mathbb{R}^3$ by pairing the Euclidean
inner product. Since the functional derivative of a function defined
on $\mathbb{R}^3$ is equal to the usual gradient of the function,
from (\ref{6.1}) we know that the Lie-Poisson bracket on
$\mathfrak{so}^\ast(3)$ take the form
\begin{equation}   \{f,g\}_\pm(\Pi)=\pm\Pi\cdot (\nabla_\Pi
f\times \nabla_\Pi g), \;\; \forall f,g\in
C^\infty(\mathfrak{so}^\ast(3)),\;\; \Pi \in
\mathfrak{so}^\ast(3).\label{6.12}
\end{equation}
The phase space of a rigid body is the cotangent bundle $T^\ast G
=T^\ast \textmd{SO}(3)$, and by using the local left trivialization,
locally, $T^\ast \textmd{SO}(3)\cong \textmd{SO}(3)\times
\mathfrak{so}^\ast(3)$, with the canonical symplectic form. Assume
that Lie group $G=\textmd{SO}(3)$ acts freely and properly by the
left translation on $\textmd{SO}(3)$, then the action of
$\textmd{SO}(3)$ on the phase space $T^\ast \textmd{SO}(3)$ is by
the cotangent lift of left translation at the identity, that is, $\Phi:
\textmd{SO}(3)\times T^\ast \textmd{SO}(3)\cong \textmd{SO}(3)\times
\textmd{SO}(3)\times \mathfrak{so}^\ast(3)\to \textmd{SO}(3)\times
\mathfrak{so}^\ast(3),$ given by $\Phi(B,(A,\Pi))=(BA,\Pi)$, for any
$A,B\in \textmd{SO}(3), \; \Pi \in \mathfrak{so}^\ast(3)$, which is
also free and proper. Assume that the action is symplectic, and admits an associated
$\operatorname{Ad}^\ast$-equivariant momentum map $\mathbf{J}:T^\ast
\textmd{SO}(3)\to \mathfrak{so}^\ast(3)$ for the left
$\textmd{SO}(3)$ action. If $\Pi \in \mathfrak{so}^\ast(3)$ is a
regular value of $\mathbf{J}$, then the regular point reduced space
$(T^\ast
\textmd{SO}(3))_\Pi=\mathbf{J}^{-1}(\Pi)/\textmd{SO}(3)_\Pi$ is
symplectically diffeomorphic to the coadjoint orbit $\mathcal{O}_\Pi
\subset \mathfrak{so}^\ast(3)$.\\

Let $I$ be the moment of inertia tensor computed with respect to a
body fixed frame, which, in a principal body frame, we may represent
by the diagonal matrix diag $(I_1,I_2,I_3)$. Let
$\Omega=(\Omega_1,\Omega_2,\Omega_3)$ be the vector of angular
velocities computed with respect to the axes fixed in the body and
$(\Omega_1,\Omega_2,\Omega_3)\in \mathfrak{so}(3)$. Consider the
Lagrangian $L(A,\Omega): \textmd{SO}(3)\times\mathfrak{so}(3)\to
\mathbb{R}$, which is given by
$L(A,\Omega)=\dfrac{1}{2}\langle\Omega,\Omega\rangle
=\dfrac{1}{2}(I_1\Omega_1^2+I_2\Omega_2^2+I_3\Omega_3^2),$ where
$A\in \textmd{SO}(3)$, $(\Omega_1,\Omega_2,\Omega_3)\in
\mathfrak{so}(3)$. If we introduce the conjugate angular momentum
$\Pi_i=\dfrac{\partial L}{\partial \Omega_i}=I_i\Omega_i$,
$i=1,2,3$, which is also computed with respect to a body fixed
frame, and by the Legendre transformation $FL:
\textmd{SO}(3)\times\mathfrak{so}(3)\to \textmd{SO}(3)\times
\mathfrak{so}^\ast(3), \; (A,\Omega)\to(A,\Pi)$, where
$\Pi=(\Pi_1,\Pi_2,\Pi_3)\in \mathfrak{so}^\ast(3)$, we have the
Hamiltonian $H(A,\Pi): \textmd{SO}(3)\times \mathfrak{so}^\ast(3)\to
\mathbb{R}$ given by
\begin{align*} H(A,\Pi)&=\Omega \cdot \Pi-L(A,\Omega)
=\frac{1}{2}(\frac{\Pi_1^2}{I_1}+\frac{\Pi_2^2}{I_2}+\frac{\Pi_3^2}{I_3}).
\end{align*}
From the above expression of the Hamiltonian, we know that
$H(A,\Pi)$ is invariant under the left $\textmd{SO}(3)$-action. For
the case $\Pi_0 =\mu \in \mathfrak{so}^\ast(3)$ is a regular value
of $\mathbf{J}$, we have the reduced Hamiltonian
$h_\mu(\Pi):\mathcal{O}_\mu\subset \mathfrak{so}^\ast(3)\to
\mathbb{R}$ given by $h_\mu(\Pi)=H(A,\Pi)|_{\mathcal{O}_\mu}$. From
the Lie-Poisson bracket on $\mathfrak{g}^\ast$, we can get the rigid
body Poisson bracket on $\mathfrak{so}^\ast(3)$, that is, for $F,K:
\mathfrak{so}^\ast(3)\to \mathbb{R}, $ we have that
$\{F,K\}_{-}(\Pi)=-\Pi\cdot(\nabla_\Pi F\times \nabla_\Pi K)$. In
particular, for $F_\mu,K_\mu: \mathcal{O}_\mu \to \mathbb{R}$, we
have that $\omega_{\mathcal{O}_\mu}^{-}(X_{F_\mu}, X_{K_\mu})=
\{F_\mu,K_\mu\}_{-}|_{\mathcal{O}_\mu}$. Moreover, for reduced
Hamiltonian $h_\mu(\Pi): \mathcal{O}_\mu \to \mathbb{R}$, we have
the Hamiltonian vector field
$X_{h_\mu}(K_\mu)=\{K_\mu,h_\mu\}_{-}|_{\mathcal{O}_\mu },$ and
hence we have that
\begin{align*} \frac{\mathrm{d}\Pi}{\mathrm{d}t}
=X_{h_\mu}(\Pi)=\{\Pi,h_\mu(\Pi)\}_{-}|_{\mathcal{O}_\mu}
=-\Pi\cdot(\nabla_\Pi\Pi\times\nabla_\Pi h_\mu)
=-\nabla_\Pi\Pi\cdot(\nabla_\Pi h_\mu\times \Pi)=\Pi\times\Omega,
\end{align*}
since $\nabla_\Pi\Pi=1$ and $\nabla_\Pi h_\mu=\Omega$. Thus, the
equations of motion for rigid body is given by
\begin{equation}
  \dfrac{\mathrm{d} \Pi}{\mathrm{d} t}=\Pi\times\Omega.\label{6.13}
\end{equation}
From Theorem 6.2 if we consider the rigid body with an external
force torque $u: T^\ast \textmd{SO}(3) \to W $, and $u\in W \cap
\mathbf{J}^{-1}(\mu)$ is invariant under the left
$\textmd{SO}(3)$-action, then the external force torque $u$ can be
regarded as a control of the rigid body, and its reduced control
$u_\mu: \mathcal{O}_\mu \to W_\mu$ is given by
$u_\mu(\Pi)=\pi_\mu(u(A,\Pi))= u(A,\Pi)|_{\mathcal{O}_\mu}, $ where
$\pi_\mu: \mathbf{J}^{-1}(\mu) \rightarrow \mathcal{O}_\mu, $ and
$W_\mu=\pi_\mu(W \cap \mathbf{J}^{-1}(\mu))$. Thus, the equations of motion for the rigid body
with external force torques $u: T^\ast \textmd{SO}(3) \to T^\ast
\textmd{SO}(3)$ are given by
\begin{equation}   \dfrac{\mathrm{d}
\Pi}{\mathrm{d} t}=\Pi\times\Omega +
\mbox{vlift}(u_\mu),\label{6.14}
\end{equation}
where $\mbox{vlift}(u_\mu)= \mbox{vlift}(u_\mu)X_{h_\mu} \in
T\mathcal{O}_\mu .$  To sum up the above discussion, we have the
following proposition.
\begin{proposition}
The 5-tuple $(T^\ast\textmd{SO}(3), \textmd{SO}(3),\omega_0,H,u )$
is a regular point reducible RCH system. For a point $\mu \in
\mathfrak{so}^\ast(3)$, the regular value of the momentum map
$\mathbf{J}:T^\ast \textmd{SO}(3)\to \mathfrak{so}^\ast(3)$, the
$R_P$-reduced system is the 4-tuple
$(\mathcal{O}_\mu,\omega_{\mathcal{O}_\mu}^{-},h_\mu,u_\mu),$ where
$\mathcal{O}_\mu \subset \mathfrak{so}^\ast(3)$ is the coadjoint
orbit, $\omega_{\mathcal{O}_\mu}^{-}$ is the orbit symplectic form on
$\mathcal{O}_\mu $, $h_\mu(\Pi)=H(A,\Pi)|_{\mathcal{O}_\mu}$,
$u_\mu(\Pi)=\pi_\mu(u(A,\Pi))= u(A,\Pi)|_{\mathcal{O}_\mu}$, and its
equation of motion is given by (\ref{6.14}).
\end{proposition}

\noindent {\bf (2). The Rigid Body with Internal Rotors.}\\

In the following we take Lie group $G= \textmd{SO}(3), \; V=
S^1\times S^1\times S^1, \; Q= G\times V $ and state the rigid body
with three symmetric internal rotors to be a regular point reducible
RCH system. We consider a rigid body (to be called the carrier body)
carrying three symmetric rotors. Denote $O$ the center of mass of
the system in the body frame and at $O$ place a set of (orthonormal)
body axes. Assume that the rotor and the body coordinate axes are
aligned with principal axes of the carrier body. The rotor spins
under the influence of a torque $u$ acting on the rotor. The
configuration space is $Q=\textmd{SO}(3)\times V$, where
$V=S^1\times S^1\times S^1$, with the first factor being rigid body
attitude and the second factor being the angles of rotors. The
corresponding phase space is the cotangent bundle $T^\ast Q$ and locally,
$T^\ast Q \cong T^\ast \textmd{SO}(3)\times T^\ast V$,
where $T^\ast V = T^\ast (S^1\times S^1\times S^1)\cong T^\ast
\mathbb{R}^3$ locally, with the canonical symplectic form $\omega_Q$.
By using the local left trivialization, locally,
$T^\ast \textmd{SO}(3)\cong \textmd{SO}(3)\times \mathfrak{so}^\ast(3)$
and $T^*\mathbb{R}^3 \cong \mathbb{R}^3 \times \mathbb{R}^{3*}$,
then we have that locally,
$T^*Q \cong \textmd{SO}(3)\times \mathfrak{so}^\ast(3)
\times \mathbb{R}^3 \times \mathbb{R}^{3*}$.
Assume that Lie
group $G=\textmd{SO}(3)$ acts freely and properly on $Q$ by the left
translation on the first factor $\textmd{SO}(3)$ and
the trivial action on the second factor $V$.
Then the action of $\textmd{SO}(3)$
on the phase space $T^\ast Q$ is by the cotangent lift of left
$\textmd{SO}(3)$ action on $Q$, that is, $\Phi:
\textmd{SO}(3)\times T^\ast \textmd{SO}(3)\times T^\ast V \cong
\textmd{SO}(3)\times \textmd{SO}(3)\times \mathfrak{so}^\ast(3)
\times \mathbb{R}^3 \times \mathbb{R}^{3*}\to \textmd{SO}(3)\times
\mathfrak{so}^\ast(3)\times \mathbb{R}^3 \times \mathbb{R}^{3*},$ given
by $\Phi(B,(A,\Pi,\alpha,l))=(BA,\Pi,\alpha,l)$, for any $A,B\in
\textmd{SO}(3), \; \Pi \in \mathfrak{so}^\ast(3), \; \alpha \in
\mathbb{R}^3, \; l \in \mathbb{R}^{3*}$, which is also free and proper.
Assume that the left $\textmd{SO}(3)$ action
is symplectic and admits an
associated $\operatorname{Ad}^\ast$-equivariant momentum map
$\mathbf{J}_Q: T^\ast Q \cong \textmd{SO}(3)\times
\mathfrak{so}^\ast(3) \times \mathbb{R}^3 \times \mathbb{R}^{3*} \to
\mathfrak{so}^\ast(3)$ for the left $\textmd{SO}(3)$ action. If $\Pi
\in \mathfrak{so}^\ast(3)$ is a regular value of $\mathbf{J}_Q$,
then the regular point reduced space $(T^\ast Q)_\Pi=
\mathbf{J}^{-1}_Q(\Pi)/\textmd{SO}(3)_\Pi$ is symplectically
diffeomorphic to the coadjoint orbit $\mathcal{O}_\Pi \times
\mathbb{R}^3 \times \mathbb{R}^{3*} \subset \mathfrak{so}^\ast(3)
\times \mathbb{R}^3 \times \mathbb{R}^{3*}$.\\

Let $I=diag(I_1,I_2,I_3)$ be the moment of inertia of the carrier
body in the principal body-fixed frame, and $J_i,\; i=1,2,3$ be the
moments of inertia of rotors around their rotation axes. Let
$J_{ik},\; i=1,2,3, \;k=1,2,3,$ be the moments of inertia of the
$i$th rotor with $i=1,2,3,$ around the $k$th principal axis with
$k=1,2,3,$ respectively, and denote by
$\bar{I}_i=I_i+J_{1i}+J_{2i}+J_{3i}-J_{ii}, \; i=1,2,3$. Let
$\Omega=(\Omega_1,\Omega_2,\Omega_3)$ be the vector of body angular
velocities computed with respect to the axes fixed in the body and
$(\Omega_1,\Omega_2,\Omega_3)\in \mathfrak{so}(3)$. Let $\alpha_i,\;
i=1,2,3,$ be the relative angles of rotors and
$\dot{\alpha}=(\dot{\alpha_1},\dot{\alpha_2},\dot{\alpha_3})$ the
vector of rotor relative angular velocities about the principal axes
with respect to a carrier body fixed frame.\\

Now, by using the local left trivialization, locally,
$T\textmd{SO}(3)\cong \textmd{SO}(3)\times \mathfrak{so}(3)$
and $T\mathbb{R}^3 \cong \mathbb{R}^3 \times \mathbb{R}^{3}$,
then we have that locally,
$TQ \cong \textmd{SO}(3)\times \mathfrak{so}(3)
\times \mathbb{R}^3 \times \mathbb{R}^{3}$.
We consider the Lagrangian
of the system $L(A,\Omega,\alpha,\dot{\alpha}):
\textmd{SO}(3)\times\mathfrak{so}(3)\times\mathbb{R}^3\times\mathbb{R}^3\to
\mathbb{R}$, which is the total kinetic energy of the rigid body
plus the total kinetic energy of rotors, given by
$$L(A,\Omega,\alpha,\dot{\alpha})=\dfrac{1}{2}[\bar{I}_1\Omega_1^2
+\bar{I}_2\Omega_2^2+\bar{I}_3\Omega_3^2
+J_1(\Omega_1+\dot{\alpha}_1)^2+J_2(\Omega_2+\dot{\alpha}_2)^2
+J_3(\Omega_3+\dot{\alpha}_3)^2],$$ where $A\in \textmd{SO}(3)$,
$\Omega=(\Omega_1,\Omega_2,\Omega_3)\in \mathfrak{so}(3)$,
$\alpha=(\alpha_1,\alpha_2,\alpha_3)\in \mathbb{R}^3$,
$\dot{\alpha}=(\dot{\alpha}_1,\dot{\alpha}_2,\dot{\alpha}_3)\in
\mathbb{R}^3$. If we introduce the conjugate angular momentum, given
by $\Pi_i= \dfrac{\partial L}{\partial \Omega_i}
=\bar{I}_i\Omega_i+J_i(\Omega_i+\dot{\alpha}_i),\quad l_i =
\dfrac{\partial L}{\partial \dot{\alpha}_i}
=J_i(\Omega_i+\dot{\alpha}_i),\quad i=1,2,3,$ and by the Legendre
transformation $FL:
\textmd{SO}(3)\times\mathfrak{so}(3)\times\mathbb{R}^3\times\mathbb{R}^3\to
\textmd{SO}(3)\times
\mathfrak{so}^\ast(3)\times\mathbb{R}^3\times\mathbb{R}^{3*},\quad
(A,\Omega,\alpha,\dot{\alpha})\to(A,\Pi,\alpha,l)$, where
$\Pi=(\Pi_1,\Pi_2,\Pi_3)\in \mathfrak{so}^\ast(3)$,
$l=(l_1,l_2,l_3)\in \mathbb{R}^{3*}$, we have the Hamiltonian
$H(A,\Pi,\alpha,l): \textmd{SO}(3)\times \mathfrak{so}^\ast
(3)\times\mathbb{R}^3\times\mathbb{R}^{3*}\to \mathbb{R}$ given by
\begin{align*} H(A,\Pi,\alpha,l)&=\Omega\cdot
\Pi+\dot{\alpha}\cdot l-L(A,\Omega,\alpha,\dot{\alpha})\\
&=\frac{1}{2}[\frac{(\Pi_1-l_1)^2}{\bar{I}_1}+\frac{(\Pi_2-l_2)^2}{\bar{I}_2}
+\frac{(\Pi_3-l_3)^2}{\bar{I}_3}+\frac{l_1^2}{J_1}+\frac{l_2^2}{J_2}+\frac{l_3^2}{J_3}].
\end{align*}
From the above expression of the Hamiltonian, we know that
$H(A,\Pi,\alpha,l)$ is invariant under the left
$\textmd{SO}(3)$-action. For the case $\Pi_0 =\mu \in
\mathfrak{so}^\ast(3)$ is the regular value of $\mathbf{J}_Q$, we
have the reduced Hamiltonian $h_\mu(\Pi, \alpha,l):\mathcal{O}_\mu
\times\mathbb{R}^3\times\mathbb{R}^3 (\subset \mathfrak{so}^\ast
(3)\times\mathbb{R}^3\times\mathbb{R}^3)\to \mathbb{R}$ given by
$h_\mu(\Pi,\alpha,l)=H(A,\Pi,\alpha,l)|_{\mathcal{O}_\mu
\times\mathbb{R}^3\times\mathbb{R}^3}.$ From (6.6), the rigid body
Poisson bracket on $\mathfrak{so}^\ast(3)$ and the Poisson bracket
on $T^\ast \mathbb{R}^3$, we can get the Poisson bracket on $T^\ast
Q$, that is, for $F,K:
\mathfrak{so}^\ast(3)\times\mathbb{R}^3\times\mathbb{R}^{3*} \to
\mathbb{R}, $ we have that
$\{F,K\}_{-}(\Pi,\alpha,l)=-\Pi\cdot(\nabla_\Pi F\times \nabla_\Pi
K)+ \{F,K\}_V(\alpha,l)$. In particular, for $F_\mu,K_\mu:
\mathcal{O}_\mu \times\mathbb{R}^3\times\mathbb{R}^{3*} \to
\mathbb{R}$, we have that $\tilde{\omega}_{\mathcal{O}_\mu
\times\mathbb{R}^3\times\mathbb{R}^{3*}}^{-}(X_{F_\mu}, X_{K_\mu})=
\{F_\mu,K_\mu\}_{-}|_{\mathcal{O}_\mu
\times\mathbb{R}^3\times\mathbb{R}^{3*} }$. Moreover, for the reduced
Hamiltonian $h_\mu(\Pi,\alpha,l): \mathcal{O}_\mu
\times\mathbb{R}^3\times\mathbb{R}^{3*} \to \mathbb{R}$, we have the
Hamiltonian vector field
$X_{h_\mu}(K_\mu)=\{K_\mu,h_\mu\}_{-}|_{\mathcal{O}_\mu
\times\mathbb{R}^3\times\mathbb{R}^{3*} },$ and hence we have that
\begin{align*}
\frac{\mathrm{d}\Pi}{\mathrm{d}t}
&=X_{h_\mu}(\Pi)(\Pi,\alpha,l)=\{\Pi,h_\mu\}_{-}(\Pi,\alpha,l)\\
&=-\Pi\cdot(\nabla_\Pi\Pi\times\nabla_\Pi h_\mu)+
\sum_{i=1}^3(\frac{\partial \Pi}{\partial \alpha_i} \frac{\partial
h_\mu}{\partial l_i}- \frac{\partial h_\mu}{\partial
\alpha_i}\frac{\partial \Pi}{\partial l_i})
=-\nabla_\Pi\Pi\cdot(\nabla_\Pi h_\mu\times \Pi)=\Pi\times\Omega,
\end{align*}
since $\nabla_\Pi\Pi=1$,  $\nabla_\Pi h_\mu=\Omega$ and
$\frac{\partial \Pi}{\partial \alpha_i}= \frac{\partial
h_\mu}{\partial \alpha_i}=0, \; i=1,2,3.$ From Theorem 6.3 if we
consider the rigid body-rotor system with a control torque $u:
T^\ast Q \to W $ acting on the rotors, and $u\in W \cap
\mathbf{J}^{-1}_Q(\mu)$ is invariant under the left
$\textmd{SO}(3)$-action, and its reduced control torque $u_\mu:
\mathcal{O}_\mu \times\mathbb{R}^3\times\mathbb{R}^{3*} \to W_\mu
\times\mathbb{R}^3\times\mathbb{R}^{3*} $ is given by
$u_\mu(\Pi,\alpha,l)= \pi_\mu(u(A,\Pi,\alpha,l))=
u(A,\Pi,\alpha,l)|_{\mathcal{O}_\mu
\times\mathbb{R}^3\times\mathbb{R}^{3*} }, $ where $\pi_\mu:
\mathbf{J}_Q^{-1}(\mu) \rightarrow \mathcal{O}_\mu
\times\mathbb{R}^3\times\mathbb{R}^{3*}, $ and $W_\mu=\pi_\mu(W \cap
\mathbf{J}^{-1}_Q(\mu)). $
Thus, the equations of motion for rigid body-rotor system with the
control torque $u$ acting on the rotors are given by
\begin{equation}
\left\{\begin{aligned}\frac{\mathrm{d}\Pi}{\mathrm{d}t}&=\Pi\times\Omega\\
\frac{\mathrm{d}l}{\mathrm{d}t}&= \mbox{vlift}(u_\mu) \end{aligned}\right.\label{6.15}\end{equation}\\
where $\mbox{vlift}(u_\mu)= \mbox{vlift}(u_\mu)X_{h_\mu} \in
T(\mathcal{O}_\mu \times\mathbb{R}^3\times\mathbb{R}^{3*}). $ To sum up
the above discussion, we have the following proposition.
\begin{proposition}
The 5-tuple $(T^\ast(\textmd{SO}(3)\times
\mathbb{R}^3),\textmd{SO}(3),\omega_Q,H,u )$ is a regular point
reducible RCH system. For a point $\mu \in \mathfrak{so}^\ast(3)$,
the regular value of the momentum map $\mathbf{J}:
\textmd{SO}(3)\times \mathfrak{so}^\ast(3) \times \mathbb{R}^3
\times \mathbb{R}^{3*} \to \mathfrak{so}^\ast(3)$, the $R_P$-reduced
system is the 4-tuple $(\mathcal{O}_\mu \times \mathbb{R}^3 \times
\mathbb{R}^{3*},\tilde{\omega}_{\mathcal{O}_\mu \times \mathbb{R}^3
\times \mathbb{R}^{3*}}^{-},h_\mu,u_\mu), $ where $\mathcal{O}_\mu
\subset \mathfrak{so}^\ast(3)$ is the coadjoint orbit,
$\tilde{\omega}_{\mathcal{O}_\mu \times \mathbb{R}^3 \times
\mathbb{R}^{3*}}^{-}$ is the orbit symplectic form on $\mathcal{O}_\mu
\times \mathbb{R}^3\times \mathbb{R}^{3*} $,
$h_\mu(\Pi,\alpha,l)=H(A,\Pi,\alpha,l)|_{\mathcal{O}_\mu
\times\mathbb{R}^3\times\mathbb{R}^{3*}}$, $u_\mu(\Pi,\alpha,l)=
\pi_\mu(u(A,\Pi,\alpha,l))= u(A,\Pi,\alpha,l)|_{\mathcal{O}_\mu
\times\mathbb{R}^3\times\mathbb{R}^{3*}}$, and its equations of motion
are given by (\ref{6.15}).
\end{proposition}

\noindent {\bf (3). Heavy Top.}\\

In the following we take Lie group $G= \textmd{SE}(3)$ and state the
heavy top to be a regular point reducible Hamiltonian system, and
hence also to be a regular point reducible RCH system without the
external force and control. We know that a heavy top is by
definition a rigid body with a fixed point in $\mathbb{R}^3$ and
moving in gravitational field. Usually, exception of the singular
point, its physical phase space is $T^\ast \textmd{SO}(3)$ and the
symmetry group is $S^1$, regarded as rotations about the z-axis, the
axis of gravity, this is because gravity breaks the symmetry and the
system is no longer $\textmd{SO}(3)$ invariant. By the semidirect
product reduction theorem (See Proposition 6.5 ), we show that the
reduction of $T^\ast \textmd{SO}(3)$ by $S^1$ gives a space which is
symplectically diffeomorphic to the reduced space obtained by the
reduction of $T^\ast \textmd{SE}(3)$ by left action of
$\textmd{SE}(3)$, that is the coadjoint orbit $\mathcal{O}_{(\mu,a)}
\subset \mathfrak{se}^\ast(3)\cong T^\ast
\textmd{SE}(3)/\textmd{SE}(3)$. In fact, in this case, we can
identify the phase space $T^\ast \textmd{SO}(3)$ with the reduction
of the cotangent bundle of the special Euclidean group
$\textmd{SE}(3)=\textmd{SO}(3)\circledS \mathbb{R}^3$ by the
Euclidean translation subgroup $\mathbb{R}^3$ and identifies the
symmetry group $S^1$ with isotropy group $G_a=\{ A\in
\textmd{SO}(3)\mid Aa=a \}=S^1$, which is Abelian and
$(G_a)_{\mu_a}= G_a =S^1,\; \forall \mu_a \in \mathfrak{g}^\ast_a$,
where $a$ is a vector aligned with the direction of gravity and
where $\textmd{SO}(3)$ acts on $\mathbb{R}^3$ in the standard way.\\

Now we consider the cotangent bundle $T^\ast G =T^\ast
\textmd{SE}(3)$, and locally, $T^\ast
\textmd{SE}(3)\cong \textmd{SE}(3)\times \mathfrak{se}^\ast(3)$,
with the canonical symplectic form. Assume that Lie group
$G=\textmd{SE}(3)$ acts freely and properly by the left translation
on $\textmd{SE}(3)$, then the action of $\textmd{SE}(3)$ on the
phase space $T^\ast \textmd{SE}(3)$ is by the cotangent lift of left
translation at the identity, that is, $\Phi: \textmd{SE}(3)\times
T^\ast \textmd{SE}(3) \cong \textmd{SE}(3)\times
\textmd{SE}(3)\times \mathfrak{se}^\ast(3)\to \textmd{SE}(3)\times
\mathfrak{se}^\ast(3),$ given by
$\Phi((B,u),(A,v,\Pi,w))=(BA,v,\Pi,w)$, for any $A,B\in
\textmd{SO}(3), \; \Pi \in \mathfrak{so}^\ast(3), \; u,v,w \in
\mathbb{R}^3$, which is also free and proper. Assume that
the action is symplectic and admits an
associated $\operatorname{Ad}^\ast$-equivariant momentum map
$\mathbf{J}:T^\ast \textmd{SE}(3)\to \mathfrak{se}^\ast(3)$ for the
left $\textmd{SE}(3)$ action. If $(\Pi,w) \in \mathfrak{se}^\ast(3)$
is a regular value of $\mathbf{J}$, then the regular point reduced
space $(T^\ast
\textmd{SE}(3))_{(\Pi,w)}=\mathbf{J}^{-1}(\Pi,w)/\textmd{SE}(3)_{(\Pi,w)}$
is symplectically diffeomorphic to the coadjoint orbit
$\mathcal{O}_{(\Pi,w)} \subset \mathfrak{se}^\ast(3)$.\\

Let $I=diag(I_1,I_2,I_3)$ be the moment of inertia of the heavy top
in the body-fixed frame, which in principal body frame. Let
$\Omega=(\Omega_1,\Omega_2,\Omega_3)$ be the vector of heavy top
angular velocities computed with respect to the axes fixed in the
body and $(\Omega_1,\Omega_2,\Omega_3)\in \mathfrak{so}(3)$. Let
$\Gamma$ be the unit vector viewed by an observer moving with the
body, $m$ be that total mass of the system, $g$ be the magnitude of
the gravitational acceleration, $\chi$ be the unit vector on the
line connecting the origin $O$ to the center of mass of the system,
and $h$ be the length of this segment.\\

By the local left trivialization, we have that locally,
$T\textmd{SE}(3)\cong \textmd{SE}(3)\times \mathfrak{se}(3)$.
We consider the Lagrangian
$L(A,v,\Omega,\Gamma): \textmd{SE}(3)\times\mathfrak{se}(3)\to
\mathbb{R}$ , which is given by $$L(A,v,\Omega,\Gamma)=
\dfrac{1}{2}(I_1\Omega_1^2+I_2\Omega_2^2+I_3\Omega_3^2)-mgh\Gamma\cdot\chi,$$
where $(A,v)\in \textmd{SE}(3)$,
$\Omega=(\Omega_1,\Omega_2,\Omega_3)\in \mathfrak{so}(3)$,
$\Gamma\in\mathbb{R}^3$, and the variable $\Gamma$ is regarded as a
parameter with respect to potential energy of the heavy top.
If we introduce the conjugate angular
momentum $\Pi_i=\dfrac{\partial L}{\partial \Omega_i}=I_i\Omega_i,
\; i=1,2,3,$ and by the Legendre transformation with the parameter $\Gamma$, that is,
$FL:\textmd{SE}(3)\times\mathfrak{se}(3)\to \textmd{SE}(3)\times
\mathfrak{se}^\ast(3),\quad (A,v,\Omega,\Gamma)\to(A,v,\Pi,\Gamma)$,
where $\Pi=(\Pi_1,\Pi_2,\Pi_3)\in \mathfrak{so}^\ast(3)$, we have
the Hamiltonian $H(A,v,\Pi,\Gamma): \textmd{SE}(3)\times
\mathfrak{se}^\ast (3)\to \mathbb{R}$ given by
\begin{align*} H(A,v,\Pi,\Gamma)&=\Omega\cdot
\Pi-L(A,\Omega,,\Gamma)
=\frac{1}{2}(\frac{\Pi_1^2}{I_1}+\frac{\Pi_2^2}{I_2}+\frac{\Pi_3^2}{I_3})+mgh\Gamma\cdot\chi
\end{align*}
From the above expression of the Hamiltonian, we know that
$H(A,v,\Pi,\Gamma)$ is invariant under the left
$\textmd{SE}(3)$-action. For the case $(\Pi_0,\Gamma_0)=(\mu,a)\in
\mathfrak{se}^\ast(3)$ is a regular value of $\mathbf{J}$, we have
the reduced Hamiltonian
$h_{(\mu,a)}(\Pi,,\Gamma):\mathcal{O}_{(\mu,a)} \to \mathbb{R}$
given by
$h_{(\mu,a)}(\Pi,\Gamma)=H(A,v,\Pi,\Gamma)|_{\mathcal{O}_{(\mu,a)}}$.
From the semidirect product bracket (\ref{6.10}), we can get the
heavy top Poisson bracket on $\mathfrak{se}^\ast(3)$, that is, for
$F,K: \mathfrak{se}^\ast(3)\to \mathbb{R}, $ we have that
$$\{F,K\}_{-}(\Pi,\Gamma)=-\Pi\cdot(\nabla_\Pi F\times\nabla_\Pi
K)-\Gamma\cdot(\nabla_\Pi F\times \nabla_\Gamma K-\nabla_\Pi K\times
\nabla_\Gamma F).$$ In particular, for $F_{(\mu,a)},K_{(\mu,a)}:
\mathcal{O}_{(\mu,a)} \to \mathbb{R}$, we have that
$$\omega_{\mathcal{O}_{(\mu,a)}}^{-}(X_{F_{(\mu,a)}},
X_{K_{(\mu,a)}})=
\{F_{(\mu,a)},K_{(\mu,a)}\}_{-}|_{\mathcal{O}_{(\mu,a)}}. $$
Moreover, for reduced Hamiltonian $h_{(\mu,a)}(\Pi,\Gamma):
\mathcal{O}_{(\mu,a)} \to \mathbb{R}$, we have the Hamiltonian
vector field
$X_{h_{(\mu,a)}}(K_{(\mu,a)})=\{K_{(\mu,a)},h_{(\mu,a)}\}_{-}|_{\mathcal{O}_{(\mu,a)}
},$ and hence we have that
\begin{align*}
\frac{\mathrm{d}\Pi}{\mathrm{d}t}&=X_{h_{(\mu,a)}}(\Pi)=\{\Pi,h_{(\mu,a)}(\Pi,\Gamma)\}_{-}\\
&=-\Pi\cdot(\nabla_\Pi\Pi\times\nabla_\Pi
h_{(\mu,a)})-\Gamma\cdot(\nabla_\Pi\Pi\times\nabla_\Gamma
h_{(\mu,a)}-\nabla_\Pi
h_{(\mu,a)}\times\nabla_\Gamma\Pi)\\
&=\Pi\times\Omega-mgh\chi\times\Gamma=\Pi\times\Omega+mgh\Gamma\times\chi,
\end{align*}
\begin{align*}
\frac{\mathrm{d}\Gamma}{\mathrm{d}t}&=X_{h_{(\mu,a)}}(\Gamma)=\{\Gamma,h_{(\mu,a)}(\Pi,\Gamma)\}_{-}\\
&=-\Pi\cdot(\nabla_\Pi\Gamma\times\nabla_\Pi
h_{(\mu,a)})-\Gamma\cdot(\nabla_\Pi\Gamma\times\nabla_\Gamma
h_{(\mu,a)}-\nabla_\Pi
h_{(\mu,a)}\times\nabla_\Gamma\Gamma)\\
&=\nabla_\Gamma\Gamma\cdot(\Gamma\times\nabla_\Pi
h_{(\mu,a)})=\Gamma\times\Omega,
\end{align*}
since $\nabla_\Pi\Pi=1, \; \nabla_\Gamma \Gamma =1, \;
\nabla_\Gamma\Pi= \nabla_\Pi\Gamma =0,$ and $\nabla_\Pi
h_{(\mu,a)}=\Omega$. Thus, the equations of motion for heavy top
is given by
\begin{equation}
\left\{\begin{aligned}&\frac{\mathrm{d}\Pi}{\mathrm{d}t}=\Pi\times\Omega+mgh\Gamma\times\chi,
\\&\frac{\mathrm{d}\Gamma}{\mathrm{d}t}=\Gamma\times\Omega.
\end{aligned} \right.\label{6.16}\end{equation}
To sum up the above discussion, we have the following proposition.
\begin{proposition}
The 4-tuple $(T^\ast \textmd{SE}(3),\textmd{SE}(3),\omega_0,H)$ is a
regular point reducible Hamiltonian system. For a point $(\mu,a)\in
\mathfrak{se}^\ast(3)$, the regular value of the momentum map
$\mathbf{J}: T^\ast \textmd{SE}(3)\to \mathfrak{se}^\ast(3)$, the
$R_P$-reduced system is the 3-tuple $(\mathcal{O}_{(\mu,a)},
\omega_{\mathcal{O}_{(\mu,a)}},h_{(\mu,a)})$, where
$\mathcal{O}_{(\mu,a)} \subset \mathfrak{se}^\ast(3)$ is the
coadjoint orbit, $\omega_{\mathcal{O}_{(\mu,a)}}$ is the orbit
symplectic form on $\mathcal{O}_{(\mu,a)}$,
$h_{(\mu,a)}(\Pi,\Gamma)=H(A,v,\Pi,\Gamma)|_{\mathcal{O}_{(\mu,a)}}$,
and its equations of motion are given by (\ref{6.16}).
\end{proposition}

\noindent {\bf (4). The Heavy Top with Internal Rotors.}\\

In the following we take Lie group $G= \textmd{SE}(3), \; V=
S^1\times S^1, \; Q= G\times V $ and state the heavy top with two
pairs of symmetric internal rotors to be a regular point reducible
RCH system. We shall first describe a heavy top with two pairs of
symmetric rotors. We mount two pairs of rotors within the top so
that each pair's rotation axis is parallel to the first and the
second principal axes of the top; see Chang and Marsden
\cite{chma00}. The rotor spins under the influence of a torque $u$
acting on the rotor. The configuration space is
$Q=\textmd{SE}(3)\times V$, where $V=S^1\times S^1$, with the first
factor being the position of the heavy top and the second factor
being the angles of rotors. The corresponding phase space is the
cotangent bundle $T^*Q$ and locally,
$T^\ast Q \cong T^\ast \textmd{SE}(3)\times T^\ast V$,
where $T^\ast V = T^\ast (S^1\times S^1)\cong T^\ast \mathbb{R}^2$ locally,
with the canonical symplectic form $\omega_Q$.
By using the local left trivialization, locally,
$T^\ast \textmd{SE}(3)\cong \textmd{SE}(3)\times \mathfrak{se}^\ast(3)$
and $T^*\mathbb{R}^2 \cong \mathbb{R}^2 \times \mathbb{R}^{2*}$,
then we have that locally,
$T^*Q \cong \textmd{SE}(3)\times \mathfrak{se}^\ast(3)
\times \mathbb{R}^2 \times \mathbb{R}^{2*}$.
Assume that Lie group
$G=\textmd{SE}(3)$ acts freely and properly on $Q$ by the left
translation on the first factor $\textmd{SE}(3)$ and
the trivial action on the second factor $V$.
Then the action of $\textmd{SE}(3)$
on the phase space $T^\ast Q$ is by the cotangent lift of the left
$\textmd{SE}(3)$ action on $Q$, that is, $\Phi:
\textmd{SE}(3)\times T^\ast \textmd{SE}(3)\times T^\ast V \cong
\textmd{SE}(3)\times \textmd{SE}(3)\times \mathfrak{se}^\ast(3)
\times \mathbb{R}^2 \times \mathbb{R}^{2*}\to \textmd{SE}(3)\times
\mathfrak{se}^\ast(3)\times \mathbb{R}^2 \times \mathbb{R}^{2*},$ given
by $\Phi((B,u)((A,v),(\Pi,w),\alpha,l))=((BA,v),(\Pi,w),\alpha,l)$,
for any $A,B\in \textmd{SO}(3), \; \Pi \in \mathfrak{so}^\ast(3), \;
u,v,w \in \mathbb{R}^3, \; \alpha \in \mathbb{R}^2,
\; l \in \mathbb{R}^{2*}$, which is also
free and proper. Assume that the action is symplectic and admits an associated
$\operatorname{Ad}^\ast$-equivariant momentum map $\mathbf{J}_Q:
T^\ast Q \cong \textmd{SE}(3)\times \mathfrak{se}^\ast(3) \times
\mathbb{R}^2 \times \mathbb{R}^{2*} \to \mathfrak{se}^\ast(3)$ for the
left $\textmd{SE}(3)$ action. If $(\Pi,w) \in \mathfrak{se}^\ast(3)$
is a regular value of $\mathbf{J}_Q$, then the regular point reduced
space $(T^\ast Q)_{(\Pi,w)}=
\mathbf{J}^{-1}_Q(\Pi,w)/\textmd{SE}(3)_{(\Pi,w)}$ is symplectically
diffeomorphic to the coadjoint orbit $\mathcal{O}_{(\Pi,w)} \times
\mathbb{R}^2 \times \mathbb{R}^{2*} \subset \mathfrak{se}^\ast(3)
\times \mathbb{R}^2 \times \mathbb{R}^{2*}$.\\

Let $I=diag(I_1,I_2,I_3)$ be the moment of inertia of the heavy top
in the body-fixed frame. Let $J_i, i=1,2$ be the moments of inertia
of rotors around their rotation axes. Let $J_{ik},\; i=1,2,\;
k=1,2,3,$ be the moments of inertia of the $i$-th rotor with $i=1,2$
around the $k$-th principal axis with $k=1,2,3,$ respectively, and
denote $\bar{I}_i=I_i+J_{1i}+J_{2i}-J_{ii}, \; i=1,2$, and
$\bar{I}_3=I_3+J_{13}+J_{23}$. Let
$\Omega=(\Omega_1,\Omega_2,\Omega_3)$ be the vector of heavy top
angular velocities computed with respect to the axes fixed in the
body and $(\Omega_1,\Omega_2,\Omega_3)\in \mathfrak{so}(3)$. Let
$\theta_i,\; i=1,2,$ be the relative angles of rotors and
$\dot{\theta}=(\dot{\theta_1},\dot{\theta_2})$ the vector of rotor
relative angular velocities about the principal axes with respect to
the body fixed frame of heavy top. Let $m$ be that total mass of the
system, $g$ be the magnitude of the gravitational acceleration and
$h$ be the distance from the origin $O$ to the center of mass of the
system.\\

Now, by the local left trivialization, locally,
$T\textmd{SE}(3)\cong \textmd{SE}(3)\times \mathfrak{se}(3)$
and $T\mathbb{R}^2 \cong \mathbb{R}^2 \times \mathbb{R}^{2}$,
then we have that locally,
$TQ \cong \textmd{SE}(3)\times \mathfrak{se}(3)
\times \mathbb{R}^2 \times \mathbb{R}^{2}$.
We consider the Lagrangian
$L(A,v,\Omega,\Gamma,\theta,\dot{\theta}):
\textmd{SE}(3)\times\mathfrak{se}(3)\times\mathbb{R}^2\times\mathbb{R}^2\to
\mathbb{R}$, which is the total kinetic energy of the heavy top plus
the total kinetic energy of rotors minus potential energy of the
system, given by
$$L(A,v,\Omega,\Gamma,\theta,\dot{\theta})
=\dfrac{1}{2}[\bar{I}_1\Omega_1^2+\bar{I}_2\Omega_2^2+\bar{I}_3\Omega_3^2
+J_1(\Omega_1+\dot{\theta}_1)^2+J_2(\Omega_2+\dot{\theta}_2)^2]-mgh\Gamma\cdot\chi,$$
where $(A,v)\in \textmd{SE}(3)$, $(\Omega,\Gamma)\in
\mathfrak{se}(3)$ and $\Omega=(\Omega_1,\Omega_2,\Omega_3)\in
\mathfrak{so}(3)$,
$\theta=(\theta_1,\theta_2)\in \mathbb{R}^2$,
$\dot{\theta}=(\dot{\theta}_1,\dot{\theta}_2)\in \mathbb{R}^2$,
$\Gamma\in\mathbb{R}^3$, and the variable $\Gamma$ is regarded as a
parameter with respect to potential energy of the system. If
we introduce the conjugate angular momentum, which is given by
$$\Pi_i= \dfrac{\partial L}{\partial
\Omega_i}=\bar{I}_i\Omega_i+J_i(\Omega_i+\dot{\theta}_i),\; i=1,2,$$
$$\Pi_3=\dfrac{\partial L}{\partial
\Omega_3}=\bar{I}_3\Omega_3,\quad l_i=\dfrac{\partial L}{\partial
\dot{\theta}_i}=J_i(\Omega_i+\dot{\theta}_i),\; i=1,2,$$ and by the
Legendre transformation with the parameter $\Gamma$, that is, $FL:
\textmd{SE}(3)\times\mathfrak{se}(3)\times\mathbb{R}^2\times\mathbb{R}^2\to
\textmd{SE}(3)\times
\mathfrak{se}^\ast(3)\times\mathbb{R}^2\times\mathbb{R}^{2*},\quad
(A,v,\Omega,\Gamma,\theta,\dot{\theta})\to(A,v,\Pi,\Gamma,\theta,l),$
where $\Pi=(\Pi_1,\Pi_2,\Pi_3)\in \mathfrak{so}^\ast(3)$,
$l=(l_1,l_2)\in \mathbb{R}^{2*}$, we have the Hamiltonian
$H(A,v,\Pi,\Gamma,\theta,l):\textmd{SE}(3)\times \mathfrak{se}^\ast
(3)\times\mathbb{R}^2\times\mathbb{R}^{2*}\to \mathbb{R}$ given by
\begin{align*} &H(A,v,\Pi,\Gamma,\theta,l)=\Omega\cdot
\Pi+\dot{\theta}\cdot l-L(A,v,\Omega,\Gamma,\theta,\dot{\theta})\\
&=\frac{1}{2}[\frac{(\Pi_1-l_1)^2}{\bar{I}_1}+\frac{(\Pi_2-l_2)^2}{\bar{I}_2}
+\frac{\Pi_3^2}{\bar{I}_3}+\frac{l_1^2}{J_1}+\frac{l_2^2}{J_2}]+mgh\Gamma\cdot\chi.
\end{align*}
From the above expression of the Hamiltonian, we know that
$H(A,v,\Pi,\Gamma,\theta,l)$ is invariant under the left
$\textmd{SE}(3)$-action. For the case $(\Pi_0,\Gamma_0)=(\mu,a)\in
\mathfrak{se}^\ast(3)$ is the regular value of $\mathbf{J}_Q$, we
have the reduced Hamiltonian
$h_{(\mu,a)}(\Pi,\Gamma,\theta,l):\mathcal{O}_{(\mu,a)}\times\mathbb{R}^2
\times \mathbb{R}^{2*} (\subset \mathfrak{se}^\ast (3)\times
\mathbb{R}^2\times \mathbb{R}^{2*}) \to \mathbb{R}$ given by
$h_{(\mu,a)}(\Pi,\Gamma,\theta,l)
=H(A,v,\Pi,\Gamma,\theta,l)|_{\mathcal{O}_{(\mu,a)}\times
\mathbb{R}^2\times \mathbb{R}^{2*}}$. From (6.6), the heavy top Poisson
bracket on $\mathfrak{se}^\ast(3)$ and the Poisson bracket on
$T^\ast \mathbb{R}^2$, we can get the Poisson bracket on $T^\ast Q$,
that is, for $F,K: \mathfrak{se}^\ast(3)\times \mathbb{R}^2\times
\mathbb{R}^{2*} \to \mathbb{R}, $ we have that
\begin{align*}\{F,K\}_{-}(\Pi,\Gamma,\theta,l)=-\Pi\cdot(\nabla_\Pi
F\times\nabla_\Pi K)-\Gamma\cdot(\nabla_\Pi F\times \nabla_\Gamma
K-\nabla_\Pi K\times \nabla_\Gamma F)+
\{F,K\}_V(\theta,l).\end{align*} In particular, for
$F_{(\mu,a)},K_{(\mu,a)}: \mathcal{O}_{(\mu,a)}\times
\mathbb{R}^2\times \mathbb{R}^{2*} \to \mathbb{R}$, we have that
$$\tilde{\omega}_{\mathcal{O}_{(\mu,a)}\times \mathbb{R}^2\times
\mathbb{R}^{2*}}^{-}(X_{F_{(\mu,a)}}, X_{K_{(\mu,a)}})=
\{F_{(\mu,a)},K_{(\mu,a)}\}_{-}|_{\mathcal{O}_{(\mu,a)}\times
\mathbb{R}^2\times \mathbb{R}^{2*}}. $$ Moreover, for the reduced
Hamiltonian $h_{(\mu,a)}(\Pi,\Gamma): \mathcal{O}_{(\mu,a)}\times
\mathbb{R}^2\times \mathbb{R}^{2*} \to \mathbb{R}$, we have the
Hamiltonian vector field
$X_{h_{(\mu,a)}}(K_{(\mu,a)})=\{K_{(\mu,a)},h_{(\mu,a)}\}_{-}|_{\mathcal{O}_{(\mu,a)}
\times \mathbb{R}^2\times \mathbb{R}^{2*}},$ and hence we have that
\begin{align*}
\frac{\mathrm{d}\Pi}{\mathrm{d}t}&=X_{h_{(\mu,a)}}(\Pi)(\Pi,\Gamma,\theta,l)
=\{\Pi,h_{(\mu,a)}\}_{-}(\Pi,\Gamma,\theta,l)=-\Pi\cdot(\nabla_\Pi\Pi\times\nabla_\Pi
h_{(\mu,a)})\\
& -\Gamma\cdot(\nabla_\Pi\Pi\times\nabla_\Gamma
h_{(\mu,a)}-\nabla_\Pi h_{(\mu,a)}\times\nabla_\Gamma\Pi)+
\sum_{i=1}^2(\frac{\partial \Pi}{\partial \theta_i} \frac{\partial
h_{(\mu,a)}}{\partial l_i}- \frac{\partial h_{(\mu,a)}}{\partial
\theta_i}\frac{\partial \Pi}{\partial l_i})\\
&=\Pi\times\Omega-mgh\chi\times\Gamma=\Pi\times\Omega+mgh\Gamma\times\chi,
\end{align*}
\begin{align*}
\frac{\mathrm{d}\Gamma}{\mathrm{d}t}&=X_{h_{(\mu,a)}}(\Gamma)(\Pi,\Gamma,\theta,l)
=\{\Gamma,h_{(\mu,a)}\}_{-}(\Pi,\Gamma,\theta,l)=-\Pi\cdot(\nabla_\Pi\Gamma\times\nabla_\Pi
h_{(\mu,a)})\\
& -\Gamma\cdot(\nabla_\Pi\Gamma\times\nabla_\Gamma
h_{(\mu,a)}-\nabla_\Pi h_{(\mu,a)}\times\nabla_\Gamma\Gamma)+
\sum_{i=1}^2(\frac{\partial \Gamma}{\partial \theta_i}
\frac{\partial h_{(\mu,a)}}{\partial l_i}- \frac{\partial
h_{(\mu,a)}}{\partial
\theta_i}\frac{\partial \Gamma}{\partial l_i})\\
&=\nabla_\Gamma\Gamma\cdot(\Gamma\times\nabla_\Pi
h_{(\mu,a)})=\Gamma\times\Omega,
\end{align*}
since $\nabla_\Pi\Pi=1, \; \nabla_\Gamma \Gamma =1, \;
\nabla_\Gamma\Pi= \nabla_\Pi\Gamma =0$, $\nabla_\Pi
h_{(\mu,a)}=\Omega$, and $\frac{\partial \Pi}{\partial \theta_i}=
\frac{\partial \Gamma}{\partial \theta_i}= \frac{\partial
h_{(\mu,a)}}{\partial \theta_i}=0, \; i=1,2$. From Theorem 6.3 if we
consider the heavy top-rotor system with a control torque $u: T^\ast
Q \to W $ acting on the rotors, and $u\in W \cap
\mathbf{J}^{-1}_Q((\mu,a))$ is invariant under the left
$\textmd{SE}(3)$-action, and its reduced control torque
$u_{(\mu,a)}: \mathcal{O}_{(\mu,a)}
\times\mathbb{R}^2\times\mathbb{R}^{2*} \to W_{(\mu,a)}
\times\mathbb{R}^2\times\mathbb{R}^{2*} $ is given by
$u_{(\mu,a)}(\Pi,\Gamma,\theta,l)=\pi_{(\mu,a)}(u(A,v,\Pi,\Gamma,\theta,l))=
u(A,v,\Pi,\Gamma,\theta,l)|_{\mathcal{O}_{(\mu,a)}
\times\mathbb{R}^2\times\mathbb{R}^{2*}}, $ where $\pi_{(\mu,a)}:
\mathbf{J}_Q^{-1}((\mu,a)) \rightarrow \mathcal{O}_{(\mu,a)}
\times\mathbb{R}^2\times\mathbb{R}^{2*}, $ and
$W_{(\mu,a)}=\pi_{(\mu,a)}(W \cap \mathbf{J}^{-1}_Q((\mu,a))). $
Thus, the equations of motion for
heavy top-rotor system with the control torque $u$ acting on the
rotors are given by
\begin{equation}
\left\{\begin{aligned}&\frac{\mathrm{d}\Pi}{\mathrm{d}t}=\Pi\times\Omega+mgh\Gamma\times\chi,
\\&\frac{\mathrm{d}\Gamma}{\mathrm{d}t}=\Gamma\times\Omega, \\
&\frac{\mathrm{d}l}{\mathrm{d}t}= \mbox{vlift}(u_{(\mu,a)}). \end{aligned}\right.\label{6.17}\end{equation}\\
where $\mbox{vlift}(u_{(\mu,a)})=
\mbox{vlift}(u_{(\mu,a)})X_{h_{(\mu,a)}} \in T(\mathcal{O}_{(\mu,a)}
\times\mathbb{R}^2\times\mathbb{R}^{2*}).$ To sum up the above
discussion, we have the following proposition.
\begin{proposition}
The 5-tuple $(T^\ast(\textmd{SE}(3)\times
\mathbb{R}^2),\textmd{SE}(3),\omega_Q,H,u)$ is a regular point
reducible RCH system. For a point $(\mu,a)\in
\mathfrak{se}^\ast(3)$, the regular value of the momentum map\\
$\mathbf{J}: \textmd{SE}(3)\times \mathfrak{se}^\ast(3)\times
\mathbb{R}^2\times \mathbb{R}^{2*} \to \mathfrak{se}^\ast(3)$, the
$R_P$-reduced system is the 4-tuple $(\mathcal{O}_{(\mu,a)}\times
\mathbb{R}^2\times \mathbb{R}^{2*},\\
\tilde{\omega}^{-}_{\mathcal{O}_{(\mu,a)}\times \mathbb{R}^2\times
\mathbb{R}^{2*}},\; h_{(\mu,a)},\; u_{(\mu,a)})$, where
$\mathcal{O}_{(\mu,a)} \subset \mathfrak{se}^\ast(3)$ is the
coadjoint orbit, $\tilde{\omega}^{-}_{\mathcal{O}_{(\mu,a)}\times
\mathbb{R}^2\times \mathbb{R}^{2*}}$ is the orbit symplectic form on
$\mathcal{O}_{(\mu,a)} \times \mathbb{R}^2\times \mathbb{R}^{2*} $,
$h_{(\mu,a)}(\Pi,\Gamma,\theta,l)=H(A,v,\Pi,\Gamma,\theta,l)|_{\mathcal{O}_{(\mu,a)}\times
\mathbb{R}^2\times \mathbb{R}^{2*}}, $ and
$u_{(\mu,a)}(\Pi,\Gamma,\theta,l)=\pi_{(\mu,a)}(u(A,v,\Pi,\Gamma,\theta,l))=
u(A,v,\Pi,\Gamma,\theta,l)|_{\mathcal{O}_{(\mu,a)}
\times\mathbb{R}^2\times\mathbb{R}^{2*}}, $ and its equations of
motion are given by (\ref{6.17}).
\end{proposition}

\noindent {\bf (5). Regular Controlled Hamiltonian Equivalence.}\\

In the following we shall state the RCH-equivalences of the rigid
body with external force torque and that with internal rotors, as
well as the heavy top and that with internal rotors. In fact, we can
choose the feedback control law such that the equivalent RCH systems
produce the same equations of motion (up to a diffeomorphism ).\\

At first, we consider the RCH-equivalence between the rigid body
with external force torque and that with internal rotors. Now let
us choose the feedback control laws such that the closed-loop
systems are Hamiltonian and retain the symmetry. If we choose the
feedback control law $u$, such that
$\mbox{vlift}(u_\mu)=p\times\Omega$, where $p$ is a constant vector,
from the equations (\ref{6.14}) of motion for the rigid body with
the $\textmd{SO}(3)$-invariant external force torque $u$, we have
that
\begin{equation}   \dfrac{\mathrm{d}
\Pi}{\mathrm{d} t}=\Pi\times\Omega+ p\times\Omega .\label{6.18}
\end{equation}
On the other hand, for the rigid body with internal rotors, we
choose the feedback control law $u$, such that
$\mbox{vlift}(u_\mu)=k(\Pi\times\Omega)$, where $k$ is a gain
parameter. From the equations (\ref{6.15}) of motion for the rigid
body with internal rotors, we have that $\frac{\mathrm{d}
l}{\mathrm{d}t}=\mbox{vlift}(u_\mu)=k\frac{\mathrm{d}\Pi}{\mathrm{d}t}$,
and by solving the integrable equation, we get that $l-k\Pi=p$,
where $p$ is a constant vector. Assuming that
$N=\Pi-l=\Pi-k\Pi-p=(1-k)\Pi-p$, then we have that
\begin{equation}\frac{\mathrm{d}N}{\mathrm{d}t}=
\frac{\mathrm{d}\Pi}{\mathrm{d}t}- \frac{\mathrm{d}l}{\mathrm{d}t} =
(1-k)\Pi\times\Omega = N\times\Omega + p\times\Omega.\label{6.19}
\end{equation}
By comparing (\ref{6.18}) and (\ref{6.19}) we know that the rigid
body with external force torque and that with internal rotors are
RCH-equivalent by a diffeomorphism $\varphi: \mathfrak{so}^\ast (3)
\rightarrow \mathfrak{so}^\ast (3), \Pi \rightarrow N $. In
particular, if we take that
$\mbox{vlift}(u_\mu)=(u_{\mu1},u_{\mu2},u_{\mu3})=(0,0,-\varepsilon
\frac{I_1-I_2}{I_1I_2}\Pi_1\Pi_2)\in \mathbb{R}^3$, we recover the
result in Bloch et al. \cite{blkrmaal92}, also see Marsden
\cite{ma92}.\\

Next, we consider the RCH-equivalence between the rigid body with
internal rotors and heavy top. If assuming that $N=\Pi +\Gamma$,
from the equations (\ref{6.16}) of motion for the heavy top, we have
that
$$ \frac{\mathrm{d}N}{\mathrm{d}t}= N\times\Omega + mgh\Gamma\times\chi
= N\times\Omega - mgh\chi\times\Gamma
$$
Thus, take that $\Gamma=\lambda\Omega$ and $p=-mgh\lambda\chi$,
where $\lambda$ is a constant, then
\begin{equation}\frac{\mathrm{d}N}{\mathrm{d}t}= N\times\Omega +
p\times\Omega.\label{6.20}
\end{equation}
In this case, by comparing (\ref{6.19}) and (\ref{6.20}) we know
that the heavy top and the rigid body with internal rotors are
RCH-equivalent. In the same way, from (\ref{6.18}) we know that the
rigid body with the external force torques and the heavy top are
also RCH-equivalent. Also see Holm and Marsden \cite{homa91}.\\

At last, we consider the RCH-equivalence between the rigid body with
internal rotors and heavy top with internal rotors. For the heavy
top with internal rotors, we choose the feedback control law $u$,
such that $\mbox{vlift}(u_{(\mu,a)})=k(\Gamma\times\Omega)$, where
$k$ is a gain parameter. From the equations (\ref{6.17}) of motion
for the heavy top with internal rotors, we have that
$\frac{\mathrm{d} \bar{l}}{\mathrm{d}t}= \mbox{vlift}(u_{(\mu,a)})
=k \frac{\mathrm{d}\Gamma}{\mathrm{d}t}$, where
$\bar{l}=(l_1,l_2,0)$, and by solving the integrable equation, we
get that $\bar{l}-k\Gamma=p_0$, where $p_0$ is a constant vector.
Assuming that $N=\Pi+\Gamma-\bar{l}=\Pi+(1-k)\Gamma-p_0$, then we
have that
$$\frac{\mathrm{d}N}{\mathrm{d}t}=
\frac{\mathrm{d}\Pi}{\mathrm{d}t}+
\frac{\mathrm{d}\Gamma}{\mathrm{d}t}-
\frac{\mathrm{d}\bar{l}}{\mathrm{d}t} = \Pi\times\Omega
+(1-k)\Gamma\times\Omega - mgh\chi\times\Gamma = N\times\Omega +
p_0\times\Omega - mgh\chi\times\Gamma.
$$
Thus, take that $\Gamma=\lambda\Omega$ and $p=p_0-mgh\lambda\chi$,
where $\lambda$ is a constant, then
\begin{equation}\frac{\mathrm{d}N}{\mathrm{d}t}= N\times\Omega +
p\times\Omega.\label{6.21}
\end{equation}
In this case, by comparing (\ref{6.19}) and (\ref{6.21}) we know
that the rigid body with internal rotors and the heavy top with
internal rotors are RCH-equivalent.

To sum up, we have the following theorem.
\begin{theorem} As two $R_P$-reduced RCH systems,

\noindent({\bf i}) the rigid body with external force torque and
that with internal rotors are RCH-equivalent;

\noindent({\bf ii}) the rigid body with internal rotors (or external
force torque) and the heavy top are RCH-equivalent;

\noindent({\bf iii}) the rigid body with internal rotors and the
heavy top with internal rotors are RCH-equivalent.
\end{theorem}

\subsection{Port Hamiltonian System with a Symplectic
Structure}

In order to understand well the abstract definition of RCH system
and the RCH-equivalence, in this subsection we will describe the RCH
system and RCH-equivalence from the viewpoint of port Hamiltonian
system with a symplectic structure. Recently years, the study of
stability analysis and control of port Hamiltonian systems and their
applications have become more and more important, and there have
been a lot of beautiful results; see Dalsmo and van der Schaft
\cite{davds99}, van der Schaft \cite{vds00,vds06}. To describe the
RCH systems well from the viewpoint of port Hamiltonian system, in
the following we first give some relevant definitions and basic
facts about the port Hamiltonian systems.

\begin{definition}
Let $(T^\ast Q, \omega)$ be a symplectic manifold and $\omega$ be
the canonical symplectic form on $T^\ast Q$. Assume that $H: T^\ast
Q \rightarrow \mathbb{R}$ is a Hamiltonian, and there exists a
subset $U\subset T^\ast Q$ and a vector field $X_H \in TT^\ast Q$ on
$T^\ast Q$ such that $i_{X_H}\omega(z)= \mathbf{d}H(z), \; \forall z
\in U $, then the triple $(T^\ast Q, \omega, H)$ is a Hamiltonian
system defined on the set $U$. Assume that $V \subset T^\ast Q$ is a
subset of $T^\ast Q$, and $P=(Y,\alpha)$, where for any $z \in V$,
$Y(z) \in T_zT^\ast Q$ and $\alpha(z) \in T^\ast _z T^\ast Q$. If $U
\cap V\neq \emptyset$, and $i_{(X_H+Y)}\omega(z)= (\mathbf{d}H
+\alpha)(z), \; \forall z \in U \cap V$, then $P=(Y,\alpha)$ is
called a port of the Hamiltonian system $(T^\ast Q, \omega, H)$
defined on the set $U$. The 4-tuple $(T^\ast Q, \omega, H, P)$ is
called a port Hamiltonian system.
\end{definition}
For the port Hamiltonian system $(T^\ast Q, \omega, H, P)$, since
$i_{X_H}\omega(z) = \mathbf{d}H(z), \; \forall z \in U$, from
$i_{(X_H+Y)}\omega(z)\\ = (\mathbf{d}H +\alpha)(z), \; \forall z \in
U \cap V$, we have that $ i_{X_H}\omega(z) +i_Y\omega(z)=
\mathbf{d}H(z) + \alpha(z).$ Thus, we can get the port balance
condition that $P=(Y,\alpha)$ is a port of the Hamiltonian system
$(T^\ast Q, \omega, H )$ as follows
\begin{equation}
i_Y\omega(z)= \alpha(z), \;\;\;\; \forall z \in U \cap
V.\label{6.22}
\end{equation}
In particular, for $U=V=T^\ast Q$, from the port balance condition
(\ref{6.22}) we know that $P=(X_H,\mathbf{d}H)$ is a trivial port of
the Hamiltonian system $(T^\ast Q, \omega, H )$.\\

Assume that $(T^\ast Q,\omega,H,F,u)$ is a RCH system with a control
law $u$. We can take that $Y= \textnormal{vlift}(F+u)\in TT^\ast Q$,
from the port balance condition (\ref{6.22}) we take that $\alpha=
i_Y \omega \in T^\ast T^\ast Q$, then $P=(Y,\alpha)$ is a
force-controlled port of the Hamiltonian system $(T^\ast Q, \omega,
H )$, and $(T^\ast Q, \omega, H, P)$ is a port Hamiltonian system
with a symplectic structure. Thus, we have the following
proposition.
\begin{proposition}
Any RCH system $(T^\ast Q,\omega,H,F,u)$ with control law $u$, is a
port Hamiltonian system with symplectic structure.
\end{proposition}
If we consider the canonical coordinates $z=(q,p)$ of the phase
space $T^\ast Q$, then $X_H= (\dot{q},\dot{p})$, and the local
expression of the RCH system is given by
\begin{equation}
\dot{q}= \frac{\partial H}{\partial p}(q,p), \;\;\;\;\;\; \dot{p}=
-\frac{\partial H}{\partial q}(q,p)
+\textnormal{vlift}(F+u)(q,p).\label{6.23}
\end{equation}
We can derive the energy balance condition, that is,
\begin{equation}
\frac{dH}{dt}=(\frac{\partial H}{\partial q})^T(q,p)\dot{q} +
(\frac{\partial H}{\partial p})^T(q,p)\dot{p}=(\frac{\partial
H}{\partial p})^T\textnormal{vlift}(F+u)(q,p)=
\dot{q}^T\textnormal{vlift}(F+u)(q,p),\label{6.24}
\end{equation}
which expresses that the increase in energy of the system is equal
to the supplied work (that is, conservation of energy). This
motivates to define the output of the system as $e=\dot{q}$, which
is considered as the vector of generalized velocities, and the
local expression of the port controlled Hamiltonian system is
given by
\begin{equation}
\dot{q}= \frac{\partial H}{\partial p}(q,p), \;\;\;\; \dot{p}=
-\frac{\partial H}{\partial q}(q,p) +B(q)f, \;\;\;\;
e=B^T(q)\dot{q}.\label{6.25}
\end{equation}
where $\textnormal{vlift}(F+u)=B(q)f$, and $f$ is an input of
system;
see van der Schaft \cite{vds00,vds06}.\\

In the following we shall state the relationships between
RCH-equivalence of RCH systems and the equivalence of port
Hamiltonian systems. We first give the definitions of equivalence of
Hamiltonian systems,  port-equivalence of port Hamiltonian systems
and equivalence of port Hamiltonian systems as follows. Assume that
$(T^\ast Q_i,\omega_i), \; i=1,2,$ are two symplectic manifolds, and
$\psi: T^\ast Q_1 \rightarrow T^\ast Q_2$ is a symplectic
diffeomorphism. Let $T\psi: TT^\ast Q_1 \rightarrow TT^\ast Q_2$ be
the tangent map of $\psi: T^\ast Q_1 \rightarrow T^\ast Q_2$, and
$\psi_\ast = (\psi^{-1})^\ast: T^\ast T^\ast Q_1 \rightarrow T^\ast
T^\ast Q_2$ be the cotangent map of $\psi^{-1}: T^\ast Q_2
\rightarrow T^\ast Q_1$. Then we can describe the equivalence of the
Hamiltonian systems as follows.
\begin{definition}
Assume that $(T^\ast Q_i, \omega_i, H_i ), \; i=1,2,$ are two
Hamiltonian systems. We say them to be equivalent, if there exists a
symplectic diffeomorphism $\psi: T^\ast Q_1 \rightarrow T^\ast Q_2$,
such that $T\psi(X_{H_1})= X_{H_2}\cdot \psi, \;
\psi_{\ast}(\mathbf{d}H_1)= \mathbf{d}H_2 \cdot \psi $, where
$i_{X_{H_i}}\omega = \mathbf{d}H_i,  \; i=1,2.$
\end{definition}
Moreover, we can describe the port-equivalence of port Hamiltonian
systems and the equivalence of port Hamiltonian systems as follows.
\begin{definition}
Assume that $(T^\ast Q_i, \omega_i, H_i, P_i ), \; i=1,2,$ are two
port Hamiltonian systems. We say them to be port-equivalent, if
there exists a diffeomorphism $\psi: T^\ast Q_1 \rightarrow T^\ast
Q_2$, such that $T\psi(Y_1)= Y_2 \cdot \psi, \;
\psi_{\ast}(\alpha_1)= \alpha_2 \cdot \psi $, where
$P_i=(Y_i,\alpha_i)$, and for any $z_i \in V_i (\subset T^\ast
Q_i)$, $Y_i(z_i) \in T_{z_i}T^\ast Q_i$ and $\alpha_i(z_i) \in
T^\ast _{z_i} T^\ast Q_i$, $i=1,2.$ Furthermore, we say two port
Hamiltonian systems $(T^\ast Q_i, \omega_i, H_i, P_i ), \; i=1,2,$
to be equivalent, if there exists a diffeomorphism $\psi: T^\ast Q_1
\rightarrow T^\ast Q_2$, such that not only two Hamiltonian systems
$(T^\ast Q_i, \omega_i, H_i ), \; i=1,2,$ are equivalent, but also
their ports are equivalent.
\end{definition}
Thus, we can obtain the following theorem.

\begin{theorem}
\noindent({\bf i}) If two RCH systems $(T^\ast
Q_i,\omega_i,H_i,F_i,W_i),\; i=1,2,$ are RCH-equivalent and their
associated Hamiltonian systems $(T^\ast Q_i,\omega_i,H_i)$, $i=1,2,$
are also equivalent, then they must be equivalent for port
Hamiltonian systems.

\noindent({\bf ii}) If two RCH systems $(T^\ast
Q_i,\omega_i,H_i,F_i,W_i),\; i=1,2,$ are RCH-equivalent, but the
associated Hamiltonian systems $(T^\ast Q_i,\omega_i,H_i),\; i=1,2,$
are not equivalent, then we can choose the control law $u_i$, such
that they are port-equivalent for port Hamiltonian systems.
\end{theorem}
{\bf Proof.} ({\bf i}) In fact, assume that two RCH systems $(T^\ast
Q_i,\omega_i,H_i,F_i,W_i), \; i=1,2,$ are RCH-equivalent, then there
exists a diffeomorphism $\varphi: Q_1\rightarrow Q_2$, such that
$\varphi^\ast: T^\ast Q_2\rightarrow T^\ast Q_1$ is symplectic, and
from Theorem 3.3 there exist two control laws $u_i: T^\ast Q_i
\rightarrow W_i, \; i=1,2, $ such that the two associated
closed-loop systems produce the same equations of motion, that is,
$X_{(T^\ast Q_1,\omega_1,H_1,F_1,u_1)}\cdot \varphi^\ast =
T\varphi^\ast
 X_{(T^\ast Q_2,\omega_2,H_2,F_2,u_2)}$. If
the associated Hamiltonian systems $(T^\ast Q_i,\omega_i,H_i),$
$i=1,2$ are also equivalent,  from $\varphi_\ast =
(\varphi^{-1})^\ast: T^\ast Q_1\rightarrow T^\ast Q_2$ is
symplectic, and $T\varphi_\ast(X_{H_1})= X_{H_2}\cdot \varphi_\ast$,
and $X_{H_i}= (\mathbf{d}H_i)^\sharp, \; i=1,2, $ we have that
$T\varphi^\ast(\mathbf{d}H_2)^\sharp = ( \mathbf{d}H_1)^\sharp \cdot
\varphi^\ast$. Note that $X_{(T^\ast Q_i,\omega_i,H_i,F_i,u_i)}
=(\mathbf{d}H_i)^\sharp+\textnormal{vlift}(F_i)+\textnormal{vlift}(u_i),
\; i=1,2, $ then,
$T\varphi^\ast(\textnormal{vlift}(F_2)+\textnormal{vlift}(u_2))=
(\textnormal{vlift}(F_1)+\textnormal{vlift}(u_1))\cdot
\varphi^\ast$. We can first take that $Y_i=
\textnormal{vlift}(F_i+u_i)\in TT^\ast Q_i, \; i=1,2, $ then we have
that $T\varphi^\ast (Y_2)= Y_1\cdot \varphi^\ast$, and hence
$T\varphi_\ast (Y_1)= Y_2\cdot \varphi_\ast$. Then we take that
$\alpha_i= i_{Y_i} \omega_i \in T^\ast T^\ast Q_i, \; i=1,2. $ Since
the map $(\varphi_\ast)_\ast= (\varphi_\ast^{-1})^\ast : T^\ast
T^\ast Q_1 \rightarrow T^\ast T^\ast Q_2$, such that
$(\varphi_\ast)_\ast (i_{Y_1} \omega_1)= i_{T\varphi_\ast
(Y_1)}(\varphi_\ast)_\ast(\omega_1)=i_{Y_2}\omega_2 \cdot
\varphi_\ast$, we have that $(\varphi_\ast)_\ast (\alpha_1)=
\alpha_2 \cdot \varphi_\ast .$ Thus, the ports $P_i=(Y_i,\alpha_i)$,
satisfying $T\varphi_\ast (Y_1)= Y_2\cdot \varphi_\ast$, and
$(\varphi_\ast)_\ast (\alpha_1)= \alpha_2 \cdot \varphi_\ast$, are
equivalent, and hence the port Hamiltonian systems $(T^\ast Q_i,
\omega_i, H_i, P_i ), \; i=1,2,$ are equivalent.\\

({\bf ii}) Assume that two RCH systems $(T^\ast
Q_i,\omega_i,H_i,F_i,W_i),\; i=1,2,$ are RCH-equivalent, but the
associated Hamiltonian systems $(T^\ast Q_i,\omega_i,H_i),\; i=1,2,$
are not equivalent, from Theorem 3.3 we can choose the control law
$u_i: T^\ast Q_i \rightarrow W_i, \; i=1,2, $ such that
$T(\varphi^\ast)\cdot X_{(T^\ast Q_2,\omega_2,H_2,F_2,u_2)} =
X_{(T^\ast Q_1,\omega_1,H_1,F_1,u_1)}\cdot \varphi^\ast$, and hence
$T(\varphi_\ast)\cdot X_{(T^\ast Q_1,\omega_1,H_1,F_1,u_1)} =
X_{(T^\ast Q_2,\omega_2,H_2,F_2,u_2)}\cdot \varphi_\ast$. We can
take that $Y_i = X_{(T^\ast Q_i,\omega_i,H_i,F_i,u_i)} =
(\mathbf{d}H_i)^\sharp
+\textnormal{vlift}(F_i)+\textnormal{vlift}(u_i)\in TT^\ast Q_i, $
and $\alpha_i= i_{Y_i} \omega_i \in T^\ast T^\ast Q_i, $ $ i=1,2 $.
Then the ports $P_i=(Y_i,\alpha_i), \; i=1,2,$ satisfy that
$T\varphi_\ast (Y_1)= Y_2\cdot \varphi_\ast$, and
$(\varphi_\ast)_\ast (\alpha_1)= \alpha_2 \cdot \varphi_\ast$, and
hence the port Hamiltonian systems $(T^\ast Q_i, \omega_i, H_i, P_i
),$ \; $i=1,2,$ are port-equivalent. \hskip 1cm $\blacksquare$\\

The theory of controlled mechanical system is a very important subject.
In this paper, we study the regular reduction theory of a controlled
Hamiltonian system with the symplectic structure and symmetry. It
is a natural problem what and how we could do, if we define a
controlled Hamiltonian system on the cotangent bundle $T^*Q$ by
using a Poisson structure, and if symplectic reduction procedure
does not work or is not efficient enough. Wang and Zhang in
\cite{wazh12} study the optimal reduction theory of a controlled
Hamiltonian system with Poisson structure and symmetry by using the
optimal momentum map. Moreover, Wang in \cite{wa13d, wa13f} study
the Hamilton-Jacobi theory of RCH system and its a variety of
reduced systems, and describe the relationship between the
RCH-equivalence for RCH systems and the solutions of corresponding
Hamilton-Jacobi equations, and apply to give explicitly the motion
equation and Hamilton-Jacobi equation of the reduced spacecraft-rotor
system on a symplectic leaf by calculation in detail, which show the
effect on controls in regular symplectic reduction and
Hamilton-Jacobi theory.\\

\noindent{\bf Acknowledgments:} The authors would like to thank the
referees for their careful review reports and comments, which are
helpful for us to revise well the manuscript. Especially grateful to
Professor Tudor S. Ratiu and MS. Wendy McKay for their support and
help. H. Wang's research was partially supported by the Natural
Science Foundation of Tianjin (05YFJMJC01200) and the Key Laboratory
of Pure Mathematics and Combinatorics, Ministry of Education, China.


\begin{thebibliography}{99}

\bibitem{abma78}
R. Abraham and J.E. Marsden, Foundations of Mechanics, 2$^{nd}$
edition, Addison-Wesley, 1978.
\bibitem{abmara88}
R. Abraham, J.E. Marsden and T.S. Ratiu, Manifolds, Tensor Analysis
and Applications,  Applied Mathematical Science, \textbf{75},
Springer-Verlag, 1988.
\bibitem{ar89}
V.I. Arnold, Mathematical Methods of Classical Mechanics, 2$^{nd}$
edition, Graduate Texts in Mathematics, \textbf{60},
Springer-Verlag, 1989.
\bibitem{blorvds02}
G. Blankenstein, R. Ortega and A.J. van der Schaft, The matching
conditions of controlled Lagrangians and IDA-passivity based
control,  Inter. J. Contr. \textbf{75}(9)(2002), 645--665.
\bibitem{blchlema01}
A.M. Bloch, D.E. Chang, N.E. Leonard and J.E. Marsden, Controlled
Lagrangians and the stabilization of mechanical systems II:
potential shaping,  IEEE Trans, Automatic Control, \textbf{46}
(2001), 1556--1571.
\bibitem{blkrmaal92}
A.M. Bloch, P.S. Krishnaprasad, J.E. Marsden, and G. S\'{a}nchez de
Alvarez, Stabilization of rigid body dynamics by internal and
external torques, Automatica, \textbf{28}(4)(1992), 745--756.
\bibitem{blle02}
A.M. Bloch and N.E. Leonard, Symmetries, conservation laws, and
control, In ``Geometry, Mechanics and Dynamics, Volume in Honour of
the 60th Birthday of J.E. Marsden" (eds. P.Newton, P.Holmes and A.
Weinstein), Springer, New York, 2002.
\bibitem{bllema00}
A.M. Bloch, N.E. Leonard and J.E. Marsden, Controlled Lagrangian and
the stabilization of mechanical systems I: the first matching
theorem,  IEEE Trans, Automatic Control, \textbf{45}(2000),
2253--2270.
\bibitem{chbllemawo02}
D.E. Chang, A.M. Bloch, N.E. Leonard, J.E. Marsden and C.A. Woolsey,
The equivalence of controlled Lagrangian and controlled Hamiltonian
systems,  ESAIM Control, Optimisation and Calculus of Variations,
\textbf{8}(2002), 393--422.
\bibitem{chma00}
D.E. Chang and J.E. Marsden, Asymptotic stabilization of the heavy
top using controlled Lagrangians, Proc. CDC, \textbf{39}(2000),
269--273.
\bibitem{chma04}
D.E. Chang and J.E. Marsden, Reduction of controlled Lagrangian and
Hamiltonian systems with symmetry, SIAM J. Control Optimization,
\textbf{43}(1)(2004), 277--300.
\bibitem{crvds87}
P.E. Crouch and A.J. van der Schaft, Variational and Hamiltonian
Control Systems, Lecture Notes in Control and Information Sciences,
Springer-Verlag, Berlin, 1987.
\bibitem{davds99}
M. Dalsmo and A.J. van der Schaft, On representations and
integrability of mathematical structures in energy-conserving
physical systems,  SIAM J. Contr. Optim. \textbf{37}(1)(1999),
54--91.
\bibitem{golasnwe83}
M.J. Gotay, R. Lashof, J. $\acute{S}$niatucki and A. Weinstein,
Closed forms on symplectic fiber bundles,  Comm. Math. Helv.
\textbf{58}(1983), 617--621.
\bibitem{homa91}
D.D. Holm and J.E. Marsden, The rotor and the pendulum, In
``Symplectic Geometry and Mathematical Physics" (eds. P. Donato, C.
Duval, J. Elhadad and G.M. Tuynman), Birkh$\ddot{a}$user, (1991),
189--203.
\bibitem{kono63}
S. Kobayashi and K. Nomizu, Foundations of Differential Geometry,
Vol.1 and Vol.2, Wiley, 1963.
\bibitem{krma87}
P.S. Krishnaprasad and J.E. Marsden, Hamiltonian structure and
stability for rigid bodies with flexible attachments,  Arch. Rat.
Mech. An. \textbf{98}(1987), 137--158.
\bibitem{lima87}
P. Libermann and C.M. Marle, Symplectic Geometry and Analytical
Mechanics, Kluwer Academic Publishers, 1987.
\bibitem{lema97}
N.E. Leonard and J.E. Marsden, Stability and drift of underwater
vehicle dynamics: mechanical systems with rigid motion symmetry,
Physica D \textbf{105}(1997), 130--162.
\bibitem{ma92}
J.E. Marsden, Lectures on Mechanics, London Mathematical Society
Lecture Notes Series, \textbf{174}, Cambridge University Press,
1992.
\bibitem{mamiorpera07}
J.E. Marsden, G. Misiolek, J.P. Ortega, M. Perlmutter and T.S.
Ratiu, Hamiltonian Reduction by Stages, Lecture Notes in
Mathematics, \textbf{1913}, Springer, 2007.
\bibitem{mara99}
J.E. Marsden and T.S. Ratiu, Introduction to Mechanics and Symmetry,
 2$^{nd}$ edition, Texts in Applied Mathematics, \textbf{17},
Springer-Verlag, 1999.
\bibitem{mawe74}
J.E. Marsden and A. Weinstein, Reduction of symplectic manifolds
with symmetry, Rep. Math. Phys., \textbf{5}(1)(1974), 121--130.
\bibitem{me73}
K.R. Meyer, Symmetries and integrals in mechanics, In ``Dynamical
Systems" (eds. M. Peixoto), Academic Press, (1973), 259--273.
\bibitem{nivds90}
H. Nijmeijer and A.J. van der Schaft, Nonlinear Dynamical Control
Systems, Springer-Verlag, New York, 1990.
\bibitem{orra04}
J.P. Ortega and T.S. Ratiu, Momentum Maps and Hamiltonian Reduction,
Progress in Mathematics, \textbf{222}, Birkh\"{a}user, 2004.
\bibitem{vds82}
A.J. van der Schaft, Hamiltonian dynamics with external forces and
observations, Mathematical Systems Theory, \textbf{15}(1982),
145--168.
\bibitem{vds86}
A.J. van der Schaft, Stabilization of Hamiltonian systems, Nonlinear
Analysis, Theory, Methods and Applications, \textbf{10}(1986),
1021--1035.
\bibitem{vdsma95}
A.J. van der Schaft and B.M. Maschke, The Hamiltonian formulation of
energy conserving physical systems with external ports,  Archiv
F$\ddot{u}$r Elektronik und $\ddot{U}$bertragungstechnik,
\textbf{49}(1995), 362--371.
\bibitem{vds00}
A.J. van der Schaft, $L_2$-Gain and Passivity Techniques in
Nonlinear Control, 2$^{nd}$ revised and enlarged edition, Comm.
Control Engrg. Ser., Springer-Verlag, London, 2000.
\bibitem{vds06}
A.J. van der Schaft, Port-Hamiltonian systems: an introductory
survey, Proceedings of the International Congress of Mathematicians,
Madrid, Spain, (2006), 1339--1365.
\bibitem{wa13d}
H. Wang, Hamilton-Jacobi theorems for a regular controlled Hamiltonian
system and its reduced systems, (arXiv: 1305.3457 ).
\bibitem{wa13f}
H. Wang, Symmetric reduction and Hamilton-Jacobi equation of rigid
spacecraft with a rotor, J. Geom. Symm. Phys., 32(2013), (arXiv:
1307.1606, a revised version ).
\bibitem{wazh12}
H. Wang and Z. X. Zhang, Optimal reduction of controlled Hamiltonian
system with Poisson structure and symmetry, Jour. Geom. Phys.,
\textbf{62}(5)(2012), 953-975.


\end{thebibliography}
\end{document}